\documentclass[12pt]{article}
\newif\ifejc
\ejctrue
\ifejc
	\usepackage{e-jc}
\fi
\usepackage[log-declarations=false]{xparse}

\ifejc\else
\fi
\usepackage[colorlinks=true,citecolor=black,linkcolor=black,urlcolor=blue]{hyperref}
\usepackage{float}

\usepackage{tikz}

\usepackage[abbrev,logical-quotes,msc-links]{amsrefs}

\ifejc\else
	\usepackage{amssymb}
	\usepackage{amsthm}

	\newtheorem{theorem}{Theorem}
	\newtheorem{proposition}[theorem]{Proposition}
	\newtheorem{lemma}[theorem]{Lemma}
	\newtheorem{corollary}[theorem]{Corollary}

\fi
\theoremstyle{plain}
\newtheorem{assumptions}[theorem]{Assumptions}

\usepackage{enumitem}
\setlist[enumerate,1]{label=(\roman*)}
		
\newenvironment{Specialindentpar}[1]%
 {\begin{list}{}%
         {\setlength{\leftmargin}{#1}}%
         \item[]%
 }
 {\end{list}}

\ifejc
	\dateline{Dec 14, 2018}{Jan 5, 2021}{Jan 29, 2021}

	\MSC{05C10, 05C83}

	\Copyright{The authors. Released under the CC BY-NC-SA license (International 4.0).}
\fi

\title{All minor-minimal apex obstructions\\
with connectivity two}

\author{Adam S. Jobson \qquad Andr\'{e} E. K\'{e}zdy\\
	\small Department of Mathematics\\[-0.8ex]
	\small University of Louisville\\[-0.8ex]
	\small Louisville, Kentucky, U.S.A.\\
	\small\tt kezdy@louisville.edu}

\begin{document}
\maketitle

\begin{abstract}
A graph is an \emph{apex graph} if it contains a vertex whose deletion leaves a planar graph.
The family of apex graphs is minor-closed and so it is characterized by a finite list of minor-minimal non-members.
The long-standing problem of determining this finite list of apex obstructions remains open.
This paper determines the $133$ minor-minimal, non-apex graphs that have connectivity two.
\end{abstract}

\section{Introduction}
A graph is an {\em apex graph} if it contains a vertex, called an {\em apex}, whose deletion leaves a planar graph.
The family of apex graphs is minor-closed and so Robertson and Seymour's graph minor theorem \cite{RS} implies that this family has a finite list of minor-minimal non-members
(also known as obstructions).
Despite considerable efforts for decades, the problem of determining the list of planar apex obstructions remains open.
In this paper we determine the $133$ minor-minimal non-apex graphs that have connectivity two.

Apex graphs play a key role in what is commonly now referred to as the ``weak structure theorem'' of the Robertson and Seymour's graph minors project
(recently this ``weak structure theorem'' has been optimized by Giannopoulou and Thilikos \cite{MR3072757}).  Apex graphs have also featured prominently in
the resolution of Hadwiger's conjecture for $K_6$-free graphs \cite{MR1238823} and the characterization of linklessly embeddable graphs.  The latter problem highlights 
the Petersen family of graphs \cite{MR1339848}, a significant collection of apex obstructions.
Even with advances in algorithmic refinements of the graph minors project that are focused on determining obstruction sets for apex families \cites{MR2487633,MR2871126},
the very general and theoretical approaches are frustratingly impractical. Disconcertingly the problem of determining all apex-planar obstructions remains open,
despite classical linear-time planarity testing algorithms and
dramatic increases in computing power.

Researchers have made progress characterizing families of graphs that are closely related to the apex graphs.
An old result of Wagner~\cite{MR979099} characterizes the ``almost planar'' graphs; a graph that is non-planar but the deletion of any vertex makes it planar.  The terms ``almost planar'' and ``nearly planar'' appear in many articles with a dizzying variety of meaning.
Gubser \cite{MR1411084} characterizes another family of graphs with another notion of ``almost planar''---a non-planar graph such that
for any edge, either the contraction or the deletion of that edge makes the graph planar (see also \cite{MR3785034}).
Ding and Dzobiak \cites{MR3461970,Dziobiak} determined the $57$ graphs that are the obstructions for the minor-closed family of apex-outerplanar graphs.

It has been a long standing open problem to determine the apex obstructions.  At a conference in $1993$, the second author discussed the apex obstruction problem with Robin Thomas.
Together with Daniel Sanders we decided to share lists of known apex obstructions; the combined list contained $123$ graphs at that time.
Twenty years later we reported $396$ known apex obstructions \cite{Kezdy}.  Our list has since grown to $401$ non-isomorphic apex obstructions.
Other groups have relayed the problem or worked to find obstructions \cites{AdlerOpen,HugoAyala,Eppstein1,Eppstein2,MR3733966,MikePierce}.  Indeed we credit
David Eppstein (see bottom entry of \cite{Eppstein2}) for finding one particularly beautiful $16$-vertex apex obstruction. He describes it this way:
``start with a cube, find a four-vertex independent set, and make three copies of each of its vertices. The resulting $16$-vertex graph has four 
$K_{3,3}$ subgraphs, one for each tripled vertex.''  We refer to this obstruction as the {\em Eppstein graph}.

In this paper we determine the $133$ minor-minimal non-apex graphs that have connectivity two.
It is straightforward to prove there are
three disconnected apex obstructions and no apex obstructions with
connectivity one (see the beginning of \autoref{sec:Simple} for details).  Determining all apex obstructions that have connectivity greater than two is the focus of future research.  Lipton et al. \cite{MR3733966} have shown that apex obstructions have connectivity at most five.
Our proof here presents the connectivity-$2$ apex obstructions in five groups (see \autoref{TreeOutline}); 
this is essentially the same argument we presented
at the AMS Meeting in $2013$ \cite{Kezdy}. \hyperref[sec:Preliminaries]{Section~\ref*{sec:Preliminaries}} introduces basic notation and definitions.  \hyperref[sec:Simple]{Section~\ref*{sec:Simple}} presents
elementary observations about apex obstructions.  Sections \ref{sec:Properties}--\ref{sec:MultipleCuts} present a general outline and resolve four of the five groups of obstructions.

\hyperref[sec:UniqueCut]{Section~\ref*{sec:UniqueCut}} considers the last group of obstructions, the $72$ connectivity-$2$ apex obstructions with
a unique $2$-cut having a planar heavy component (defined in the paragraph before Lemma~\ref{lightComponent}).  These obstructions are considerably more difficult to characterize because they are very close to $3$-connected.
One aspect of connectivity-2 apex obstructions is that they may contain vertices that
are not branch vertices of any Kuratowski subdivision.  There are seven obstructions exhibiting this 
phenomenon, three appear
at the bottom of \autoref{NonPlanarC} and four appear at the bottom of \autoref{Intersecting2Cuts}.
This is not a phenomenon encountered in general $3$-connected graphs that contain a subdivision of $K_{3,3}$ (see statement (6.2) of \cite{MR1339848}). 
The final $72$ connectivity-$2$ apex obstructions behave similarly to the $3$-connected graphs that contain a subdivision of $K_{3,3}$:
all vertices are branch vertices of a Kuratowski subdivision.
Proposition~\ref{ExistenceOfKuratowskis} is the main tool used to prove this claim; it highlights the significance of ``close'' Kuratowski subdivisions (which we do not define in this paper),
a notion remotely reminiscent to `communicating Kuratowski subgraphs' (introduced in \cite{MR1339847}).
The most important consequence of Proposition~\ref{ExistenceOfKuratowskis} is that the apex obstructions having 
a unique $2$-cut with a planar heavy component  
contain three Kuratowski subdivisions whose branch vertices cover the entire vertex set.
This property is key in our proof that the obstruction list is complete.  It appears likely that the study of apex obstructions with higher connectivity will
require similar collections of  Kuratowski subdivisions whose branch vertices 
cover the entire vertex set, though the Eppstein graph shows that
as many as four such Kuratowski subdivisions are needed.  
Unfortunately connectivity appears to be a poor proxy 
for a still missing notion related to close complexes of Kuratowski subdivisions.
After reducing to small graphs the computation of the $72$ connectivity-$2$ apex obstructions having 
a unique $2$-cut with a planar heavy component, the final subsection discusses the computer work applied to show the list
of connectivity-$2$ apex obstructions is complete (see \autoref{graph6Appendix} for graph6 presentations of all $133$ connectivity-$2$ apex obstructions).

Because it may be of independent interest, \hyperref[sec:DoubleApex]{Section~\ref*{sec:DoubleApex}} 
presents another interpretation of the characterization of connectivity-$2$ apex obstructions
in terms of double apex graphs.

\section{Preliminaries}
\label{sec:Preliminaries}

All graphs in this paper are finite, simple, and undirected.
The set of vertices of the graph $G$ is denoted $V(G)$, which is often abbreviated to $V$. 
Similarly $E(G)$, or simply $E$, denotes the set of edges of $G$.
Edges are unordered pairs of vertices, but following standard notation,
the edge $e = \{a,b\}$ is abbreviated $ab$.  The vertices $a$ and $b$ are the {\em endpoints} of
the edge $e=ab$.  The edge $e=ab$ is {\em incident} to $a$ and $b$.  Incidence is often
written as membership, as when $a \in e$ signals that edge $e$ is incident with vertex $a$.
The {\em neighbors} of a vertex $u$ in the graph $G$ is the set $N(u) = \{v \in V : uv \in E\}$,
which is sometimes denoted $N_G(u)$ to emphasize the graph.
For any $v \in V(G)$, the set of edges incident to $v$ is denoted $E_v = \{e \in E(G): v \in e\}$.
Following standard notation, $d_G(v)$, $\delta(G)$, $\kappa(G)$ denote the {\em degree} of the vertex $v$ in $G$,
the {\em minimum degree} of a vertex in $G$, and the vertex {\em connectivity} of $G$, respectively.
A $uv$-path in $G$ is a path whose endpoints are $u$ and $v$.
If $S$ is collection of vertices, then $G[S]$ denotes the subgraph of $G$ induced by $S$.

Often, when the meaning is clear, mathematical elements are recast and operations are overloaded for notational convenience.  So, for example,
the graph induced by $V(G) - S$ is denoted simply $G-S$, which replaces the more syntactically accurate but 
cumbersome $G[V(G) - S]$, thereby overloading the subtraction operator.  
Furthermore, in the case that $v$ is a vertex,
$G - \{v\}$ is abbreviated to $G - v$, which implicitly recasts a vertex as a set. 
Recasting graphs as their vertex sets (and vice versa) appears frequently.  For example, if $A$ and $B$ are
subgraphs of $G$, then $A \subseteq V(G) - B$ tacitly recasts $A$ and $B$ as their vertex sets and 
is short for $V(A) \subseteq V(G) - V(B)$ since it does not make sense to subtract
a graph from a vertex set.  But note carefully that $A \subseteq V(G) - B$ is very different in meaning
from $A \subseteq G - B$; in the latter expression no recasting takes place 
because all variables represent graphs and subtraction makes sense.
Other operators are overloaded naturally as needed.  For example, following standard notation, the addition operator
can be applied to (possibly recast) sets of edges
 or (possibly recast) sets of vertices. For example, $G+e$ represents the graph resulting from the addition of the edge $e$ to the graph $G$.
Similar conveniences appear throughout.

The number of components in $G-S$ is $c(G-S)$.
A vertex set $S$ is a {\em cutset} for a graph $G$
if $c(G)<c(G-S)$.  A cutset with $k$ vertices is called a {\em $k$-cut}.
A cutset $S$ for $G$ is said to {\em separate} vertices $u$ and $v$ if $u$ and $v$
are in the same component of $G$ but different components of $G-S$.  
A path of $G$ {\em crosses} a cutset $S$ if it contains vertices from two different components of $G-S$.
If $S=\{v\}$ is a cutset for $G$, then $v$ is a {\em cut vertex}, also known as a $1$-cut.
If $H$ is a subgraph of a component of $G-S$, then
the {\em $H$-side} (of $G-S$) is the component of $G-S$ that contains $H$.

A subdivision of $K_5$ or $K_{3,3}$ is called a {\em Kuratowski subgraph}.
Kuratowski proved that a graph is non-planar if and only if it contains a Kuratowski subgraph.
A vertex of a Kuratowski subgraph is a {\em branch} vertex if its degree in the Kuratowski subgraph is at least $3$.
A vertex of a Kuratowski subgraph that is not a branch vertex is also called a {\em subdividing vertex}
or a {\em non-branch} vertex.

The contraction of an edge $e$ in the graph $G$ produces a graph denoted $G/e$.
A graph $G$ contains $H$ as a {\em minor}, denoted $H \leq_m G$, if 
a subgraph isomorphic to $H$
can be obtained from a subgraph of $G$ by a sequence of edge contractions.
Observe that the minor order is transitive.  A family ${\cal F}$ is {\em minor-closed}
if $G \in {\cal F}$ and $H \leq_m G$ implies that $H \in {\cal F}$.
If ${\cal F}$ is a minor-closed family, then the minor-minimal graphs that are not in ${\cal F}$ are
called {\em obstructions}; so an obstruction is a graph $G \notin {\cal F}$ such that $H \in {\cal F}$, 
for all $H \leq_m G$ and $H \neq G$.

The disjoint union of graphs $G$ and $H$ is denoted $G+H$. Also $2G$ abbreviates $G+G$.

\section{Simple observations}
\label{sec:Simple}

A graph is an {\em apex graph} if it contains a vertex, called an {\em apex}, whose deletion produces a planar graph.
Every planar graph is an apex graph.
The family of apex graphs is minor-closed and so it has finite list of minor-minimal non-members, also known as forbidden minors or obstructions.
Let ${\cal F}$ denote the finite set of obstructions for apex graphs.

We begin with a few simple observations about graphs in ${\cal F}$.

\begin{lemma} \label{MinDegree3} If $G \in {\cal F}$, then $\delta(G) \geq 3$.
\end{lemma}
\begin{proof} Suppose, to the contrary, that there exists $G \in {\cal F}$ and $v \in V(G)$ with $d_G(v) \leq 2$.  If $d_G(v) \leq 1$, then
$G-v$ is apex with an apex, $w$.  Note that $G-v-w$ is planar which implies that $G-w$ is also planar, contradicting $G \in {\cal F}$.
If $d_G(v) = 2$, then consider an edge $e$ incident to $v$.  The contraction of $e$ produces an apex graph $G/e$ with an apex vertex, $w$.
Now $G/e - w$ is planar which implies that $G-w$ is planar, contradicting $G \in {\cal F}$.
\end{proof}

\begin{lemma} \label{Disconnected} The disconnected graphs in ${\cal F}$ are $2K_5, 2K_{3,3},$ and $K_5 + K_{3,3}$.
\end{lemma}
\begin{proof} The reader can easily verify that $2K_5, 2K_{3,3},$ and $K_5 + K_{3,3}$ are disconnected
graphs in ${\cal F}$. It suffices to show there are no others.  Consider a disconnected graph $G \in {\cal F}$.  Each component of $G$ must be
an apex graph since removing any one edge from $G$ produces an apex graph.  
There must be at least two non-planar components since otherwise the whole graph is an apex graph.
It follows that there are exactly two components and each is non-planar apex.
Consider an arbitrary component of $G$, call it $C$.  Let $H$ be a Kuratowski subgraph in $C$.
If there is an edge $e \in E(C)- E(H)$, then $G-e$ is not apex, contradicting that $G \in {\cal F}$;
hence, $E(H)=E(C)$.  Lemma~\ref{MinDegree3} guarantees $\delta(G) \geq 3$, from which it follows that all the vertices of
$C$ are branch vertices of $H$; consequently, $C \cong K_5$ or $C \cong K_{3,3}$.
\end{proof}

\begin{lemma} \label{OneConnected} Every connected graph in ${\cal F}$ has connectivity at least two.
\end{lemma}
\begin{proof} Suppose, to the contrary, that $G \in {\cal F}$ and $\kappa(G) = 1$.  Let $v$ be a cut vertex of $G$.
Now $G-v$ is non-planar; therefore, there exists a Kuratowski subgraph, $H$, in one of the components, say component $C$, of $G-v$.
Define $T = \{uv \in E(G) : u \in V(C) \}$.  Because $G-T$ is apex with apex vertex in $V(H) \subseteq V(C)$, it follows that
$G-C$ is planar.  Let $e$ be an edge incident to $v$ that is also incident to a vertex in a component of $G-v$ other than $C$.  Because $G-e$ is apex,
there is an apex of $G-e$, say $w$, that is in $V(H) \subseteq V(C)$.  Because $G[V(C)\cup \{v\}]$ is apex (with apex $w$) and $G-C$ is planar, it follows
that $G-w$ is planar: embed $G[V(C)\cup \{v\}] - w$ in the plane with $v$ on the exterior face, separately embed $G-C$ in the plane with $v$ on the exterior face, and identify $v$'s in these embeddings, producing a planar embedding
of $G-w$.  This contradicts that $G$ is not an apex graph.
\end{proof}

Because of Lemmas \ref{Disconnected} and \ref{OneConnected}, we now consider only obstructions in ${\cal F}$ that 
have connectivity at least two.  Indeed, in this paper we determine all graphs $G \in {\cal F}$ such
that $\kappa(G) =2$; there are precisely $133$ of them.

For $G \in {\cal F}$ and $v \in V(G)$, let $H_v$ denote a Kuratowski subgraph in $G-v$.  Note that $H_v$ witnesses that $G$ is not apex; it 
is a {\em Kuratowski witness} for $v$.  

If $S \subseteq V$ is a cutset for $G$ and $C$ is a component of $G-S$, then the \emph{augmentation of $C$} is
\[C^+ = G[V(C) \cup S] + \{uv : u,v \in S\}.\]
This is also referred to as the {\em augmented component} of $G-S$ obtained from $C$.

The next lemma gathers several elementary properties of Kuratowski witnesses.

\begin{lemma} \label{Elementary}
Suppose $G \in {\cal F}$ and $\kappa(G) \ge 2$.  For all $u,v,w \in V(G)$, and all Kuratowski witnesses
$H_u$, $H_v$ and $H_w$ for $u$, $v$ and $w$, respectively,
\begin{itemize}
\item[i)] $ V(H_u \cap H_v) \neq \varnothing$,
\item[ii)] if $V(H_u \cap H_v \cap H_w) = \varnothing$, then $E(G) = E(H_u \cup H_v \cup H_w)$,
\item[iii)] if $V(H_u \cap H_v \cap H_w) = \varnothing$ and $x \in V(H_u) - V(H_v \cup H_w)$, then
$x$ is a branch vertex of $H_u$,
\item[iv)] if $\kappa(G) = 2$, $S$ is a $2$-cut of $G$, and $s \in S$, then for any
Kuratowski witness $H_s$, there
exists an augmented component $C^+$ of $G-S$ such that $V(H_s) \subseteq V(C^+)$.
\end{itemize}
\end{lemma}
\begin{proof} {\it i)} If $V(H_u \cap H_v) = \varnothing$, then $G = H_u + H_v$ is a disconnected obstruction in ${\cal F}$.
{\it ii)} Assume $V(H_u \cap H_v \cap H_w) = \varnothing$.  If there were an edge $e \in E(G) - E(H_u \cup H_v \cup H_w)$, then $G-e$ would have no apex.
{\it iii)} Suppose $V(H_u \cap H_v \cap H_w) = \varnothing$ and $x \in V(H_u) - V(H_v \cup H_w)$.  Assume, to the contrary, that $x$ is not a branch vertex of $H_u$; so $d_H(x)=2$.  Lemma~\ref{MinDegree3} implies $d_G(x)\geq 3$. Consequently there is at least one edge $e$ incident to $x$ that does not belong to $E(H_u \cup H_v \cup H_w)$, contradicting that $G-e$ has an apex.
{\it iv)} Assume $G \in {\cal F}$, $\kappa(G) = 2$, $S$ is a $2$-cut of $G$, and $s \in S$.  Consider a Kuratowski witness $H_s$.  Observe that some augmented component of $G-S$ must contain all the vertices in $V(H_s)$ because no path of the $2$-connected subgraph $H_s$ can cross the $1$-cut $S- s$.
\end{proof}

\begin{lemma} \label{RemoveEdge} Suppose that $G \in {\cal F}$, $\kappa(G) = 2$, and $S=\{a,b\}$ is a $2$-cut of $G$. 
Every augmented component $C^+$ of $G-S$ contains an edge $e \in E(C^+) - ab$ such that
there exists an $ab$-path in $C^+ - e - ab$.
\end{lemma}
\begin{proof}  Consider an arbitrary component of $G-S$, call it $C$, and its augmentation $C^+$.
Let $P$ be an $ab$-path in $C^+ - ab$ ($P$ exists because $S$ is a minimum cutset of $G$).  
The path $P$ contains a vertex $v \not\in S$.  Let $e$ be an edge incident to
$v$ that is not in $P$; the existence of such an edge follows from Lemma~\ref{MinDegree3}.
Now $e$ is the desired edge in $E(C^+) - ab$.
\end{proof}

Next we present a technical lemma that will be applied in several upcoming arguments.

\begin{lemma} \label{BranchVerticesOfHw} Suppose that $G \in {\cal F}$, $\kappa(G) = 2$, $S=\{a,b\}$ is a $2$-cut of $G$,
and $C$ is a component of $G-S$.
If $e \in E(C^+) - ab$ and
there exists an $ab$-path in $C^+ - e - ab$, then for any apex $w$ of $G - e$ and
any Kuratowski subgraph $H_w$ in $G$ avoiding $w$,
the branch vertices of $H_w$ are all in $V(C^+)$.
\end{lemma}
\begin{proof}  
Because $G-e-w$ is planar, it follows that $e \in E(H_w)$.  All of the branch vertices
of $H_w$ must be in the same augmented component of $G-S$ because $|S|=2$ but $H_w$ is a subdivision of $K_5$ or $K_{3,3}$ which are $3$-connected.  If the branch vertices of $H_w$ are
not in $V(C^+)$, then $e$ appears in $H_w$ only as an edge along a path connecting $a$ and $b$.
Consequently, a Kuratowski subgraph in $G-w-e$ would exist by replacing this $ab$-path in $C^+$ by
another $ab$-path from $C^+ - e - ab$ (which exists by assumption), contradicting that $w$ is
an apex vertex for $G-e$.
So the branch vertices of $H_w$ must be in $V(C^+)$.
\end{proof}

Recall that, if $S \subseteq V$ is a cutset for $G$, then $c(G-S)$ is equal to the number of components in $G-S$.

\begin{lemma} \label{TwoComponents} If $G \in {\cal F}$, $\kappa(G) = 2$, and $S$ is a $2$-cut of $G$, 
then $c(G-S) = 2$.
\end{lemma}
\begin{proof} Assume, to the contrary, that $c(G-S)>2$.  Let $S=\{a,b\}$ and let $C_1,C_2,C_3$ be three
components of $G-S$.  Applying Lemma~\ref{Elementary} part {\it iv)}, we may assume that $V(H_a) \subseteq C_1^+$.
It follows from Lemma~\ref{Elementary} parts {\it i)} and {\it iv)} that $V(H_b) \subseteq C_1^+$.  
Lemma~\ref{RemoveEdge} guarantees there exists 
an edge $e$ in $E(C_3^+) - ab$ such that an $ab$-path remains in $C_3^+ -ab-e$.  Because $G \in {\cal F}$, the graph $G-e$
is an apex graph; let $w$ be an apex for $G-e$.  
Consider $H_w$, a Kuratowski subgraph avoiding $w$.  By Lemma~\ref{BranchVerticesOfHw} the branch vertices of $H_w$ must be in $V(C_3^+)$.  

Now Lemma~\ref{RemoveEdge} guarantees there exists 
an edge $f \in E(C_2^+) - ab$ such that an $ab$-path remains in $C_2^+ -ab-f$.
Notice that if $f$ is an edge in $H_w$ then, because the branch vertices of $H_w$ must be in $V(C_3)$,
the edge $f$ appears in $H_w$ only as an edge along an $ab$-path in $C_2^+ -ab$.
Consequently there is a Kuratowski subgraph $H_w^*$ in $G[V(C_2 \cup C_3) \cup \{a,b\}]$
that avoids the edge $f$; it is obtained from $H_w$ by replacing, if necessary, any
$ab$-path in $H_w \cap (C_2^+ -ab)$ by an $ab$-path in $C_2^+ - ab  - f$.
This implies that $H_a$, $H_b$ and $H_w^*$ share
no common vertex, contradicting $G-f$ is apex.
\end{proof}

\begin{lemma} \label{AugmentComponent} If $G \in {\cal F}$, $\kappa(G) = 2$, and $S$ is any $2$-cut of $G$, 
then the augmentation of any component of $G-S$ is non-planar.
\end{lemma}
\begin{proof} Let $S=\{a,b\}$ be an arbitrary $2$-cut of $G$.  Lemma~\ref{TwoComponents} guarantees
that $G-S$ has exactly two components, call them $C_1$ and $C_2$.
Applying Lemma~\ref{Elementary} part {\it iv)}, we may assume that $H_a \subseteq C_1^+$.
It follows from Lemma~\ref{Elementary} parts {\it i)} and {\it iv)} that $H_b \subseteq C_1^+$.
It suffices to prove that $C_2^+$ is non-planar.  
Lemma~\ref{RemoveEdge} guarantees there exists 
an edge $e \in E(C_2^+) -ab$ such that an $ab$-path remains in $C_2^+ -ab-e$.  Because $G \in {\cal F}$, the graph $G-e$
is an apex graph; let $w$ be an apex for $G-e$ and let $H_w$ be a Kuratowski subgraph of $G$ avoiding $w$.  
By Lemma~\ref{BranchVerticesOfHw} the branch vertices of $H_w$ must be in $V(C_2^+)$.  
This implies that $H_w \subseteq C_2^+$, so $C_2^+$ is non-planar.
\end{proof}

Consider $G \in {\cal F}$ with $\kappa(G) = 2$ and $S=\{a,b\}$ any $2$-cut of $G$.
By Lemma~\ref{TwoComponents}, there are only two components
of $G-S$, call them $C_1$ and $C_2$.
We may assume that $H_a \subseteq C_1^+$ and $H_b \subseteq C_1^+$.
We call $C_1$ the {\em heavy} component of $G-S$ because $H_a,H_b \subseteq C_1^+$;
$C_2$ is the {\em light} component.  
A $2$-cut of $G$ is {\em basic} if the vertex set of its heavy component is minimal (with respect to set inclusion); that is, no $2$-cut produces
a heavy component whose vertex set is properly contained in this one.

\begin{lemma} \label{lightComponent} If $G \in {\cal F}$, $\kappa(G) = 2$, and $S$ is any basic $2$-cut of $G$, 
then the augmentation of the light component of $G-S$ is isomorphic to $K_5$, $K_{3,3}$, or $K_{3,3}+e$.
\end{lemma}
\begin{proof} Let $S=\{a,b\}$ be an arbitrary basic $2$-cut of $G$.  
Let $C_1$ and $C_2$ be the two components of $G-S$,
with $C_1$ the heavy component and $C_2$ the light one.  
Lemma~\ref{RemoveEdge} guarantees there exists 
an edge $e \in E(C_2^+) -ab$ such that an $ab$-path remains in $C_2^+ -ab-e$.  Because $G \in {\cal F}$, the graph $G-e$
is an apex graph.  Let $w$ be an apex for $G-e$ and let $H_w$ be a Kuratowski subgraph of $G$ avoiding $w$.  
Observe $w \in V(H_a \cap H_b) \subseteq V(C_1)$.
By Lemma~\ref{BranchVerticesOfHw} the branch vertices of $H_w$ must be in $V(C_2^+)$.
Furthermore, there must be some part of $H_w$ that is not in $C_2$ since
$V(H_w \cap H_a) \neq \varnothing$ and $V(H_w \cap H_b) \neq \varnothing$, so in particular,
$\{a,b\} \subseteq V(H_w)$.

\vspace{\baselineskip}
\noindent\textsc{Claim}: $w$ is an apex for $G - f$ and $G/f$, for all $f \in E(C_2^+)-\{ab\}$.

\begin{Specialindentpar}{\parindent}
Suppose, to the contrary, there is an edge $f \in E(C_2^+)-\{ab\}$ such that
$w$ is not an apex for $G-f$ or $G/f$.
There must be
an apex vertex for $G-f$ and $G/f$.  Let $z$ be an apex for
$G-f$ if $w$ is not an apex for it; otherwise let $z$ be an apex for $G/f$.  
Clearly $z \in V(H_a \cap H_b \cap H_w) \subseteq V(C_1)$ and
$z \notin \{a,b,w\}$.
Because the branch vertices of $H_w$ are in $C_2^+$ and $z \notin \{a,b\}$, it follows that 
$H_w \cap G[V(C_1) \cup \{a,b\}]$ is
an $ab$-path through $z$ and $ab \notin E(G)$.  Indeed $z$ must be  
a cut vertex in $G[V(C_1) \cup \{a,b\}]$ that separates $a$ from $b$ since otherwise an $ab$-path
in $G[V(C_1) \cup \{a,b\}]$ that avoids $z$ could substitute in $H_w$ to create a Kuratowski subgraph
in $G-f-z$ or $G/f-z$.

Now consider the two components of $G[V(C_1) \cup \{a,b\}] - z$; call them $U_a$ and $U_b$, where
$a \in V(U_a)$ and $b \in V(U_b)$.
Set $U_a^+ = G[V(U_a) \cup \{z\}]$ and $U_b^+ = G[V(U_b) \cup \{z\}]$.

Note that $H_a$ must be a subgraph of $U_a^+$ or $U_b^+$, since $H_a$ is a subgraph of $G[V(C_1) \cup \{a,b\}]$
and $z$ is a cutvertex for this graph.  If $H_a$ is a subgraph of $U_a^+$, then $H_z$ is also 
a subgraph of $U_a^+$ since Lemma~\ref{Elementary} part (i) shows that $V(H_a \cap H_z) \neq \varnothing$.
This implies that $\{z,a\}$ is a $2$-cut of $G$ with a heavy component properly contained in $C_1$, contradicting
that $S$ is basic.
So $H_a$ is a subgraph of $U_b^+$.
Similarly $H_b$ is a subgraph of $U_a^+$.
Lemma~\ref{Elementary} part (i) states that $V(H_a \cap H_b) \neq \varnothing$, so
$V(U_a \cap U_b) = \{z\}$ implies $V(H_a \cap H_b) = \{z\}$.
But recall that $w \in V(H_a \cap H_b)$;
so $w=z$, a contradiction.

\end{Specialindentpar}

Because $w$ is an apex for $G-f$, for all $f \in E(C_2^+)-\{ab\}$, it follows that
$E(C_2^+)-\{ab\} \subseteq E(H_w)$.  Therefore there are no vertices in $V(C_2)$ 
that have degree two in $H_w$
since $\delta(G)\geq 3$ would otherwise guarantee an edge in $E(C_2^+)-\{ab\} - E(H_w)$.  Consequently,
all the vertices of $C_2$ are branch vertices of $H_w$.
If $a$ is not a branch vertex of $H_w$, then 
it is a subdividing vertex of $H_w$ so there is an edge of $H_w$ incident to $a$ whose other endpoint is in $C_2$ 
(because
the branch vertices of $H_w$ are all in $C_2 \cup \{a,b\}$); call this edge $f$.  
Now $G/f$ contains $H_w/f$, contradicting the claim that established $w$ is an apex for $G/f$.
So $a$ is a branch vertex of $H_w$.
Symmetrically $b$ is a branch vertex of $H_w$.
Because all of the vertices in $C_2 \cup \{a,b\}$ are branch vertices of $H_w$ and 
the branch vertices of $H_w$ are all in $C_2 \cup \{a,b\}$, it now follows that
$C_2^+$ is isomorphic to $K_5$, $K_{3,3}$, or $K_{3,3}+e$.
\end{proof}

The next lemma, a consequence of Lemma~\ref{AugmentComponent} and Lemma~\ref{lightComponent},  is a powerful tool to analyze obstructions with more than one $2$-cut.
\begin{lemma} \label{All2CutsAreBasic}
Suppose that $G \in {\cal F}$ and $\kappa(G) = 2$. Every $2$-cut of $G$ is basic.
\end{lemma}
\begin{proof}
Suppose that $S = \{a,b\}$ is an arbitrary $2$-cut of $G \in {\cal F}$.    Assume, to the contrary, that $S$ is not basic.  Let $C_S$ (resp. $L_S$) denote the heavy (resp. light)
component of $G-S$.  Because $S$ is not basic,
there exists a basic $2$-cut $Q$ with heavy (resp. light) component $C_Q$ (resp. $L_Q$) 
such that $V(C_Q) \subsetneq V(C_S)$, or
equivalently, $\overline{V(C_S)} \subsetneq \overline{V(C_Q)}$.
Now
$$
V(L_S^+) = V(L_S) \cup S = \overline{V(C_S)} \subsetneq \overline{V(C_Q)} = V(L_Q) \cup Q = V(L_Q^+),
$$
so there exists a vertex $z \in V(L_Q^+) - V(L_S^+).$
Because $Q$ is basic, Lemma~\ref{lightComponent} implies that 
$L_Q^+$ is isomorphic to $K_5$, $K_{3,3}$, or $K_{3,3}+e$.  But, $V(L_S^+) \subsetneq V(L_Q^+)$ and 
$z \in V(L_Q^+) - V(L_S^+)$
means that $L_S^+$ is a subgraph of $K_5 - z$, $K_{3,3} - z$, or $K_{3,3}+e - z$.
Each of these possibilities contradicts Lemma~\ref{AugmentComponent} which shows that $L_S^+$ is non-planar.
\end{proof}

\section{The five types of obstructions and some shared properties}
\label{sec:Properties}

\begin{figure}[ht]
	\centering
	\includegraphics[page=1]{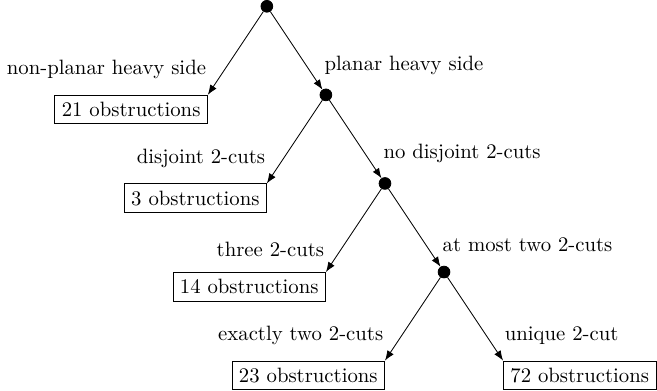}%
	\caption{Partition of the 133 connectivity-$2$ apex obstructions into five types.}
	\label{TreeOutline}
\end{figure}

The connectivity-$2$ obstructions to the apex family are arranged into five groups.
\autoref{TreeOutline} shows a partition of these $133$ obstructions according to whether the heavy component 
of a $2$-cut induces a planar graph and further properties of $2$-cuts.  This partition follows the outlines of our characterization of these graphs.

\begin{theorem} \label{edgecover} Suppose that $G \in {\cal F}$, $\kappa(G) = 2$, $S=\{a,b\}$ is any $2$-cut of $G$, and $C$ is the heavy component of $G-S$. 
If $H_a$ and $H_b$ are Kuratowski subgraphs of $G$ avoiding $a$ and $b$, respectively, then
$$E(C) \subseteq E(H_a\cup H_b).$$
\end{theorem}
\begin{proof} 
We argue by contradiction.  Assume that $f \in E(C) - E(H_a \cup H_b)$.
By Lemma~\ref{All2CutsAreBasic}, $S$ is a basic $2$-cut.
Lemma~\ref{lightComponent} guarantees that the augmentation of the light 
component, $L$, of $G-S$ is isomorphic to $K_5$, $K_{3,3}$, or $K_{3,3}+e$.
Now $G-f$ is apex with apex $z$, say.  Observe that $z \in V(H_a \cap H_b) \subseteq V(G) - (V(L) \cup \{a,b\})$.
There can therefore be no $ab$-path in $G[C \cup \{a,b\}] - f -z$ since otherwise the contraction of this
path would, together with $L$, realize a Kuratowski subgraph in $G-f - z$.
If the endpoints of $f$ are in same component (we may assume the component containing $a$) of $G[C \cup \{a,b\}]-z$ (see \autoref{Case1}), 
then, because $H_z$ must contain $f$ and $H_a \cap H_z \neq \varnothing$, both $H_a- z$ and $H_z$ are also in this component.
Consequently the $2$-cut $\{a,z\}$ has a heavy component properly
contained in $C_1$, contradicting that $S$ is basic.

\begin{figure}[H]
	\centering
	\includegraphics[page=2]{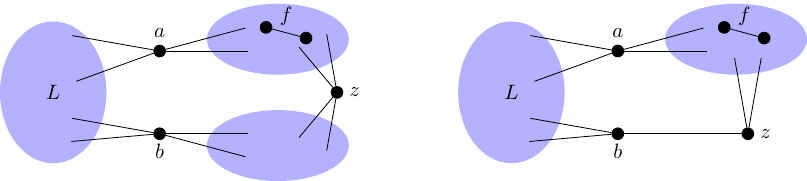}%
	\caption{Cases in which endpoints of $f$ are in same component of $G-L-S-z$.}
	\label{Case1}
\end{figure}

Observe that $f$ is not incident to $z$.
So we may assume that the endpoints of $f$, $u$ and $v$, are in opposite components of $G[C \cup \{a,b\}]-z$.
Let us label these components $A$ and $B$, where $A$ (resp. $B$) denotes the component containing neighbors
of $a$ (resp. $b$).  Without loss of generality, $u \in A$ and $v \in B$ (see \autoref{Case2}).

\begin{figure}[ht]
	\centering
	\includegraphics[page=3]{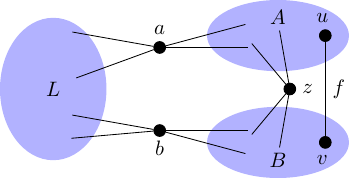}%
	\caption{Endpoints of $f$ are in different components of $G-L-S-z$.}
	\label{Case2}
\end{figure}

Because $f \not\in E(H_a\cup H_b)$, it follows that $|\{u,v\} \cap V(H_a)| \leq 1$
and $|\{u,v\} \cap V(H_b)| \leq 1$; consequently $H_a$ and $H_b$ exist in $G/f$, the graph obtained
from $G$ by contracting $f$.

\vspace{\baselineskip}
\noindent\textsc{Claim}: There are two internally vertex-disjoint paths from $a$ to $\{u,z\}$ in $G-b$.

\begin{Specialindentpar}{\parindent}
Suppose, to the contrary, that there are not 
two internally vertex-disjoint paths from $a$ to $\{u,z\}$ in $G-b$.  Menger's theorem then guarantees
that there is a vertex $w$ in the subgraph of $G$ induced by $A \cup \{a,z\}$ separating $a$ from $\{u,z\}$
(see \autoref{A2Fan}).

\begin{figure}[ht]
	\centering
	\includegraphics[page=4]{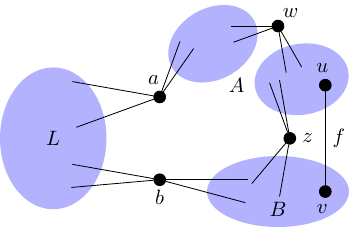}%
	\caption{No two internally vertex-disjoint paths implies a vertex $w$.}
	\label{A2Fan}
\end{figure}

Recall that $f \in H_z$ since $z$ is an apex for $G-f$.  Therefore $H_z$ in on the $f$-side of the $2$-cut
$\{b,w\}$.  Because $H_b \cap H_z \neq \varnothing \neq H_w \cap H_z$, it follows that $H_b$ and $H_w$ are
also on the $f$-side of the $2$-cut $\{b,w\}$.  Consequently the $2$-cut
$\{b,w\}$ has a heavy component properly contained in $C_1$, contradicting that $S$ is basic.
\end{Specialindentpar}

By symmetry, there are two internally vertex-disjoint paths from $b$ to $\{v,z\}$ in $G-a$.
This means that there are two internally vertex-disjoint $ab$-paths in $C \cup \{a,b\}$ that
remain internally vertex-disjoint after contracting the edge $f$.
Now consider $G/f$; it is an apex graph with an apex vertex $z^\prime$.
As noted earlier, $H_a$ and $H_b$ exist in $G/f$ so $z^\prime \in V(H_a \cap H_b) \subseteq V(C)$.
But the deletion of the vertex $z^\prime$ leaves an $ab$-path on the $C$-side of the $2$-cut
$\{a,b\}$ in $G/f$.  In particular, this $ab$-path can be contracted to the edge $ab$ which, 
together with $L$ produces a Kuratowski subgraph in $G/f - z^\prime$, contradicting that
$z^\prime$ is an apex vertex.
\end{proof}

\section{The 21 connectivity-2 obstructions with a non-planar heavy component}
\label{sec:HeavyNonplanar}

We now turn to characterizing apex obstructions in which some $2$-cut
has a non-planar heavy component.

\begin{figure}[ht]
	\centering
	\includegraphics[page=5,scale=0.95]{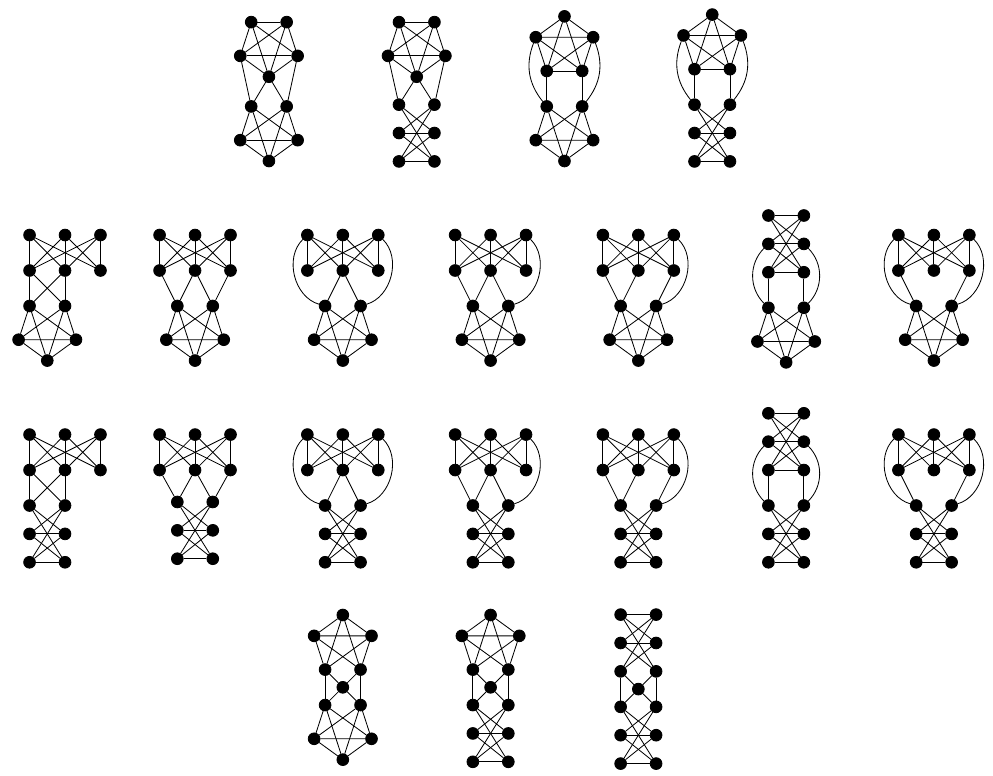}%
	\caption{The $21$ connectivity-$2$ apex obstructions with a $2$-cut having a non-planar heavy component.}
	\label{NonPlanarC}
\end{figure}

\begin{theorem}
Suppose that $G \in {\cal F}$ and $\kappa(G) = 2$.
If some $2$-cut of $G$ has a non-planar heavy component, 
then $G$ is isomorphic to a graph in \autoref{NonPlanarC}.
\end{theorem}
\begin{proof}  
Let $S=\{a,b\}$ be a $2$-cut of $G$ with non-planar heavy component 
$C$ of $G-S$.
By Lemma~\ref{All2CutsAreBasic}, $S$ is a basic $2$-cut.
Lemma~\ref{lightComponent} guarantees that the augmentation of the light 
component, $L$, of $G-S$ is isomorphic to $K_5$, $K_{3,3}$, or $K_{3,3}+e$.
Observe that $ab \not\in E(G)$ since otherwise $G$ contains two disjoint Kuratowski subgraphs,
one in $G[V(L) \cup S]$ and one in $C$,
contradicting that $G \in {\cal F}$ and $\kappa(G) = 2$.
Similar reasoning shows also that $L^+ \not\cong K_{3,3}+e$,
a notable characteristic of this collection of obstructions.

Let $K$ be a minimum order Kuratowski subgraph of $C$. 
Let $H_a = K = H_b$ be the two Kuratowski subgraphs avoiding $a$ and $b$.
Theorem~\ref{edgecover} implies that $E(G) \subseteq E(L) \cup E_a \cup E_b \cup E(K)$,
where recall that $E_v = \{e \in E(G): v \in e\}$ denotes the edges of $G$ incident to $v$.
This means that $C = K$ and, because $\delta(G) \geq 3$, the only possible non-branch vertices of $K$ are neighbors
of either $a$ or $b$.

\vspace{\baselineskip}
\noindent\textsc{Step 1}: For $v \in S$, $|N(v) \cap V(K)| \geq 2$.
\begin{Specialindentpar}{\parindent}
   We may assume $v = a$.
   Suppose, to the contrary, that $|N(a) \cap V(K)| \leq 1$.  If $|N(a) \cap V(K)| = 0$, then
   $\kappa(G) < 2$ so we may assume that $N(a) \cap V(K) = \{w\}$, for some vertex $w \in V(K)$.
   Now consider the $2$-cut $\{b,w\}$ of $G$.  Because $H_w$ can not be disjoint from $K = H_b$, 
	 it follows that $H_w$ is a subgraph of $G[V(K) - w + b]$.  Consequently
	 the $2$-cut $\{b,w\}$ has a heavy component that is properly contained in $C_1$,
	 contradicting that $S$ is basic. 
\end{Specialindentpar}

\vspace{0.5\baselineskip}
\noindent\textsc{Step 2}: For $v \in S$, $|N(v) \cap V(K)| \leq 2$.
\begin{Specialindentpar}{\parindent}
   Without loss of generality $v = a$.
   Suppose, to the contrary, that $|N(a) \cap V(K)| \geq 3$.  Choose $w \in N(a) \cap V(K)$
   and set $e = aw$.  Consider $G-e$; it is apex with an apex vertex $z$, say.
   Now $z \in V(K)$ since $K \subseteq G-e$.  
	 Because $|N(a) \cap V(K)| \geq 3$, the vertex $a$ has at least one neighbor in $V(K) -\{w,z\}$.
	 By Step 1, the vertex $b$ has
   at least one neighbor in $G-e-z$. Since $K-z$ is connected, there is an $ab$-path in $G[V(K) \cup \{a,b\}] - e - z$ 
	 that produces a Kuratowski subgraph corresponding to the augmentation
   of $L$ in $G-e-z$, a contradiction. 
\end{Specialindentpar}

\vspace{0.5\baselineskip}
\noindent\textsc{Step 3}: Any subdividing vertex of $K$ is adjacent to both $a$ and $b$.
\begin{Specialindentpar}{\parindent}
   Suppose that $w$ is a subdividing vertex of $K$.  As noted prior to Step 1, $w \in N(a) \cup N(b)$ since otherwise
   the minimum degree three of $w$ implies an edge incident to $w$ that is not covered by $K$, contradicting
   Theorem~\ref{edgecover}.  Assume now, contrary to the claim, that 
   $|N(w) \cap \{a,b\}| = 1$.  Without loss of generality, $aw \in E(G)$.
   Choose $w^\prime \in (N(w) \cap K) - N(a)$; such a vertex exists
   since $w$ has two neighbors in $K$ but $a$ only has one more neighbor in $K$ besides $w$.  
   Consider $G/ww^\prime$; it is apex with apex $z$, say.  Note that contracting $ww^\prime$
   preserves a version of $K$ and both $a$ and $b$ still have two neighbors
   in this version of $K$ (because $a$ is not adjacent to $w^\prime$ and $b$ is not adjacent to $w$) implying that there is an
	 $ab$-path in $G[V(K) \cup \{a,b\}]/ww^\prime-z$ which determines a
	 Kuratowski subgraph corresponding to the augmentation
   of $L$ in $G/ww^\prime-z$, a contradiction.
\end{Specialindentpar}

\vspace{0.5\baselineskip}
\noindent\textsc{Step 4}: $K$ has at most one subdividing vertex.
\begin{Specialindentpar}{\parindent}
   Observe that, by Steps 1--3, there are at most two subdividing vertices of $K$.
   Suppose, to the contrary, that $x$ and $y$ are subdividing vertices of $K$.  
   Prior steps guarantee $N(a) \cap V(K) = \{x,y\} = N(b) \cap V(K)$.  Let $w \in (V(K) \cap N(x)) - y$.
	 Consider these three
	 Kuratowski subgraphs of $G/xw$: the subgraph induced by $V(L) \cup \{a,b,x\}$, the subgraph induced by $V(L) \cup \{a,b,y\}$, 
	 and $K/xw$.
	 These three Kuratowski subgraphs have no common vertex, contradicting that $G/xw$ is apex.
\end{Specialindentpar}

\vspace{0.5\baselineskip}
\noindent\textsc{Step 5}: If $K$ has a subdividing vertex $w$, then $N(w) \cap V(K) \subseteq N(a) \cup N(b)$.
\begin{Specialindentpar}{\parindent}
   Let $x$ and $y$ be the neighbors of $w$ in $K$.  If $x \not\in N(a) \cup N(b)$, then 
   contracting $wx$
   preserves a version of $K$ that must contain an apex $z$ for $G/wx$.  However $a$ and $b$ still have two neighbors
   to this version of $K$ implying that there is a Kuratowski subgraph corresponding to the augmentation
   of $L$ in $G/wx-z$, a contradiction.  Therefore, $x \in N(a) \cup N(b)$.  By symmetry, $y \in N(a) \cup N(b)$.
\end{Specialindentpar}

So, in summary, $L^+$ is isomorphic to $K_5$ or $K_{3,3}$, $C$ is a Kuratowski subgraph, $K$, 
with at most one subdividing vertex, and $|N(a) \cap V(K)| = 2 = |N(b) \cap V(K)|$.  
If $K$ has a subdividing vertex $w$ with $N(w) \cap V(K) = \{x,y\}$, then $xy$ is not an edge of $G$ 
(by minimality of $K$) and,
without loss of generality, $N(a) \cap V(K) = \{w,x\}$ and $N(b) \cap V(K) = \{w,y\}$.
The reader can now easily verify that the graphs in \autoref{NonPlanarC} enumerate the possible apex obstructions
with these properties.
\end{proof}

\section{The 3 connectivity-2 apex obstructions with planar heavy components and two disjoint 2-cuts}
\label{sec:Disjoint}

\begin{figure}[ht]
	\centering
	\includegraphics[page=6]{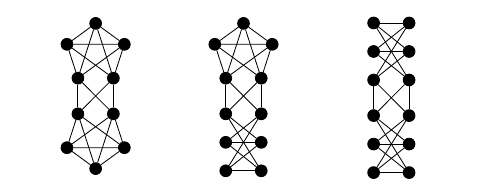}%
	\caption{The $3$ connectivity-$2$ apex obstructions in which there are disjoint $2$-cuts and every $2$-cut has a planar heavy component.}
	\label{Disjoint2Cuts}
\end{figure}

To complete the characterization of connectivity-$2$ apex obstructions, it suffices to
consider connectivity-$2$ apex obstructions in which every $2$-cut has a planar heavy component.
We next consider such connectivity-$2$ apex obstructions with two disjoint $2$-cuts.
There are three of these (see \autoref{Disjoint2Cuts}).

\begin{theorem}
Suppose $G \in {\cal F}$, $\kappa(G) = 2$, and every $2$-cut has a planar heavy component.
If $G$ has disjoint $2$-cuts, then $G$ is isomorphic to a graph in \autoref{Disjoint2Cuts}.
\end{theorem}
\begin{proof}  
Suppose that $S$ and $T$ are disjoint $2$-cuts of $G$.
Let $S=\{a,b\}$, and let $C_S$ be the heavy component of $G-S$.
Lemma~\ref{All2CutsAreBasic} yields that $S$ and $T$ are basic $2$-cuts.
Lemma~\ref{lightComponent} guarantees that the augmentation, $L_S^+$, of the light 
component, $L_S$, of $G-S$ is isomorphic to $K_5$, $K_{3,3}$, or $K_{3,3}+e$, each of which is $3$-connected.
Similarly $L_T^+$, the augmentation of the light component of $T$, is isomorphic to $K_5$, $K_{3,3}$, or $K_{3,3}+e$.

\vspace{\baselineskip}
\noindent\textsc{Case 1}: $|T \cap V(L_S)| = 2$.
\begin{Specialindentpar}{\parindent}
   
  The deletion of $T$ from $L_S^+$ is still connected.  It follows that $G-T$ is connected, a contradiction.
\end{Specialindentpar}

\vspace{0.5\baselineskip}
\noindent\textsc{Case 2}: $|T \cap V(L_S)| = 1$.
\begin{Specialindentpar}{\parindent}
   The deletion of a vertex of $T$ from $L_S^+$ is still connected.  It follows that $G-T$ can only be disconnected
   if a vertex of $G$ is a cut vertex of $G$, contradicting that $\kappa(G)=2$.
\end{Specialindentpar}

\vspace{0.5\baselineskip}
\noindent\textsc{Case 3}: $|T \cap V(L_S)| = 0$.
\begin{Specialindentpar}{\parindent}
   Let $T=\{u,v\} \subset C_S$ with $L_T$ (resp. $C_T$) the light (resp. heavy) component of $G-T$.  
	 Because $G-T$ has two components and one contains $L_S \cup S$, it follows that
	 $L_S  \cup S \subset C_T$ because $T$ is basic.  Symmetrically, $L_T \cup T \subset C_S$.
	 
	 We claim that there exists an $av$-path in $G-\{b,u\}-V(L_S)-V(L_T)$.
	 It suffices to find an $av$-path in $G-\{b,u\}$.  If $G-\{b,u\}$ is connected, then
	 there is nothing to prove; so assume that $R=\{b,u\}$ is a $2$-cut of $G$.
	 Because $S$ is basic, $V(C_R) \cap (V(L_S) \cup \{a\}) \neq \varnothing$ so $a \in V(C_R)$.
	 Symmetrically, $v \in V(C_R)$.  Consequently a path in $C_R$ connects $a$ and $v$ in $G - R$.
	
	 By symmetry there is an $au$-path in $G-\{b,v\}-V(L_S)-V(L_T)$,
	 a $bv$-path in $G-\{a,u\}-V(L_S)-V(L_T)$, and a $bu$-path in $G-\{a,v\}-V(L_S)-V(L_T)$.
	 All of these paths must be internally vertex disjoint since otherwise there would be
	 an $ab$-path avoiding $\{u,v\}$ (which completes $L_S$ to $L_S^+$ in $C_T$) or a $uv$-path avoiding $\{a,b\}$ 
	 (which completes $L_T$ to $L_T^+$ in $C_S$), contradicting that
	 the heavy components of $S$ and $T$ are planar.

\begin{figure}[H]
	\centering
	\includegraphics[page=7]{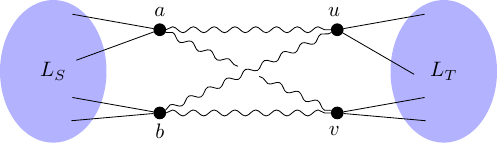}%
	\caption{Four internally vertex-disjoint paths connect $\{a,b\}$ to $\{u,v\}$.}
	\label{NoUVPathInX}
\end{figure}
	
Contract all four of these paths into edges connecting $\{a,b\}$ to $\{u,v\}$ (see \autoref{NoUVPathInX}). The graphs in \autoref{Disjoint2Cuts} are obstructions and one of them must be obtained via these contractions since $L_S^+-ab$ and $L_T^+-uv$ must contain $K_5-e$ or $K_{3,3}-e$. Because the original graph is a
minor-minimal apex obstruction, it must be one of these.\qedhere
\end{Specialindentpar}
\end{proof}

\begin{figure}[ht]
	\centering
	\includegraphics[page=8,scale=0.9]{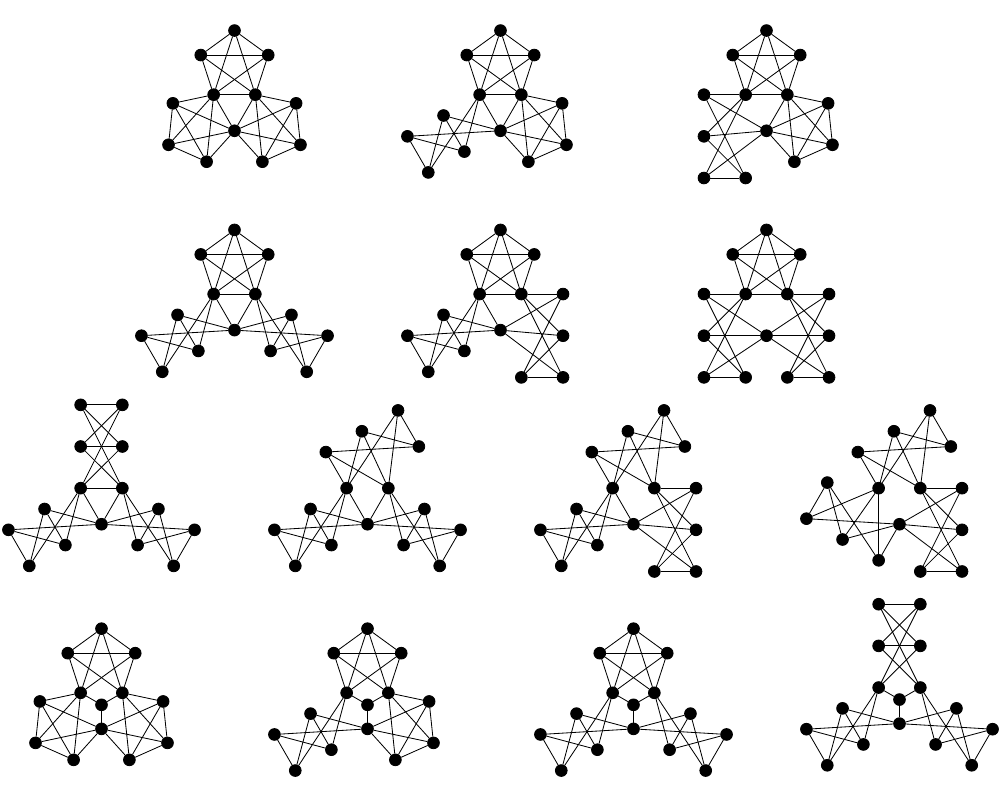}%
	\caption{The $14$ connectivity-$2$ apex obstructions with planar heavy components and at least three $2$-cuts, but no disjoint $2$-cuts.}
	\label{Intersecting2Cuts}
\end{figure}

\section{The 37 connectivity-2 apex obstructions with planar heavy components and more than one, but no disjoint, 2-cuts}
\label{sec:MultipleCuts}

We next consider the remaining connectivity-$2$ apex obstructions with with planar heavy components, more than one $2$-cut, and every pair of $2$-cuts intersect.  There are $37$ of these obstructions, $14$ of which have more than two $2$-cuts (see \autoref{Intersecting2Cuts}) and $23$ of which have exactly two $2$-cuts (see \autoref{MoreIntersecting2Cuts}).  

\begin{theorem} \label{Intersecting2CutTheorem}
Suppose $G \in {\cal F}$, $\kappa(G) = 2$, and every $2$-cut has a planar heavy component.
If $G$ has more than one $2$-cut and every pair of $2$-cuts intersect,
then $G$ is isomorphic to a graph in \autoref{Intersecting2Cuts} or \autoref{MoreIntersecting2Cuts}.
\end{theorem}
\begin{proof}
Suppose $S=\{a,b\}$ and $T=\{a,x\}$ are $2$-cuts of $G$. 
By Lemma~\ref{All2CutsAreBasic}, all $2$-cuts of $G$ are basic.
Let $C_S$ and $C_T$ be the heavy components of 
$G-S$ and $G-T$, respectively; let $L_S$ and $L_T$ be the corresponding light components.
Lemma~\ref{lightComponent} guarantees that $L_S^+$ and $L_T^+$ are
isomorphic to $K_5$, $K_{3,3}$, or $K_{3,3}+e$.
If $x$ were a vertex in $L_S$ then $G- \{a,x\}$ would be connected because $L_S^+$ is $3$-connected;
so $x \not\in V(L_S)$.  Symmetrically
$b \not\in V(L_T)$.  
Because $G-T$ has two components and one component contains $L_S \cup \{b\}$, it follows that
$L_S \cup \{b\} \subset C_T$ because $T$ is basic.  Symmetrically, $L_T \cup \{x\} \subset C_S$.
Consequently $L_S \cap L_T = \varnothing$.

\vspace{\baselineskip}
\noindent\textsc{Case 1}: $R=\{b,x\}$ is a $2$-cut of $G$.
\begin{Specialindentpar}{\parindent}  

Let $L_R$ denote the light component of $G-R$.  Reasoning as in the paragraphs above,
we find that $L_R^+$ is isomorphic to $K_5$, $K_{3,3}$, or $K_{3,3}+e$, and $L_R \cap L_S = \varnothing = L_R \cap L_T$.
Consider $H_a, H_b$, and $H_x$.  Note $b \in V(H_a)$ and $a \in V(H_b)$
because $C_S$ is planar.  Symmetrically $a \in V(H_x)$ and $x \in V(H_a)$ because $C_T$ is planar; and $b \in V(H_x)$ and $x \in V(H_b)$
because $C_R$ is planar.  In this case, $H_x \cap L_T = \varnothing = H_x \cap L_R$ because $\{a,b\} \subset V(H_x)$.  A symmetric
statement applies to $H_a$ and $H_b$.

Now let $P = V(G) - V(L_S) - V(L_T) - V(L_R) - \{ a,b,x\}$. 
 
\vspace{0.5\baselineskip}
\noindent\textsc{Case 1a}: $P = \varnothing$.
\begin{Specialindentpar}{18pt}

Because $H_x \cap L_T = \varnothing = H_x \cap L_R$ and $P = \varnothing$,
we find $V(H_x) \subset V(L_S) \cup S$.
This implies $G[V(L_S) \cup S] \cong L_S^+$.
It follows that $G[V(L_S) \cup S]$ is isomorphic to one of $K_5$, $K_{3,3}$, or $K_{3,3}+e$.
Symmetrically $G[V(L_T) \cup T]$ and $G[V(L_R) \cup R]$ are also isomorphic to one of these graphs.
There are ten non-isomorphic graphs that result from these possibilities,
and they each produce an obstruction, as shown in the top ten graphs of \autoref{Intersecting2Cuts}.
\end{Specialindentpar}

\vspace{0.5\baselineskip}
\noindent\textsc{Case 1b}: $P \neq \varnothing$.
\begin{Specialindentpar}{18pt}
   Let $W$ be a component of $P$.  Observe that $W$ must have vertices adjacent to all three vertices $a,b,$ and $x$
	 since otherwise the removal from $G$ of one of the $2$-cuts $R$, $S$, or $T$ would have more than two components,
	 contradicting Lemma~\ref{TwoComponents}.  Now consider the graph $H$ 
	 obtained from $G$ by contracting $W$ to a vertex $w$.  The vertex $w$ has neighborhood $\{a,b,x\}$ in $H$.
	  
	 If $ab, ax, bx \in E(G)$, then $H-w$ (hence $G-W$) has a minor isomorphic to one top ten graphs in \autoref{Intersecting2Cuts},
	 a contradiction.
	 Indeed if any one these three edges is in $E(G)$, then such a minor can be produced: for example,
	 if $ab \in E(G)$, then contract $H/wx$.
	 So we may assume $ab, ax, bx \notin E(G).$
	 We claim that $H[V(L_S) \cup S]$ is isomorphic to $K_5-e$ or $K_{3,3}-e$.
	 It suffices to prove that $H[V(L_S) \cup S]+ab \not\cong K_{3,3}+e$.  
	 If $H[V(L_S) \cup S]+ab \cong K_{3,3}+e$, then 
	 the contraction $H/wx$ again contains a minor isomorphic to one top ten graphs in \autoref{Intersecting2Cuts},
	 a contradiction. 
	
	 By symmetry we conclude that 
   $H[L_S \cup S]$, $H[L_T \cup T]$ and $H[L_R \cup R]$ are all isomorphic to $K_5-e$ or $K_{3,3}-e$.  
	 There are four possible graphs constructed using
   such components and they are shown at the bottom of \autoref{Intersecting2Cuts}.  Because these graphs are obstructions, 
   it follows that $H$ must be one of these graphs (and $P=W=\{w\}$); and $G \cong H$. 
\end{Specialindentpar}
\end{Specialindentpar}

\vspace{0.5\baselineskip}
\noindent\textsc{Case 2}: $R=\{b,x\}$ is not a $2$-cut of $G$.
\begin{Specialindentpar}{\parindent}
   As in Case 1, $b \in V(H_a)$ and $a \in V(H_b)$ since otherwise $C_S$ is non-planar. 
	 Also $x \in V(H_a)$.
	 Now let $P = V(G) - V(L_S) - V(L_T) - \{ a,b,x\}$.
   Note that $V(H_a) \cap V(L_T) = \varnothing$ since $\{b,x\} \subset V(H_a)$ and $\{a,x\}$ separates $b$ from $L_T$.
	 Similarly, $V(H_a) \cap V(L_S) = \varnothing$.
   Consequently $V(H_a) \subseteq V(P) \cup \{b,x\}$ and $P \neq \varnothing$.
  
	 Every component of $G[P]$ has a neighbor of $a$; otherwise $\{b,x\}$ is
   a $2$-cut (which is Case 1). 
	 So, there exists an edge $e = aw$ with $w$ in $P$.  

\begin{figure}[ht]
	\centering
	\includegraphics[page=9]{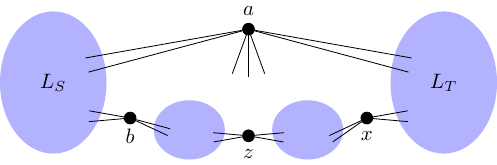}%
	\caption{Augmentations of $L_S$ and $L_T$ to Kuratowski subgraphs.}
	\label{BXNotA2Cut}
\end{figure}
   We claim that $w \in V(H_a)$.  If $w \notin V(H_a)$, then $G/e$ must have an apex, say $v$,
	 such that $v \in V(H_a) \subset P \cup \{b,x\}$.  Note that $v \notin \{b,x\}$ because
	 $G/e-b$ contains an augmentation of $L_T$ with $ax$-path through $P - b$, and similarly
	 $G/e-x$ contains an augmentation of $L_S$ with $ab$-path through $P - x$.  If $v \in P$
	 then $v$ would separate $b$ from $x$ in $H_a$ otherwise $L_S$ can be completed to a Kuratowski subgraph in $G/e-v$ with a $v$-avoiding $bx$-path
	 in $G[P \cup \{b,x\}]$ continuing to $a$ through $L_T$.  But $H_a$ is $2$-connected, so $v$ can not separate $b$ from $x$ in $H_a$.
	 So $G/e$ has no apex, a contradiction.
	
   Next we claim that $N(a) \cap (P \cup \{b,x\}) = \{w\}$.  Suppose, on the contrary, there is
   a vertex $v \in P \cup \{b,x\}$ such that $v \neq w$ and $f = av \in E(G)$.  Consider
   an apex $u$ for $G-f$.  Notice that $u \in V(H_a) \subseteq P \cup \{b,x\}$.
   Now $u\neq x$ because a Kuratowski subgraph exists by augmenting $G[V(L_S) \cup S]$
   with an $ab$-path from $a$ to $w$ through $P$ to $b$.  
   Similarly, $u \neq b$; so
   $u \in P$.  Because $u \notin \{b,x\}$, $u$ must separate $b$ from $x$ in $G-f$;
	 otherwise $L_S$ can be completed to a Kuratowski subgraph in $G-f-u$ with a $u$-avoiding $bx$-path
	 continuing to $a$ through $L_T$.
   However, $b$ and $x$ are in $H_a$ which is a subgraph $G[P \cup \{b,x\}]$.
	 Because $H_a$ is $2$-connected, $u$ can not separate $b$ from $x$, a contradiction.
	 So $N(a) \cap (P \cup \{b,x\}) = \{w\}$ (see \autoref{BXNotA2Cut2}).
	
\begin{figure}[H]
	\centering
	\includegraphics[page=10]{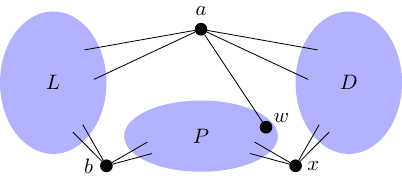}%
	\caption{$N(a) \cap (P \cup \{b,x\}) = \{w\}$.}
	\label{BXNotA2Cut2}
\end{figure}
	
	 Similar reasoning as the last paragraph shows that every edge of $G[P \cup \{b,x\}]$  must be an edge of $H_a$;
	 otherwise deleting this edge produces a graph with no apex.
	
	 Now suppose that there is an edge $f \in E(H_a)$ such that $|f \cap \{b,w,x\}| \leq 1$ and $f$ is incident
	 to a non-branch vertex of $H_a$.
	 Consider $G/f$; it contains $b,w,x$.  As before, $b$ and $x$ cannot be apex vertices
	 for $G/f$.  Furthermore, no vertex in $V(H_a)$ separates $b$ from $x$, so they can not be apex vertices.
	 Hence $G/f$ has no apex,
	 a contradiction.  We conclude that every edge of $H_a$ has branch vertices of $H_a$ as
	 endpoints or its endpoints are both in $\{b,w,x\}$.  
   
	 We claim that every vertex of $P \cup \{b,x\}$ is a branch vertex of $H_a$.
	 Because every edge of $G[P \cup \{b,x\}]$ is an 
	 edge of $H_a$ and every edge of $H_a$ has branch vertices of $H_a$ as
	 endpoints or its endpoints are both in $\{b,w,x\}$, it suffices to prove that
	 $b$, $w$, and $x$ are branch vertices of $H_a$.
	
	 If $b$ is not a branch vertex of $H_a$, then the two neighbors of $b$ in $H_b$ must be $w$ and $x$.
	 Because every edge of $G[P \cup \{b,x\}]$ is an 
	 edge of $H_a$, it follows that $b$ has only $w$ and $x$ as neighbors in $G[P \cup \{b,x\}]$.
	 This implies $\{w,x\}$ is a $2$-cut of $G$ separating $P-\{w\}$ from $b$; this $2$-cut is disjoint from $\{a,b\}$,
	 contradicting that $G$ has no disjoint $2$-cuts.  Hence $b$ must be a branch vertex of $H_a$.
	 By symmetry, $x$ is also a branch vertex of $H_a$.
	
	 Apply similar reasoning to show $w$ is a branch vertex.  If $w$ is not a branch vertex of $H_a$, then the two neighbor of $w$
	 in $H_a$ (and hence in $G[P \cup \{b,x\}]$) are $b$ and $x$.  This implies that $\{b,x\}$ is a $2$-cut of $G$,
	 (which is Case 1).
	 
   So every vertex of $P \cup \{b,x\}$ is a branch vertex of $H_a$.  This means that 
	 there are $23$ possible graphs shown in \autoref{MoreIntersecting2Cuts}.
   The three vertices in the center of these figures are, clockwise from the top of the three vertices,
   $a,x,b$.  The light component, $L_S$, appears to the left, and $L_T$ appears to the right.  The unsubdivided
   Kuratowski subgraph $H_a$ appear at the bottom of each graph drawing.
	 \qedhere

\end{Specialindentpar}
\end{proof}

\begin{figure}[htb]
	\centering
	\includegraphics[page=11,scale=0.9]{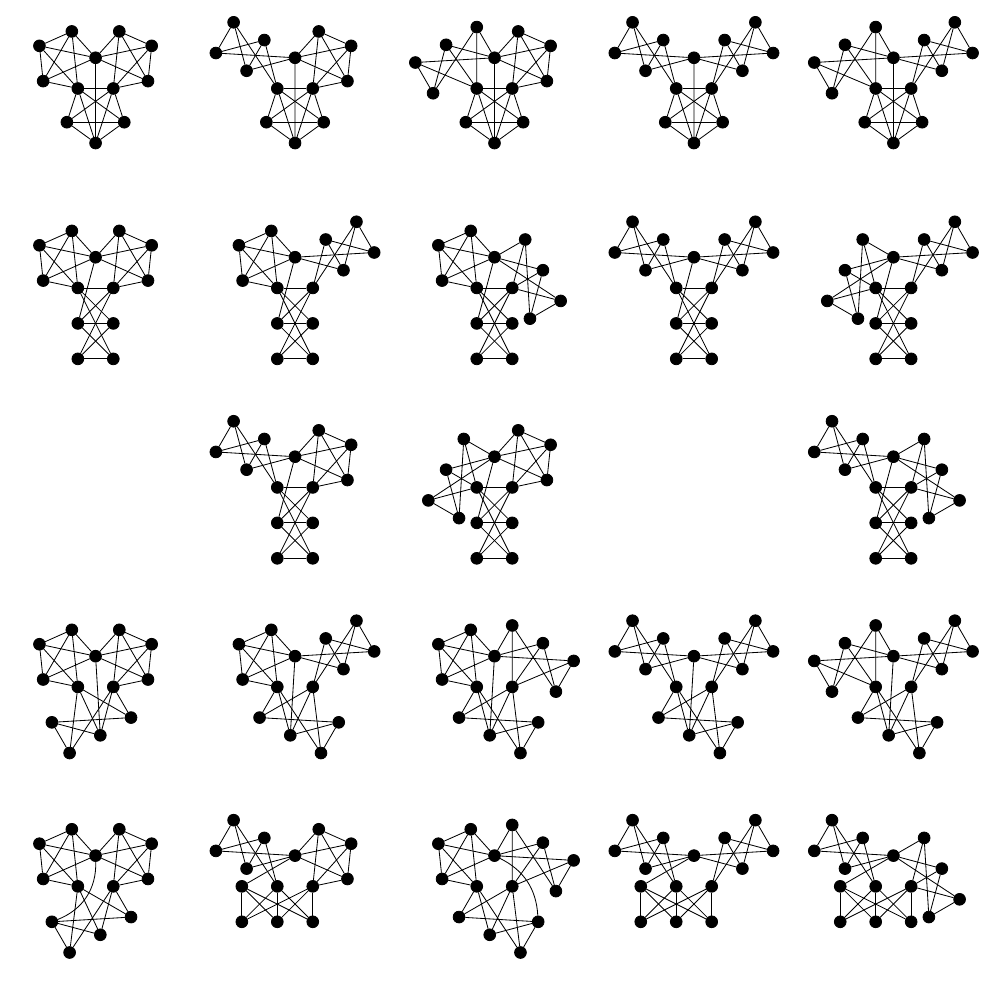}%
	\caption{The $23$ connectivity-$2$ apex obstructions with planar heavy components and exactly two $2$-cuts that intersect.}
	\label{MoreIntersecting2Cuts}
\end{figure}

\section{The \texorpdfstring{$72$}{72} connectivity-\texorpdfstring{$2$}{2} apex obstructions with a unique \texorpdfstring{$2$}{2}-cut having a planar heavy component}
\label{sec:UniqueCut}
We next consider the final collection of connectivity-$2$ apex obstructions, the obstructions with a unique $2$-cut whose removal produces
a planar heavy component.
These obstructions are considerably more difficult to characterize than the prior types, so we must present several structural results before tackling the final characterization.  

In this section we let $C$ denote the heavy component of an apex
obstruction after deleting the unique $2$-cut, with $C$ representing either this subgraph or its vertex set, as convenient. 
Similar expedient recastings also apply to other subgraphs and their vertex sets,
with context
alleviating any ambiguities.  

There are $72$ apex obstructions determined in this section.  Thirty-three of these obstructions (see \autoref{ThirtyThree}) have a heavy component $C$ and a unique $2$-cut $\{a,b\}$ such that $G[C \cup \{a,b\}]$ has a $2$-cut separating $a$ from $b$; the remaining thirty-nine (see \autoref{ThirtyNine}) have no such $2$-cut.  

The critical intermediate goal is to prove that there exist two Kuratowski subgraphs whose branch vertices cover every vertex of the unique $2$-cut and every vertex of the heavy component (Theorem~\ref{WeakWTheorem}). Adding the light-component Kuratowski subgraph to these two makes three Kuratowski subgraphs whose branch vertices cover the entire vertex set of the apex obstruction,
a property that reduces the problem of characterizing all of these remaining obstructions to a very small
number of cases which we complete in the final subsection.

Many preliminary lemmas are needed for this final analysis.  To simplify the presentation, we adopt the following notation and assumptions in this section.

\begin{assumptions} Standard assumptions for this section:
\begin{itemize}[nosep]
\item $G$ is a minor-minimal non-apex graph.
\item $G$ has connectivity two and a unique $2$-cut $S = \{a,b\}$.
\item $G-S$ has two components, the heavy component $C$ and the light component $L$.
\item $C$ is a planar graph.
\end{itemize}
\label{assumptions}
\end{assumptions}

\subsection{Some preliminary lemmas for the final case}

This subsection presents several lemmas establishing basic structure.

\begin{lemma}\label{LastCaseStructure}
Under Assumptions~\ref{assumptions}:
\begin{enumerate}
\item If $H_a$ and $H_b$ are any Kuratowski subgraphs of $G$ avoiding $a$ and $b$, respectively, then
$b \in V(H_a)$ and $a \in V(H_b)$.
\item $ab \not\in E(G)$.
\item  If $H_a$ and $H_b$ are any Kuratowski subgraphs of $G$ avoiding $a$ and $b$, respectively, then
\[E(G[C \cup \{a,b\}]) \subseteq E(H_a \cup H_b).\]
\item Any $2$-cut of $G[C \cup \{a,b\}]$ that separates $a$ from $b$ does not induce an edge.
\end{enumerate}
\end{lemma}
\begin{proof} (i) Because $C$ is planar, any Kuratowski subgraph of $G[C \cup \{a,b\}]$ must contain at least one vertex
from $\{a,b\}$.  Consequently any Kuratowski subgraph avoiding $a$ must contain $b$ and vice versa.

\noindent
(ii) If $e=ab$ were an edge of $G$, then $G-e$ would have an apex $z \in V(H_a) \cap V(H_b) \subseteq C$.  Indeed $z$ would have to separate $a$ from $b$ in $C$; otherwise an $ab$-path in $G-e$ could complete $L$ to $L^+$.  Consequently, either $\{a,z\}$ or $\{b,z\}$ would be another $2$-cut in $G$.

\noindent
(iii)  
By Theorem~\ref{edgecover}, $E(C) \subseteq E(H_a) \cup E(H_b)$.  Since (i) states that $ab \not\in E(G)$, to prove (iii) it suffices (by symmetry of $a$ and $b$) to prove that every edge
$aw$ with $w \in V(C)$ is also an edge in $H_a$ or $H_b$. We argue by contradiction.

Assume that there is an edge $f=aw$, with $w \in V(C)$ such that
$f$ is not in $E(H_a) \cup E(H_b)$.  Consider $G-f$; it has an apex vertex $z$.  Note that
$H_a$ and $H_b$ are still intact in $G-f$, so $z \in V(H_a) \cap V(H_b) \subset C$; in particular,
$z \not\in \{a,b\}$.  Now $z$ must separate $a$ from $b$ in 
$G[C \cup \{a,b\}] - f$, otherwise an $ab$-path avoiding $z$ and $f$ could
be used to augment $L$ to a Kuratowski subgraph $L^+$, contradicting that $G-f-z$ is planar.
Define $A$ to be the set of vertices that are connected to $a$ via a $z$-avoiding path
in $G[C \cup \{a,b\}] - f$.  
Because $a \in V(H_b)$ it has at least two neighbors in $H_b$; at least one of these is not $z$.
Consequently, $A \neq \varnothing$.
Now observe that  $\{a,z\}$ is $2$-cut of $G$ separating $b$ from vertices in $A$, contradicting
that $G$ has a unique $2$-cut. 

\noindent
(iv)  Assume, to the contrary, that a $2$-cut, $\{w,x\}$ of $G[C \cup \{a,b\}]$ separates $a$ from $b$ and $f = wx \in E(G)$.
Partition the vertices of $C-\{w,x\}$ into two sets:
\[
   A = \{ v \in C - \{w,x\} : \mbox{there is a $va$-path avoiding $\{w,x\}$ in $G[C \cup \{w,x,a\}]$}\}\]
and
\[
   B = \{ v \in C - \{w,x\} : \mbox{there is a $vb$-path avoiding $\{w,x\}$ in $G[C \cup \{w,x,b\}]$}\}\]
(see \autoref{Branch2}).

\begin{figure}[H]
	\centering
	\includegraphics[page=12]{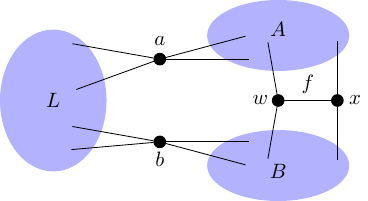}%
	\caption{$V(H_b) \subseteq A \cup \{a,w,x\}$ and $V(H_a)\subseteq B \cup \{b,w,x\}$.}
	\label{Branch2}
\end{figure}

Let $H_a$ and $H_b$ be minimum-sized Kuratowski subgraphs of $G$ avoiding $a$ and $b$, respectively.
Because $V(H_b) \subset C \cup \{a\}$ and the branch vertices of $H_b$ must all be on one side of $\{w,x\}$, it follows from $wx \in E(G)$ and the minimum size of $H_b$ that $V(H_b) \subseteq A \cup \{a,w,x\}$.  A similar argument shows $V(H_a) \subseteq B \cup \{b,w,x\}$.  In particular, $A \neq \varnothing \neq B$.

Observe that $w$ must have neighbors in $A$ and $B$ since $\{a,x\}$ and $\{b,x\}$ are
not $2$-cuts of $G$; similarly $x$ must have neighbors in $A$ and $B$.
Note too that there must be a $wx$-path in $G[A \cup \{w,x\}] - f$; otherwise $\{a,w\}$ or $\{a,x\}$ would
be another $2$-cut of $G$.  Symmetrically there must be a $wx$-path in $G[B \cup \{w,x\}] - f$.
Finally there must be two vertex-disjoint paths from $a$ to $\{w,x\}$ in $G[A \cup \{a,w,x\}]$
since $\{a,b\}$ is the only $2$-cut of $G$.  Similarly, 
there must be two vertex-disjoint paths from $b$ to $\{w,x\}$ in $G[B \cup \{b,w,x\}]$.

Now we claim that $G-f$ has no apex. Suppose, to the contrary, that $y$ is an apex for $G-f$.  
If $y \in A \cup \{w,x\}$, then there remains an $a$ to $\{w,x\}$-path that can be completed
to an $ab$-path in $G[C \cup \{a,b\}]- f -y$, thus augmenting $L$ to $L^+$ in $G-f-y$.
A similar argument applies if $y \in B \cup \{w,x\}$.
So it remains to consider $y \in \{a,b\} \cup L$.
If $y \in \{a\} \cup L$, then consider $H_a$ in $G-f-y$.  Recall that $V(H_a) \subseteq B \cup \{b,w,x\}$, so
 $H_a$ is only possibly missing the edge $f$ in $G-f-y$.
Now a $wx$-path whose internal vertices are entirely in $A$ can substitute for the edge $f$ in $H_a$.  Hence no
vertex in $\{a,b\} \cup L$
is an apex for $G-f$.  A similar argument show that $y=b$
is not an apex in $G-f$.
Hence $G-f$ has no apex, a contradiction.
\end{proof}


The next lemma embodies fundamental arguments to which we shall often appeal.
The lemma may be interpreted as stating that $J = G[C \cup \{a,b\}]$ is a minor-minimal
double apex graph with roots $a$ and $b$ (see Theorem \ref{DoubleApexFormulation}).

\begin{lemma}\label{BasicLemma} Assume Assumptions~\ref{assumptions} and $J = G[C \cup \{a,b\}]$.
For any edge $e \in E(J)$,
\begin{enumerate}
\item\label{BasicLemma.delete} $J-e-a$ is planar or $J-e-b$ is planar (or both), and
\item\label{BasicLemma.contract} $J/e-a$ is planar or $J/e-b$ is planar (or both).
\end{enumerate}
In part \ref{BasicLemma.contract} the notation means that if $e$ is incident to vertex $a$ (resp. $b$) in $J$, then $e$ is contracted in $J/e$
 to form a new vertex also labeled $a$ (resp. $b$).
\end{lemma}
\begin{proof}
   \ref{BasicLemma.delete} Assume, to the contrary, that $J-e-a$ and $J-e-b$ are both non-planar. 
	Consider Kuratowski subgraphs $H_a$ and $H_b$ in $J-e$ avoiding $a$ and $b$, respectively.  
Now observe that $e \notin E(H_a) \cup E(H_b)$, contradicting Lemma~\ref{LastCaseStructure} part~(iii).

\ref{BasicLemma.contract} Assume, to the contrary, that $J/e-a$ and $J/e-b$ are both non-planar. 
Consider $G/e$; it has an apex $z$.  Because $J/e-a$ and $J/e-b$ are both non-planar, $z \notin \{a,b\}$.
Now $z$ must separate $a$ from $b$ in $J/e$, otherwise a $z$-avoiding path in $G/e$ could complete $L$ to $L^+$,
contradicting that $z$ is an apex for $G/e$.  Because $\{a,z\}$ and $\{b,z\}$ are not $2$-cuts of $G$,
it follows that $z$ is the contracted vertex created by contracting the edge $e$.  Consequently the vertices
in $e$ form a $2$-cut in $J$ separating $a$ from $b$, contradicting Lemma~\ref{LastCaseStructure} part (iv).
\end{proof}

The next lemma gives information about edges shared by Kuratowski subgraphs missing $a$ and $b$.

\begin{lemma} \label{BothBranch}
Assume Assumptions~\ref{assumptions}. Suppose $H_a$ and $H_b$ are any Kuratowski subgraphs avoiding $a$ and $b$ respectively.
If $uv = e \in E(H_a) \cap E(H_b)$, then
\begin{enumerate}
\item $u$ and $v$ are both branch vertices of $H_a$, or
\item $u$ and $v$ are both branch vertices of $H_b$,  
\end{enumerate}
or both.
\end{lemma}
\begin{proof} Let $J = G[C \cup \{a,b\}]$.  If $u$ is a branch vertex of neither $H_a$ nor $H_b$, then
$H_a$ and $H_b$ remain non-planar after contracting $e$.  Consequently $J/e-a$ and $J/e-b$ are both non-planar,
contradicting Lemma~\ref{BasicLemma} part~\ref{BasicLemma.contract}. So $u$ is a branch vertex of $H_a$ or $H_b$;
similarly $v$ is a branch vertex of $H_a$ or $H_b$.  If $u$ or $v$ is a branch vertex of both, then (i) or (ii) follows.
So we may assume, to the contrary, that $u$ is a branch vertex of $H_a$ but not a branch vertex of $H_b$ and, symmetrically,
$v$ is a branch vertex of $H_b$ but not a branch vertex of $H_a$.  In this case, again 
$H_a$ and $H_b$ remain non-planar after contracting $e$ so $J/e-a$ and $J/e-b$ are both non-planar,
again contradicting Lemma~\ref{BasicLemma} part~\ref{BasicLemma.contract}.
\end{proof}

The next lemma give a powerful tool to prove that vertices are branch vertices of Kuratowski subgraphs
avoiding $a$ and $b$.

\begin{lemma} \label{Proposition3}
Assume Assumptions~\ref{assumptions}. Suppose $H_a$ and $H_b$ are any Kuratowski subgraphs avoiding $a$ and $b$ respectively.
\begin{enumerate}
\item If $v \in V(H_a) - (V(H_b) \cup \{a,b\})$, then $v$ and all of its neighbors are branch vertices of $H_a$.
\item If $v \in V(H_b) - (V(H_a) \cup \{a,b\})$, then $v$ and all of its neighbors are branch vertices of $H_b$.
\end{enumerate}
\end{lemma}
\begin{proof} By symmetry it suffices to prove (i).  Let $v \in V(H_a) - (V(H_b) \cup \{a,b\})$.  Assume, to the
contrary, that $w$ is $v$ or a neighbor of $v$ but $w$ is not a branch vertex of $H_a$.
If $w=v$, then because $w$ has only degree two in $H_a$, there is an edge incident to $w$ (since $\delta(G) \geq 3$, Lemma~\ref{MinDegree3}) that is not in $E(H_a) \cup E(H_b)$, contradicting Lemma~\ref{LastCaseStructure} part~(iii).
So $w \neq v$.  Let $e = wv$ and $J= G[C \cup \{a,b\}]$.

Because $w$ is not branch of $H_a$, the Kuratowski $H_a$ survives in $J/e$.  Also, $H_b$ survives in $J/e$ because $v \notin 
V(H_b)$.  Consequently $J/e-a$ and $J/e-b$ are non-planar, contradicting Lemma~\ref{BasicLemma} part~\ref{BasicLemma.contract}.
\end{proof}


\subsection{Vertices in the unique 2-cut must be branch vertices}

The main goal in this subsection is to prove that vertex $a$ must be a branch vertex of any Kuratowski subgraph avoiding $b$ and,
symmetrically, $b$ must be a branch vertex of any Kuratowski subgraph avoiding $a$ (statement of Corollary~\ref{BranchVertices}).
To simplify statements in this subsection we focus just on proving the latter claim since the former claim follows from symmetry.
So the reader should keep in mind the symmetrical consequences of the results that follow.

First we focus on the connectivity of $C$.

\begin{lemma} \label{TwoConnectivityC}
Assume Assumptions~\ref{assumptions} and $H_a$ is any Kuratowski subgraph avoiding $a$.  
If $b$ is not a branch vertex of $H_a$, then $C$ is $2$-connected.
\end{lemma}
\begin{proof}  Let $H_b$ be any Kuratowski subgraph missing $b$.
If $b$ is not a branch vertex of $H_a$, then it has exactly two neighbors in $H_a$.
Since all edges of $G[C \cup \{a,b\}]$ must be covered by $E(H_a) \cup E(H_b)$
(Lemma~\ref{LastCaseStructure} part~(iii)), it follows that $N(b) \cap C = \{x,y\}$, for some $x,y \in C$.
Suppose, to the contrary, that $c$ is a cut vertex of $C$.  If $c=x$, $c=y$, or $x$ and $y$ are in the same component of $C - c$,
then $\{a,c\}$ is another $2$-cut of $G$, a contradiction.  Therefore $c$ separates $x$ from $y$ in $C$.
In $G[C \cup \{b\}]$, the set $\{b,c\}$ is a $2$-cut so all branch vertices of $H_a$ must be on one side of this $2$-cut.
If one side of this $2$-cut contains no branch vertices of either $H_a$ or $H_b$, then (since each side is non-empty) there is an edge $e$ incident to 
$c$ that can be contracted that preserves $H_a$ and $H_b$.
Consequently, $G[C \cup  \{a,b\}]/e -a $ and  $G[C \cup  \{a,b\}]/e -b $ are non-planar, contradicting
Lemma~\ref{BasicLemma} part~\ref{BasicLemma.contract}.  So, each side of $G[C \cup \{b\}] - c$ has branch vertices.
That is, the branch of $H_a$ are all in one side, and the branch vertices of $H_b$, except possibly $a$, are on the other side.
Consequently, the vertex $a$ can have at most one neighbor, call it $z$, in one of the sides, which without loss of generality, contains $y$.  If there is a $cz$-path avoiding $y$, then $by$ is an edge that can be contracted preserving $H_a$ and $H_b$;
that is, $G[C \cup  \{a,b\}]/by -a $ and  $G[C \cup  \{a,b\}]/by -b $ are non-planar, contradicting
Lemma~\ref{BasicLemma} part~\ref{BasicLemma.contract}.  So every $cz$-path contains $y$.  If $y = z$, then $\{c,y\}$ or $\{b,y\}$ is another $2$-cut of $G$.  If $y\neq z$, then $\{a,y\}$ is a $2$-cut separating $c$ from $z$.
\end{proof}

Now we are ready to prove a claim that plays a major role in the final characterization of the 
connectivity-$2$ apex obstructions.  The claim states
that, under the Assumptions~\ref{assumptions}, 
the vertices in the unique $2$-cut, $a$ and $b$, must be branch vertices of all of their Kuratowski witnesses.
This is the statement of Corollary~\ref{BranchVertices}, which will be a consequence of Proposition~\ref{ExistenceOfKuratowskis}.  
Though Proposition~\ref{ExistenceOfKuratowskis} proves a seemingly weaker existential claim, its significance is
indicated by its long proof in which several $3$-connected apex obstructions play a vital role, including the Petersen-family graphs $M$, $Y^-$, $P_7$ and the Petersen graph itself.

\begin{proposition}\label{ExistenceOfKuratowskis}
If Assumptions~\ref{assumptions} are satisfied, then
there exists a Kuratowski subgraph avoiding vertex $a$ in which $b$ is a branch vertex.
\end{proposition}
\begin{proof} Suppose, to the contrary, that $b$ is not a branch vertex of any 
Kuratowski subgraph avoiding $a$. Let $H_a$ be a Kuratowski subgraph of $G$ avoiding $a$.
Necessarily $b \in V(H_a)$ (Lemma~\ref{LastCaseStructure} part~(i)), so $b$ has at least two neighbors in $H_a$.
If $|N_G(b) \cap C| > 2$, then
some edge in $G[C \cup \{b\}]$ incident to $b$ does not appear in $H_a$ or any Kuratowski subgraph of $G$ avoiding $b$, 
contradicting Lemma~\ref{LastCaseStructure} part~(iii). 
Consequently we may assume $|N_G(b) \cap C| = 2$; say $N_G(b) \cap C = \{x,y\}$.

Fix a plane embedding
of $C$.  
Let $K=V(H_a) - b$.  The plane embedding of $C$ includes
a plane embedding of $G[K]$.  Because $H_a-b$ is a subdivision of $K_5-e$ or $K_{3,3}-e$, a plane embedding
of $H_a-b$ is unique; consequently, a plane embedding of $G[K]$ is unique.  
Define $A$ to be the component of $G[C \cup \{a\} - K] $ that contains the vertex $a$,
and let $\overline{a}$ be the vertex obtained by contracting $A$ in this graph to a single vertex.

Clearly $G[C \cup \{a\}]$ is not planar.  It is possible that this graph becomes planar after contracting $A$ to $\overline{a}$.
We focus on $H_a - b$, a subgraph of $G[C \cup \{a\}]$. 
If it is possible to append $\overline{a}$ (and its incident edges to neighbors in $H_a - b$) to the plane embedding of $H_a-b$
without introducing a crossing in the plane (i.e. extend the plane embedding of $H_a-b$ to a plane
embedding of $H_a-b + \overline{a}$), then we say that {\em $\overline{a}$ hits only one face}; otherwise
{\em $\overline{a}$ hits multiple faces}.

The following steps provide contradictions that 
complete the proof of the proposition.

\vskip 0.35cm
\noindent
\fbox{\textsc{Step} 1: For any $H_a$, $\overline{a}$ hits only one face of the plane embedding of $H_a - b$.}

Suppose, on the contrary, there is a choice of $H_a$ such 
that $\overline{a}$ hits multiple faces of the plane embedding of $H_a - b$.
There are two cases according to whether $H_a-b$ is a subdivision of $K_5-e$ or $K_{3,3}-e$.

Consider first the case in which $H_a-b$ is a subdivision of $K_5-e$.
The planar embedding of this subdivision is unique; it is shown in \autoref{K5MinusEdgeInPlane}.
\begin{figure}[H]
	\centering
	\includegraphics[page=13]{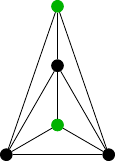}%
	\caption{The case when $H_a-b$ is a subdivision of $K_5-e$; neighbors of $b$ are green.}
	\label{K5MinusEdgeInPlane}
\end{figure}
Lemma~\ref{TwoConnectivityC} implies that $\overline{a}$ must have at least two neighbors in $H_a-b$.
Suppose that $\overline{a}$ has only two such neighbors.  There are seven non-isomorphic ways  (see \autoref{ABarAttachesToK5Minus}) that
$\overline{a}$ has exactly two neighbors in this subdivision of $K_5-e$ so that the resulting graph is not planar.
\begin{figure}[H]
	\centering
	\includegraphics[page=14]{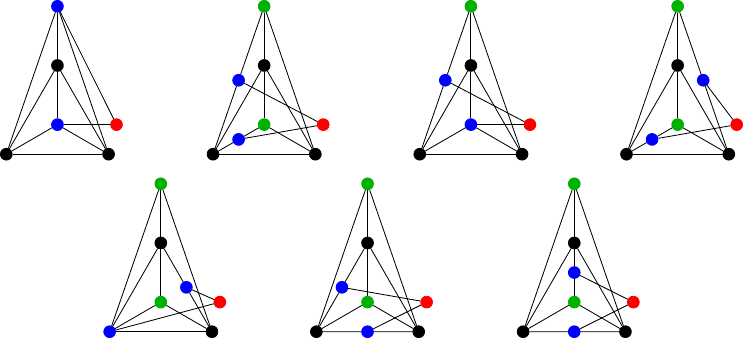}%
	\caption{In the case $H_a-b$ is a subdivision of $K_5-e$, there are seven ways that $\overline{a}$ has two neighbors in $H_a-b$ not both of which are in the same face.}
	\label{ABarAttachesToK5Minus}
\end{figure}
The three right-most cases on the top row of \autoref{ABarAttachesToK5Minus} all include neighbors of $\overline{a}$ (shown in blue)
that can be contracted along edges of $H_a - b$ so that one is contracted to $x$ and the other is contracted to $y$.  These
contractions produce a proper minor of $J = G[C \cup \{a,b\}]$ in which the contracted $H_a - b$, still a subdivision of $K_5 - e$, extends
to a subdivision of $K_5$ by adding either $b$ or $\overline{a}$.  
This contradicts Lemma~\ref{BasicLemma} part~\ref{BasicLemma.contract}, so these cases cannot occur.

The left-most case on the top row of \autoref{ABarAttachesToK5Minus} implies that $\{x,y\}$ is another $2$-cut
of $G$, a contradiction.  As an aside, this last configuration actually
has disjoint $2$-cuts, $\{x,y\}$ and $\{a,b\}$, showing that this case would produce
obstructions shown in \autoref{Disjoint2Cuts}. 

The cases shown on the bottom row of \autoref{ABarAttachesToK5Minus} can be dismissed as follows.
Consider adding the vertex $b$ to each graph along with the edge $\overline{a}b$ (which corresponds to contracting 
$L$ in $G$ to the edge $ab$).  \autoref{CompleteK5Minus} shows the resulting graphs (where now vertex $b$ is shown in red).
The leftmost graph of \autoref{CompleteK5Minus} is a $3$-connected apex obstruction 
from the Petersen family; it is usually called $M$.  
Because $M$ cannot appear as a proper minor of $G$, this case does not occur.  The other two graphs of \autoref{CompleteK5Minus}
contain $M$ as a proper minor, so those cases too cannot occur.

\begin{figure}[H]
	\centering
	\includegraphics[page=15]{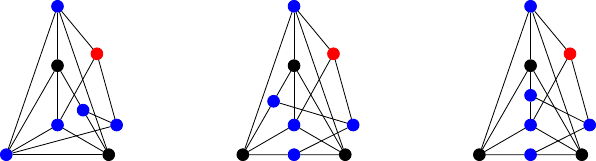}%
	\caption{The leftmost graph is $M$, a Petersen family apex obstruction.  The other two graphs properly contain $M$ as a minor. The vertex $b$ is shown in red.}
	\label{CompleteK5Minus}
\end{figure}

To complete the analysis of the case in which $H_a-b$ is a subdivision of $K_5-e$, 
suppose now that $\overline{a}$ has at least three neighbors in $H_a-b$.
Observe that $\overline{a}$ cannot have two neighbors that attach as in \autoref{ABarAttachesToK5Minus} since
each of these cases produces a contradiction to Lemma~\ref{BasicLemma} (parts \ref{BasicLemma.delete} and \ref{BasicLemma.contract}) or a minor of $M$.
Consequently, we may assume that the three
neighbors of $\overline{a}$ are not all on one face, but any two of them are on a single face.
There are four non-isomorphic ways to select such neighbors for $\overline{a}$ (see \autoref{K5Plus3}).
\begin{figure}[H]
	\centering
	\includegraphics[page=16]{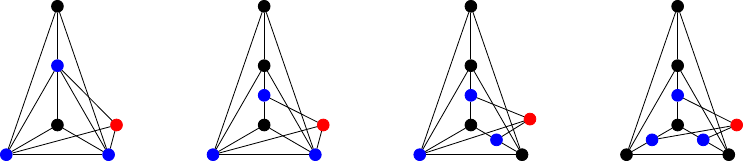}%
	\caption{The four non-isomorphic ways that $\overline{a}$ (shown in red) has three neighbors (shown in blue) in $H_a-b$ (which is a subdivision of $K_5-e$) so that all three neighbors are not on the same face, but any two are on a single face.}
	\label{K5Plus3}
\end{figure}
One way is to choose 
three neighbors of $\overline{a}$ to be the three degree four vertices of $H_a-b$.
The remaining three cases can be contracted to this case.  Now
adding to this graph the vertex $b$ and the
edge $\overline{a}b$, as before, produces the graph shown
in \autoref{YK5Minus}.  It is an apex obstruction from the Petersen family; it is called $Y^-$.
It cannot occur as a minor of $G$, so this completes the analysis of Step 1 
when $H_a-b$ is a subdivision of $K_5-e$

\begin{figure}[H]
	\centering
	\includegraphics[page=17]{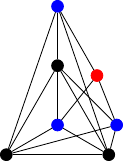}%
	\caption{The graph $Y^-$, a Petersen family apex obstruction.  Vertex $b$ is shown in red.}
	\label{YK5Minus}
\end{figure}

Consider next the case in which $H_a-b$ is a subdivision of $K_{3,3}-e$.
The planar embedding of this subdivision is unique; it is shown in \autoref{K33MinusEdgeInPlane}.
\begin{figure}[H]
	\centering
	\includegraphics[page=18]{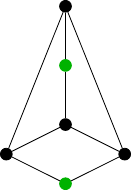}%
	\caption{The case when $H_a-b$ is a subdivision of $K_{3,3}-e$; neighbors of $b$ are green.}
	\label{K33MinusEdgeInPlane}
\end{figure}

The analysis now follows along similar reasoning as the case in which $H_a-b$ is a subdivision of $K_5-e$.
Suppose that $\overline{a}$ has only two neighbors in $H_a-b$.  There are four non-isomorphic ways  (see \autoref{ABarAttachesToK33Minus}) that
$\overline{a}$ has exactly two neighbors in this subdivision of $K_{3,3}-e$ so that the resulting graph is not planar.
\begin{figure}[H]
	\centering
	\includegraphics[page=19]{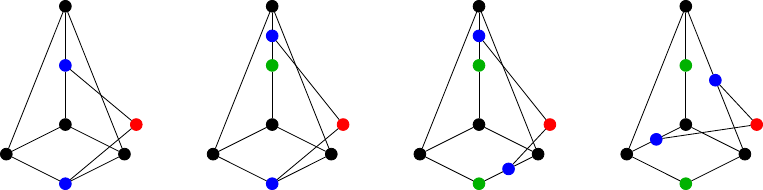}%
	\caption{In the case $H_a-b$ is a subdivision of $K_{3,3}-e$, there are four ways that $\overline{a}$ has two neighbors in $H_a-b$ not both of which are in the same face.}
	\label{ABarAttachesToK33Minus}
\end{figure}
In the left three graphs of \autoref{ABarAttachesToK33Minus} one can reason, as before in the $K_5-e$ case, 
that either $\{x,y\}$ is another $2$-cut of $G$ or
one neighbor 
of $\overline{a}$ can be contracted to $x$ and the other to $y$ so that 
a contradiction to Lemma~\ref{BasicLemma} part~\ref{BasicLemma.contract} occurs.  So $\overline{a}$ has no two neighbors of this type,
leaving only the exclusion of the graph at the right of \autoref{ABarAttachesToK33Minus}.
This graph extends with $b$ and the edge $\overline{a}b$ to the graph shown
in \autoref{PetersenGraph}.  This graph is isomorphic to the Petersen graph, an apex obstruction, so cannot be a proper minor of $G$.
\begin{figure}[htb]
	\centering
	\includegraphics[page=20]{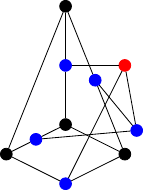}%
	\caption{In the case $H_a-b$ is a subdivision of $K_{3,3}-e$, a graph isomorphic to the Petersen graph emerges. Vertex $b$ is red.}
	\label{PetersenGraph}
\end{figure}


So it remains to consider the case in which $\overline{a}$ has three (or more) neighbors in $H_a-b$ such that,
not all are in a single face, but any two of them are on a single face.  
There are ten non-isomorphic ways $\overline{a}$ has three such neighbors (see \autoref{K33Plus3}).

\begin{figure}[H]
	\centering
	\includegraphics[page=21]{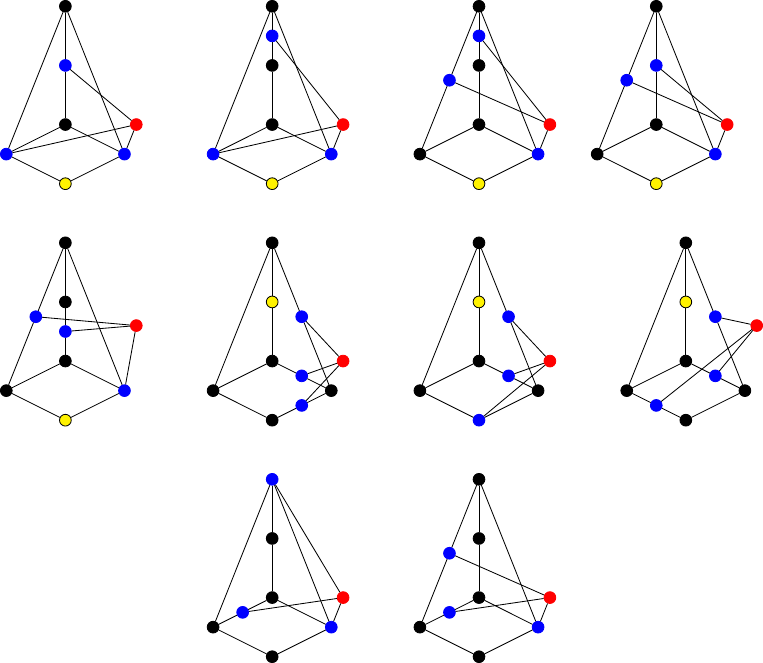}%
	\caption{The ten non-isomorphic ways $\overline{a}$ (shown in red) has three neighbors in $H_a-b$ (which is a subdivision of $K_{3,3}-e$) so that all three neighbors (shown in blue) are not on the same face, but any two are on a single face.}
	\label{K33Plus3}
\end{figure}

In the eight graphs in the top two rows of \autoref{K33Plus3} a vertex, either $x$ or $y$, can be deleted and the remaining graph remain non-planar; in each case the deletable vertex is shown in yellow. Consequently in these graphs either $bx$ or $by$, according to whether $x$ or $y$ is deletable, is an edge that contradicts Lemma~\ref{BasicLemma} (part~\ref{BasicLemma.contract}). So these eight graphs cannot occur.

The bottom right graph of \autoref{K33Plus3} can be contracted to the bottom left graph, so it remains to dismiss the bottom left graph. As before, add vertex $b$ and the edge $\overline{a}b$ to this graph produces the graph shown in \autoref{P7}.
This graph is a graph isomorphic to the graph $P_7$, a Petersen-family minor. Because it is an apex obstruction it cannot occur as a proper minor of $G$. This completes the proof of Step $1$.

\begin{figure}[htb]
	\centering
	\includegraphics[page=22]{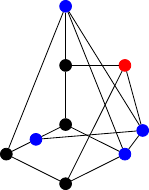}%
	\caption{In the case $H_a-b$ is a subdivision of $K_{3,3}-e$, a graph isomorphic to the graph $P_7$, a Petersen family apex obstruction, emerges. Vertex $b$ is red.}
	\label{P7}
\end{figure}


By Step 1, we may assume that $\overline{a}$ has neighbors in only one face of the plane embedding of $H_a - b$.
Let $F$ be a face of the plane embedding of $H_a - b$ such that $\overline{a}$ has only neighbors in it.
There could be two such faces, a technicality addressed later (see Step 5).
Without loss of generality, the face $F$ is a face not incident to $y$ (see \autoref{WLOGTheFaceF}).
Recall that $K=V(H_a) - b$ and $A$ is the component of $G[C \cup \{a\} - K] $ that contains the vertex $a$.
Let $int(F)$ denote all of the vertices of $C - K$ that appear in the interior of $F$.  More precisely, $int(F)$ consists of
the vertices in the interior of the region of
the plane avoiding $y$ that is bounded by $F$ in our fixed plane embedding of $C$.

Choose $H_a$ that minimizes the number of vertices that are on the boundary or the interior of $F$.  
We may also assume that embedding has been fixed to minimize the number of crossings
produced when $a$ is reinserted into the face $F$.

\begin{figure}[htb]
	\centering
	\includegraphics[page=23]{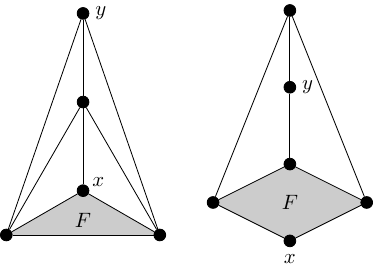}%
	\caption{Whether $H_a - b$ is a subdivision of $K_5-e$ or $K_{3,3}-e$, we may assume $F$ is a face that is not incident to $y$ in the unique plane embedding of $H_a - b$.}
	\label{WLOGTheFaceF}
\end{figure}

Define $C_1,\ldots,C_t$ to be the components of $G[C \cup \{a\} - K] $ that intersect $int(F) \cup \{a\}$.
Note that $A$ is one of these components; we assume $C_1=A$.
Observe that  $|N(C_i) \cap F| \geq 3$, for $i>1$; otherwise there would
be another $2$-cut in $G$.  This argument does not apply to the component $C_1$ because it contains the vertex $a$ which
has neighbors in $L$.  However, Lemma~\ref{TwoConnectivityC} implies that $\overline{a}$ must have at least two neighbors in $H_a-b$; hence $C_1$ must have at least two neighbors on $F$.

Define $S$ to be the set of vertices on $F$ with at least one neighbor in $int(F) \cup \{a\}$.  
  Note that $|S|\geq 2$ since if $S=\{v\}$, then $\{b,v\}$ is another $2$-cut
of $G$.  An important observation is that, by Lemma~\ref{Proposition3}, all the vertices in $int(F)$ and their neighbors on $F$
are all branch vertices of $H_b$.  So all vertices in $S$ are branch vertices of $H_b$, except possibly vertices whose
only neighbor in $int(F) \cup \{a\}$
is $a$.

\vskip 0.35cm
\noindent
\fbox{\textsc{Step} 2: $|C_i|=1,$ for all $i=1,\ldots,t$.}

Assume, to the contrary, that $|C_i|>1$ for some $i=1,\ldots,t$.  
Let $u \in C_i - \{a\}$ and $v \in C_i$ be chosen so that $uv$ is an edge.  Let $e=uv$.
Lemma~\ref{Proposition3} implies that all of the vertices in $C_i$ and the neighbors of $u$
must be branch vertices of $H_b$.

Recall that $J=G[C \cup \{a,b\}]$.  
It suffices to show that there exists a Kuratowski minor in $G[C \cup \{a\}]$ avoiding $y$,
since then 
$yb$ can be contracted so that $J/yb-a$ and $J/yb-b$ are both non-planar,
contradicting Lemma~\ref{BasicLemma} part~\ref{BasicLemma.contract}.
Consider $J/e$ and let $w$ be the composite vertex in $J/e$ that
results from identifying $u$ and $v$. 
Because $H_a$ is untouched by the contraction of $e$, 
Lemma~\ref{BasicLemma} part~\ref{BasicLemma.contract} implies that $J/e-b$ is planar.  Because $H_a-b$ has a unique planar embedding and it is a 
subgraph of $J/e-b$, the cycle $F$ must separate $y$ from $w$ in any plane embedding of $J/e-b$.
In particular, the faces containing $w$ do not contain $y$.  
Let $W$ be the plane graph formed by the union of the faces containing $w$ in a plane embedding of $J/e-b$.
Now it suffices to find a Kuratowski minor in the graph obtained from $W$ by 
splitting $w$ back into $u$ and $v$.

\vspace{0.5\baselineskip}
\noindent\textsc{Case 1}: $u$ and $v$ have a common neighbor.

Let $w_1$ be the common neighbor of $u$ and $v$. Because $w_1$ is also a branch vertex
of $H_b$ and all edges of $C$ are covered by $E(H_a) \cup E(H_b)$, it follows that
all edges of the triangle $G[\{u,v,w_1\}]$ are in $H_b$; so $H_b$ is a subdivision of $K_5$.
In this case, $u$ and $v$ must have degree exactly four in $G[C \cup \{a\}]$.  
Suppose $N(u)=\{v,w_1,w_2,w_3\}$.
Note that $u$ and its four neighbors,
are the five branch vertices of $H_b$. 

Suppose that $\{w_1,w_2,w_3\} \subset N_J(v)$.  This implies $N_J(u)-\{v\} = N_J(v)-\{u\} = \{w_1,w_2,w_3\}$.
Splitting $w$ back into $u$ and $v$ in $W$ produces a $K_5$ subdivision in $G[C \cup \{a\}]$ avoiding $y$ (see \autoref{W4toK5}.

\begin{figure}[htb]
	\centering
	\includegraphics[page=24]{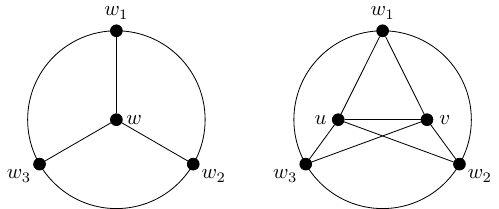}%
	\caption{The faces of the plane embedding of $J/e-b$ containing $w$ on left.  On the right, a subdivision of $K_{5}$ is found in $J-\{b,y\}$ if $\{w_1,w_2,w_3\} \subset N(u) \cap N(v)$.}
	\label{W4toK5}
\end{figure}

If $\{w_2,w_3\} \not\subset N(v)$, then $v=a$ (because $v\neq a$ implies that all of $v$'s neighbors are branch vertices).
Without loss of generality, $w_3$ is not a neighbor of $v$.

\begin{figure}[htb]
	\centering
	\includegraphics[page=25]{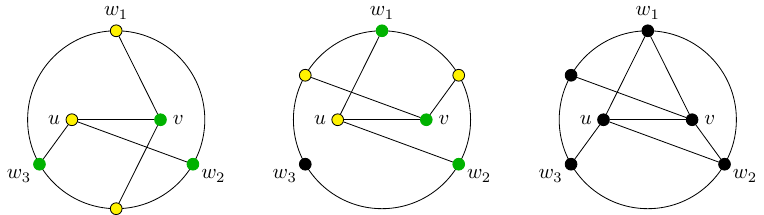}%
	\caption{Splitting $w$ back into $u$ and $v$ produces a Kuratowski minor in $G[C \cup \{a\}]$ avoiding $y$. Cases \ref{fig25.1}, \ref{fig25.2}, and \ref{fig25.3} in this order from left to right.}
	\label{W4toKuratowski}
\end{figure}

The vertices $w_1$, $w_2$, and $w_3$ partition the cycle $W-w$ into into three segments.
Splitting $w$ back into $u$ and $v$ must produce a crossing since $G[C \cup \{a\}]$ is non-planar.
In particular, splitting $w$ back into $u$ and $v$ produces a Kuratowski subdivision or minor in $G[C \cup \{a\}]$ avoiding $y$ 
in each of the remaining cases:
\begin{enumerate}[nosep]
	\item\label{fig25.1} $v$ has a neighbor in the segment between $w_2$ and $w_3$ ($K_{3,3}$ subdivision),
	\item\label{fig25.2} $v$ has neighbors in the segments between $w_1$ and $w_3$ and between $w_2$ and $w_3$ ($K_{3,3}$ subdivision), and
	\item\label{fig25.3} $v$ is adjacent to $w_1$ and $w_2$ so has a neighbor between $w_1$ and $w_3$ ($K_5$ minor).
\end{enumerate}
See \autoref{W4toKuratowski} which gives drawings for cases 1), 2), and 3) in this order from left to right.

\vspace{0.5\baselineskip}
\noindent\textsc{Case 2}: $u$ and $v$ have no common neighbor.

Suppose first that $v \neq a$.  Under this supposition, all of the neighbors of $v$ are also branch vertices of $H_b$ (Lemma~\ref{Proposition3}).
Also the degree of $u$ and $v$ are at least three, so $u$ and $v$ have four different neighbors implying 
that $H_b$ is a subdivision of $K_{3,3}$.  It follows that 
$u$ and $v$ both have degree three in $G[C \cup  \{a\}]$ and have no common neighbors.
Set $N(u) - \{v\} = \{\alpha,\beta\}$ and $N(v) - \{u\} = \{\gamma,\delta\}$.
In $W$ the neighbors of $w$ (the four vertices $\alpha,\beta,\gamma,\delta$) must appear along the face created
by deleting $w$ so that the $\alpha$ and $\beta$ are not consecutive (see left of \autoref{CaseK33MissingY}),
otherwise the plane embedding
of $J/e-b$ could be extended to $J-b$, contradicting that $G[C \cup \{a\}]$ is non-planar.
This produces a subdivision of $K_{3,3}$ 
that does not contain $y$ (see right of \autoref{CaseK33MissingY}). 

\begin{figure}[htb]
	\centering
	\includegraphics[page=26]{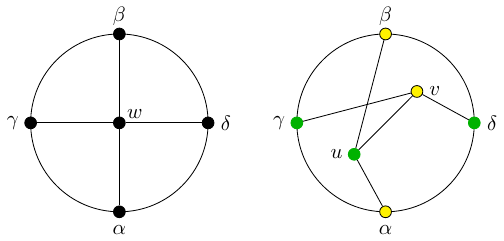}%
	\caption{The faces of the plane embedding of $J/e-b$ containing $w$ on left.  On the right, a subdivision of $K_{3,3}$ is found in $J-\{b,y\}$.}
	\label{CaseK33MissingY}
\end{figure}

So $v = a$.  We may assume that, except for $v$, all of the neighbors of $u$ are on $F$
(otherwise we could substitute one of them for $v$ instead).  

If $C_1$ contains three vertices, say $C_1=\{u,v,z\}$ then, because $u \neq a \neq z$, the closed neighborhoods
of $u$ and $z$ are all branch vertices of $H_b$ (Lemma~\ref{Proposition3}).  Because all of the neighbors of $u$ are on $F$,
$u$ and $z$ are not adjacent so 
$H_b$ is a subdivision of $K_{3,3}$.  All edges of $C$ are covered by $E(H_a) \cup E(H_b)$, so
$u$ and $z$ must have exactly three neighbors.
Let $w_1$ and $w_2$ be the two neighbors of $u$ on $F$.  This implies $N_J(u) = \{v,w_1,w_2\} = N_J(z)$.
Let $\alpha$ be a neighbor of $v$ on $F$.  A subdivision of $K_{3,3}$ exists in
$G[F \cup \{u,v,z\}]$ (with branch vertices $u,v,z,w_1,w_2,\alpha$).  Because this subdivision of $K_{3,3}$ avoids $y$, 
we may assume that $C_1=\{u,v\}$.

If $H_b$ is a subdivision of $K_{5}$, then $u$ has three neighbors $w_1,w_2$,and $w_3$ on $F$.
Because $u$ and $v$ have no common neighbors and the neighbors of $v$ must be on $F$,
at least two of the neighbors of $v$
must be in $F$ in different interior faces, otherwise this planar embedding of $C$ could be
extended to a planar embedding of $G[C \cup \{a\}]$.  This implies a $K_{3,3}$ subdivision avoiding $y$, a contradiction (see left and center panel of \autoref{W4toKuratowski}).

If $H_b$ is a subdivision of $K_{3,3}$, then $u$ and $v$ have degree three but no common neighbors in $G[C \cup  \{a\}]$.
In this case $u$ has two neighbors $w_1,w_2$ on $F$.
The two neighbors of $v$ on $F$
must appear in different interior faces, otherwise this planar embedding of $C$ could be
extended to a planar embedding of $G[C \cup \{a\}]$.  This implies a $K_{3,3}$ subdivision avoiding $y$, a contradiction (see center panel of \autoref{W4toKuratowski}).

\vskip 0.35cm
\noindent
\fbox{\textsc{Step} 3: $t\leq 2$.}

We shall prove Step 3 by contradiction: assume that $t>2$.
Step 2 implies $C_2$ and $C_3$ each contain exactly one vertex.  These two vertices are nonadjacent branch vertex of $H_b$.
The neighbors $C_2$ and $C_3$ on $F$ must be branch vertices of $H_b$ too. 
So $H_b$ is necessarily a subdivision of $K_{3,3}$;
$C_1$ and $C_2$ have the same three neighbors in $S$.
Consequently, the graph $G[F \cup C_2 \cup C_3]$ contains a subdivision of $K_5-e$ 
that has an embedding
in the plane (inherited from the planar embedding of $C$) in which $F$ contains the 
degree-four branch vertices of the subdivided $K_5-e$ and the other branch vertices are in $int(F)$.  This is a contradiction
because
$K_5-e$ has a unique planar embedding in which the triangle connecting the degree-four vertices is a Jordan curve separating the degree-three vertices. 
So this proves Step 3 and we may assume that $t \leq 2$.

\vskip 0.35cm
\noindent
\fbox{\textsc{Step} 4: $t\leq 1$.}

For the rest of the proof we may assume, by Step 2, that $C_1=\{a\}$.  Assume, to the contrary, that $t=2$.
Let $w$ denote the one vertex in $C_2$.
Note that all of the vertices in the closed neighborhood of $w$
are branch vertices of $H_b$ (by Lemma~\ref{Proposition3}) and, $d_G(w)=3$ or $d_G(w) = 4$ 
(because $E(C) \subset E(H_a) \cup E(H_b)$ by Lemma~\ref{LastCaseStructure} part~(iii)).

Recall that $K=V(H_a) - b$.
An edge in $G[K]$ that is not used
by the subdivision $H_a$ is called a {\em chord}.  To the given fixed, plane embedding of $C$
add the vertex $a$ (and its incident edges), 
placing $a$ inside $F$ in general position so as to minimize the number of crossings.
The next definition of crossing chords and related arguments are now with respect to this embedding of $G[C \cup \{a\}]$.
A chord $e=uv$ is {\em crossing} if it is a chord in the face $F$ (that produces a crossing involving $a$ in the plane embedding
of $H_a-b$ inherited from the plane embedding of $C$)
and $a$ has neighbors in different components of $F -\{u,v\}$.
Recall the strong Tutte-Hannani theorem: a graph is planar if and only if it has an embedding so that no
two vertex-disjoint edges cross an odd number of times.  Applying the Tutte-Hannani to the fixed embedding of $G[C \cup \{a\}]$, 
it follows
that there must be two disjoint edges that cross exactly once.

Lemma~\ref{MinDegree3} yields $d_G(w)\geq 3$.  
Let $w_1,w_2$ be the two neighbors of $w$ in $S$
chosen consecutively in the circular ordering of the neighbors of $w$
determined by the plane embedding of $C$ so that $a$ has been placed in the face, $F_1$, using the edges $w_1w$ and $ww_2$.
By the crossing-minimizing placement of $a$, there must be a neighbor of $a$ in $F_1$, call it $a_1$
(possibly $a_1 \in \{w_1,w_2\}$).
Because there is a crossing, there must be a neighbor of $a$ outside of $F_1$, call it $a_2$. 
We may choose $a_2$ to be the first neighbor of $a$ clockwise from $w_2$ (away from $w_1$) along $F$ that is outside of $F_1$.
Let $w_3$ be the next neighbor (clockwise) of $w$ that is at least as far as $a_2$.
If $a$ has a neighbor in $F_1 - \{w_1,w_2\}$, then delete all edges incident to $w$ except $w_1w$ and $ww_2$
and then contract these two edges to produce a crossing chord in this minor of $G$.
Otherwise we may assume all neighbors of $a$ in $F_1$ are in $\{w_1,w_2\}$, in particular, we may assume $a_1=w_1$.
If $a_2\neq w_3$, then contracting $w_1w$ and $ww_3$ produces a crossing chord in a minor of $G$.
So assume $\{a_1,a_2\} = \{w_1,w_3\}$. 
If $N_J(w) = \{w_1,w_2,w_3\}$ and $N_J(a) = \{w_1,w_3\}$, then a crossing-free drawing of $G[C \cup \{a\}]$ could
be made by placing $a$ in the face defined by the edges $w_1w$ and $ww_3$.
If $w$ has another neighbor, $w_4$, then contracting $w_2w$ and $ww_4$ produces a crossing chord in a minor of $G$.
Otherwise $N_J(w) = \{w_1,w_2,w_3\}$, and all of $a$'s neighbors are in the face determined by $w_1w$ and $ww_3$.

So the only case in which we cannot create a minor with a crossing chord is when
$a$ and $w$ are both branch of $H_b$, which is a subdivided $K_{3,3}$, and
$N(a) \cap F = \{w_1,w_2,w_3\} = N(w) \cap F$.  

If $H_a$ is a subdivided $K_5$, then $F$ has three
branch vertices of $H_a$ which partition $F$ into three paths.  If $w_i$, $w_j$ with $1 \leq i < j \leq 3$
are on the same path of $H_a$, then rerouting this path through $w_i$, $w_j$ and $w$ 
produces another $H_a$ that has smaller $|int(F)|$.  So the only remaining case is if
$w_1$, $w_2$, and $w_3$ are interior vertices on the three different paths of $H_a$ on $F$.
This is shown on the left panel of \autoref{HardK5Case}.

\begin{figure}[htb]
	\centering
	\includegraphics[page=27]{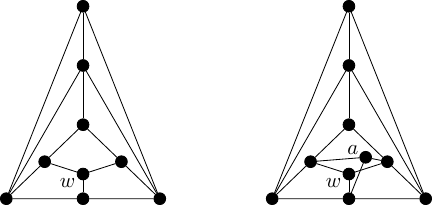}%
	\caption{$H_a$ is a subdivided $K_5$ and the neighbors of $w$ are interior vertices on the three paths of $H_a$ on $F$ (left). Adding $a$, produces graph shown in center panel which remains non-planar after deleting $y$ (top vertex).}
	\label{HardK5Case}
\end{figure}

Adding $a$ (and its incident edges to the neighbors of $w$) produces graph (shown in right panel of \autoref{HardK5Case}) 
which remains non-planar after deleting $y$, a contradiction.

Consider now the case in which $H_a$ is a subdivided $K_{3,3}$. 
If $w_i$, $w_j$ with $1 \leq i < j \leq 3$
are on the same path of $H_a$, then rerouting this path through $w_i$, $w_j$ and $w$ 
produces another $H_a$  with smaller $|int(F)|$.  So the only remaining case is if
$w_1$, $w_2$, and $w_3$ are interior vertices on the three different paths of $H_a$ on $F$.
The non-isomorphic cases are shown on the left panel of \autoref{HardK33Case}.

\begin{figure}[htb]
	\centering
	\includegraphics[page=28]{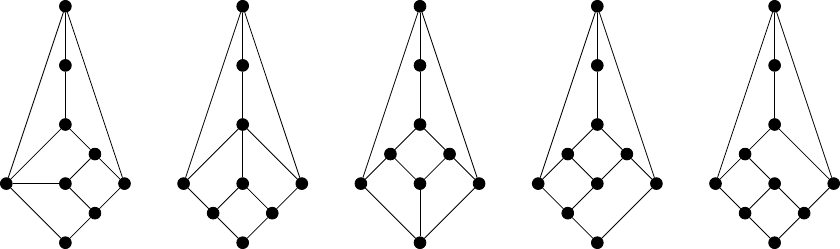}%
	\caption{$H_a$ is a subdivided $K_{3,3}$ and the neighbors of $w$ are interior vertices on the three different paths of $H_a$ on $F$. The non-isomorphic cases are shown here.}
	\label{HardK33Case}
\end{figure}

In each case it is possible to find a new subdivision of $K_{3,3}$ for $H_a$ (in which $x$ and $y$ are branch vertices connected by a path through $b$) with smaller $|int(F)|$ as shown in 
\autoref{HardK33CaseResolved}.

\begin{figure}[htb]
	\centering
	\includegraphics[page=29]{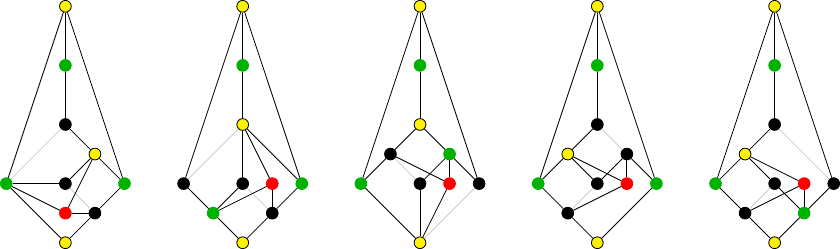}%
	\caption{In each case a new subdivision of $K_{3,3}$ can be chosen for $H_a$ with smaller $|int(F)|$. The branch vertices of the new subdivision of $K_{3,3}$ are shown in green and yellow.  Vertex $a$ is shown in red. The non-branch vertex $b$ of this $K_{3,3}$ subdivision, which connects $x$ and $y$, is not shown.}
	\label{HardK33CaseResolved}
\end{figure}


  If $a$ has neighbors in $C$ that are all in a single path
of $H_a-b$, then a technical difficulty arises: 
$\overline{a}$ indeed hits only one face (as Step 1 guarantees), but there are two choices for this one face --- 
either of the two faces incident to this path (in the unique plane
embedding of $H_a-b$). The next step excludes this possibility, resolving this difficulty.

A path in $C$ whose internal vertices are disjoint from $H_a-b$ is called an {\em external} path.

\vskip 0.35cm
\noindent
\fbox{\textsc{Step} 5: The neighbors of $a$ in $C$ are not all on one path of $H_a - b$.}

Assume, to the contrary, that the neighbors of $a$ in $C$ are all on one path of $H_a - b$.
Because $t=1$, every crossing pair of edges involves a crossing chord and an edge incident to $a$.
Keep in mind that, by the minimality of $H_a$, the endpoints of any crossing chord cannot both be on
the same path connecting branch vertices of $H_a$.

Let $F_1$ and $F_2$ be the two faces (in the unique plane
embedding of $H_a-b$)
incident to the path $P$ (connecting branch vertices $H_a - b$) that contains all of the neighbors of $a$ in $C$.
It is very important to note that Steps 2-4 apply to both $F_1$ and $F_2$. 
In particular, this means that no vertices of $C$ are in the interior
of $F_1$ or $F_2$; so these faces of the plane embedding of $H_a-b$ are actually faces of the plane embedding of $C$.
As in prior steps, we now consider cases depending on whether $H_a - b$ is a subdivision of $K_5-e$ or
$K_{3,3}-e$.  

First consider the easier case in which $H_a - b$ is a subdivision of $K_5-e$.
Suppose that $P$ is a path ending at $x$.   Let $a_1$ and $a_2$ be the extreme neighbors that $a$ has
on $P$; that is, $a_1$ is the furthest from $x$ and $a_2$ is the closest to $x$. 
There must be a chord in $F_1$ that produces a crossing with
$a$.  This chord must have one endpoint between $a_1$ and $a_2$.  The other endpoint cannot  be on $P$ by the minimality of $H_a$.
There must be a chord in $F_2$ with similar properties.  These chords are drawn in color in \autoref{K5Case3}.
Regardless of where the other endpoints of these chords occur, 
contractions produce the graph on the right of \autoref{K5Case3}, a $K_5$ minor that avoids $y$, a contradiction.

\begin{figure}[htb]
	\centering
	\includegraphics[page=30]{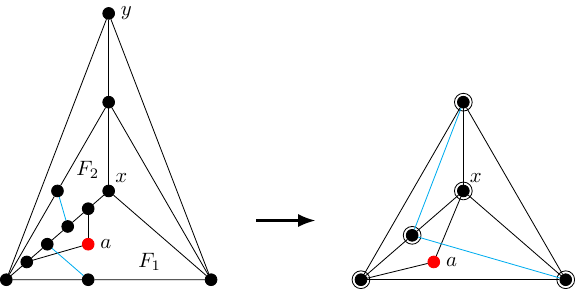}%
	\caption{$H_a-b$ is a subdivision of $K_5-e$ and $a$ has neighbors in $C$ only along a path of $H_a-b$ ending at $x$.  A $K_5$ minor avoiding $y$ appears.}
	\label{K5Case3}
\end{figure}

So suppose that $P$ is a path of $H_a-b$ that does not end at $x$. 
Again $F_1$ and $F_2$ have chords that produce crossing with $a$ (see left of \autoref{K5Case4}).
Regardless of where the other endpoints of these chords occur, contractions produce
the graph shown in the center of \autoref{K5Case4}.  Adding the vertex $b$ and 
the edge $ab$ (which corresponds to contracting the light component) produces the apex obstruction $M$ as a minor of $G$,
a contradiction.

\begin{figure}[htb]
	\centering
	\includegraphics[page=31]{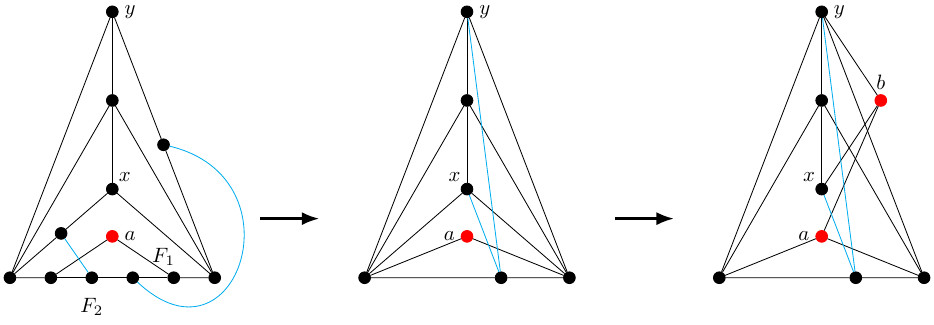}%
	\caption{$H_a-b$ is a subdivision of $K_5-e$ and $a$ has neighbors in $C$ only along a path of $H_a-b$ that does not end at $x$. An $M$ minor appears after adding $b$ and the edge $ab$ (which corresponds to contracting the light component).}
	\label{K5Case4}
\end{figure}

Next consider the case in which $H_a - b$ is a subdivision of $K_{3,3}-e$.
Two subcases occur depending on whether the path $P$ (containing all neighbors of $a$ in $C$) is incident
to $x$ or $y$, or neither $x$ nor $y$.  

First consider the subcase in which $P$ is incident to neither $x$ nor $y$.
Again define $a_1$ and $a_2$ to be the extreme neighbors of $a$ along $P$; that is, the neighbors of $a$ closest
to the branch vertices of $H_a-b$ at the ends of $P$.
As before, $F_1$ and $F_2$ must contain chords producing crossings with edges incident to $a$  (see \autoref{K33CasesBegin})
These chords must have ends in the interval of $P$ between $a_1$ and $a_2$; call these ends $\alpha$ and $\beta$.  The vertices at the other ends
of these chords are of three possible types, 

\begin{itemize}
\item[type (I):] can be contracted to $x$ or $y$ along paths in $H_a-b$ (without passing through a branch 
vertex of $H_a-b$), 
\item[type (II):] appear in the interior of the path in $H_a-b$ opposite $x$ or $y$ (along the face $F_1$ or $F_2$),
\item[type (III):] branch 
vertex of $H_a-b$ (adjacent to $x$ or $y$ in $H_a-b$).
\end{itemize}

\begin{figure}[H]
	\centering
	\includegraphics[page=32]{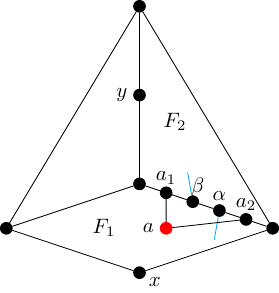}%
	\caption{The case in which $H_a-b$ is a subdivision of $K_{3,3}-e$ and $a$ has neighbors in $C$ only along a path of $H_a-b$ that ends neither at $x$ nor $y$.}
	\label{K33CasesBegin}
\end{figure}

First we argue that $\alpha = \beta$.  If $\alpha \neq \beta$,
then we claim that the edges along $P$ can be contracted, preserving $H_a$, so that a Kuratowski avoiding $b$ exists.
This is a contradiction because then neither $a$ nor $b$ is an apex in $G[C \cup \{a,b\}]$ after contracting $\alpha$ to $\beta$.
All nine cases (two end vertices each have independently any of three types) can be dismissed by examining graphs shown in \autoref{K33CasesEndsAgree}.  All graphs in this figure have $a_1$ and $a_2$ contracted (if necessary) to the ends of the path $P$ and $\alpha$ contracted to $\beta$.
If both ends of the crossing chords are type (I), the leftmost graph shows that there is a minor of $K_{3,3}$ avoiding $b$.
If the ends of the crossing chords are both type (III), the right graph shows that there is a minor of $K_{5}$ avoiding $b$.
The right-most graph also works if one chord is of type (I) and the other is of type (III).
If one end of a crossing chords is type (I) and other is type (II), the left-center graph shows that there is a minor of $K_{3,3}$ avoiding $b$.  If one end is of type (II) and other is type (II) or (III), then the right-center graph of \autoref{K33CasesEndsAgree} shows
a minor of $K_{3,3}$ avoiding $b$.
 
\begin{figure}[H]
	\centering
	\includegraphics[page=33]{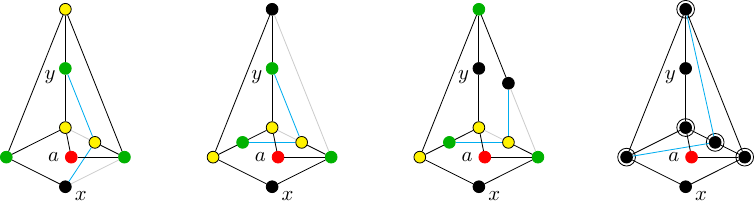}%
	\caption{If $\alpha \neq \beta$, then edges along $P$ can be contracted, preserving $H_a$, so that a Kuratowski subgraph avoiding $b$ emerges.}
	\label{K33CasesEndsAgree}
\end{figure}

So we may assume that $\alpha = \beta$.


If both chords have ends that are of type (I), then an edge $e$ exists (shown dotted in \autoref{K33Case1Final})
such that $J-e-a$ and $J-e-b$ are non-planar, contradicting Lemma~\ref{BasicLemma} part~\ref{BasicLemma.delete}.
This can be seen in \autoref{K33Case1Final} where a subdivision of $K_{3,3}$ is found in both $J-e-a$ (new $H_a$), and $J-e-b$ (new $H_b$).

\begin{figure}[H]
	\centering
	\includegraphics[page=34]{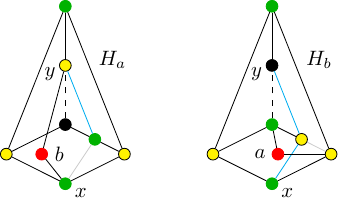}%
	\caption{If $\alpha = \beta$ and both crossing chords have ends that are of type (I),
	then a choice of $H_a$ and $H_b$ are shown that avoid an edge of $C$ (shown dotted).}
	\label{K33Case1Final}
\end{figure}


If both chords have ends that are of type (III), then an edge $e$ exists (shown in pink in \autoref{K33Case2Final})
such that $J/e-a$ and $J/e-b$ are non-planar, contradicting Lemma~\ref{BasicLemma} part~\ref{BasicLemma.contract}.  This figure
also shows a contradiction if one end of a crossing chord has type (III).

\begin{figure}[H]
	\centering
	\includegraphics[page=35]{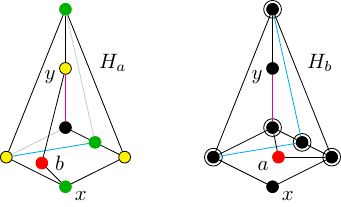}%
	\caption{If $\alpha = \beta$ and one crossing chord has type (III), then a choice of $H_a$ and $H_b$ are shown that a share a contractible edge (shown in pink).}
	\label{K33Case2Final}
\end{figure}


If the end of one crossing chord has type (I) and the other is of type (II), 
then an edge $e$ exists (shown dotted in \autoref{K33Case3Final})
such that $J-e-a$ and $J-e-b$ are non-planar, contradicting Lemma~\ref{BasicLemma} part~\ref{BasicLemma.delete}.

\begin{figure}[H]
	\centering
	\includegraphics[page=36]{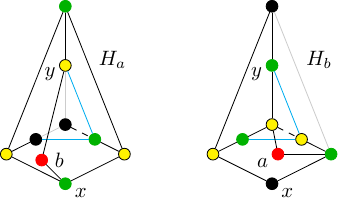}%
	\caption{If $\alpha = \beta$ and one crossing chord has type (I) and the other has type (II),
	then a choice of $H_a$ and $H_b$ are shown that a share a deletable edge (shown dotted).}
	\label{K33Case3Final}
\end{figure}


If both ends of the crossing chords have type (II), 
then an edge $e$ exists (shown in pink in \autoref{K33Case4Final})
such that $J/e-a$ and $J/e-b$ are non-planar, contradicting Lemma~\ref{BasicLemma} part~\ref{BasicLemma.delete}. 

\begin{figure}[H]
	\centering
	\includegraphics[page=37]{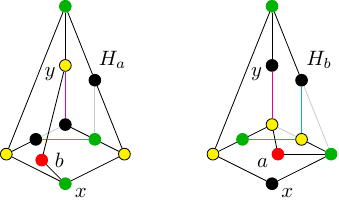}%
	\caption{If $\alpha = \beta$ and both ends of the crossing chords have type (II), then a choice of $H_a$ and $H_b$ are shown that a share a contractible edge (shown in pink).}
	\label{K33Case4Final}
\end{figure}

Finally consider the subcase in which $a$ has neighbors only along the path $P$ of $H_a-b$ and $P$ 
is incident to $x$ or $y$.  Without loss of generality,
the path $P$ is incident to $x$.  There are two possibilities, shown in \autoref{K33LastCases}, according to whether
$a$ has neighbors on one side of $x$ (along $P$) or both sides of $x$.  \autoref{K33LastCases} introduces the labeling of the remaining branch vertices of $H_a$ as
$p,q,r,s$ as well as the extreme neighbors, $a_1$ and $a_2$, of $a$ along $P$, as shown.

\begin{figure}[H]
	\centering
	\includegraphics[page=38]{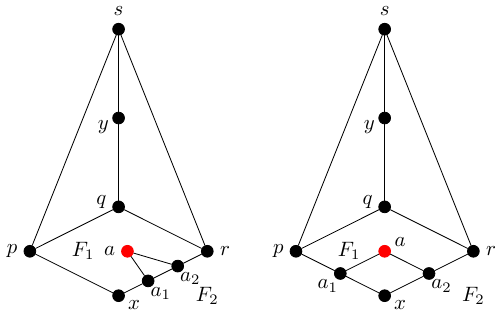}%
	\caption{The two possibilities when $H_a-b$ is a subdivision of $K_{3,3}-e$, the path $P$ is incident to $x$, and $a$ has neighbors on one side of $x$ (left) or both sides of $x$ (right).}
	\label{K33LastCases}
\end{figure}

Recall that there must be a crossing chord in both faces $F_1$ and $F_2$ preventing $a$ from being placed
into each face.  By the minimality of $H_a$, these chords must have ends on distinct paths of $H_a$.

Regarding the leftmost graph in \autoref{K33LastCases}, we may assume that
$a_1=x$ and $a_2=r$, since contracting edges along the $[x,a_1]$ and $[a_2,r]$ segments of $P$
cannot destroy crossing chords.  The crossing chords must have one endpoint between $a_1$ and $a_2$ along $P$
and another endpoint outside this interval.  For the $F_1$ face, the other endpoint is contractible to $p$ or to $q$.
 For the $F_2$ face, the other endpoint is contractible to $p$ or to $s$.
The four resulting cases are depicted in \autoref{K33LastCasesTABLE}.

\begin{figure}[H]
	\centering
	\includegraphics[page=39]{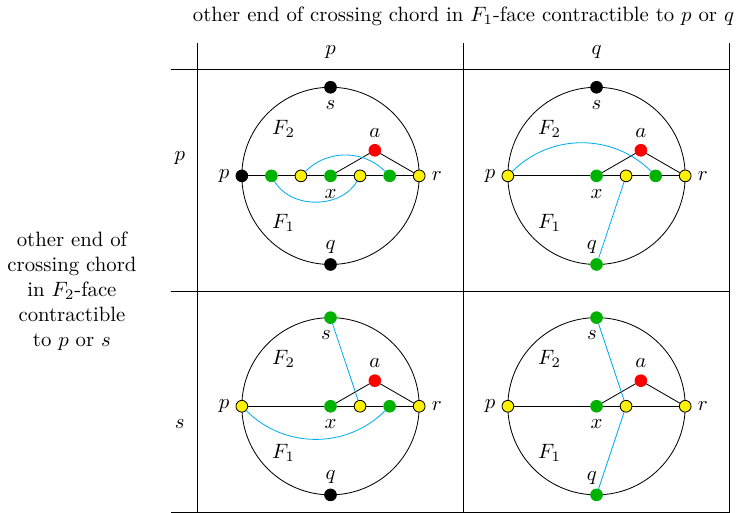}%
	\caption{The four possible combinations for ends of crossing chords (shown in light blue).}
	\label{K33LastCasesTABLE}
\end{figure}

In all cases, except the bottom right case of \autoref{K33LastCasesTABLE}, the ends of the crossing chords
must be distinct because the chords must cross each other as well cross edges incident to $a$.  If the chords did not cross each other,
then they could be brought into a single face.  In the bottom right case, we may contract the ends of the crossing chords in the $(x,r)$-interval of $P$ to the same vertex.
As the figure shows, all four cases result in a $K_{3,3}$ subdivision avoiding $y$, a contradiction.

In the final case of Step 5, consider the rightmost graph in \autoref{K33LastCases}.
This graph can be redrawn as shown in \autoref{FinalCaseStep5}; the ends of the crossing chords are left undecided in this drawing.

\begin{figure}[H]
	\centering
	\includegraphics[page=40]{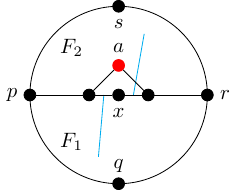}%
	\caption{A redrawing of a portion of the rightmost graph from \autoref{K33LastCases} with crossing chords shown in light blue.}
	\label{FinalCaseStep5}
\end{figure}

If the ends of the chords can be contracted to $s$ and $q$, then a minor isomorphic to $M$ (a Petersen-family graph) in $G$ emerges
(\autoref{FinalCaseStep5-Memerges}) after adding $b$ and the edge $ab$, the latter of which corresponds to contracting the light component of $G$.

\begin{figure}[H]
	\centering
	\includegraphics[page=41]{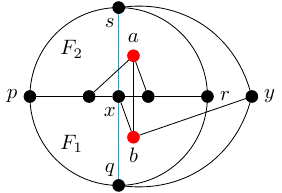}%
	\caption{If the crossing chords have ends that can be contracted to $s$ and $q$, then a minor of $M$ emerges.}
	\label{FinalCaseStep5-Memerges}
\end{figure}

In the remaining possible placements of the ends of the crossing chords, each possibility can be contracted to one of the two shown 
in \autoref{FinalCaseStep5-Last}.  Observe that the crossing chords must again cross each other (as well as edges incident to $a$)
to avoid being placed into a single face.  In both cases, a $K_{3,3}$ avoiding $y$ is shown, a contradiction.

\begin{figure}[H]
	\centering
	\includegraphics[page=42]{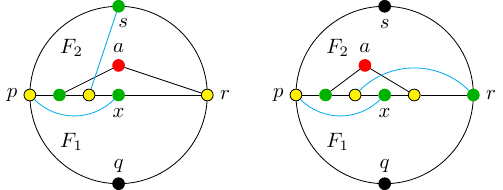}%
	\caption{Crossing chords must cross each other in the final two possible placements of the ends of these crossing chords (shown in blue). The cases can each be contracted to the graphs depicted, or similar graphs. A $K_{3,3}$ avoiding $y$ arises in both cases.}
	\label{FinalCaseStep5-Last}
\end{figure}


\vskip 0.35cm
\noindent
\fbox{\textsc{Step} 6: The final contradiction to complete the proof.}

In this final step, we may assume that $a$ has neighbors in only one face of the plane embedding
of $H_a-b$ (Step 1).  Furthermore there are no other vertices of $C$ in this face (Steps 2-4).
Finally we may assume that all of $a$'s neighbors in $C$ are not along a single path of $H_a-b$ (Step 5).
Consequently, there is exactly one face, $F$, containing all of the neighbors of $a$ in $C$ and there must be a crossing chord in $F$.

If $H_a$ is a subdivision of $K_5$, then there are four remaining non-isomorphic positions
for a crossing chord.  These cases are shown in \autoref{ReduceK5}.

\begin{figure}[H]
	\centering
	\includegraphics[page=43]{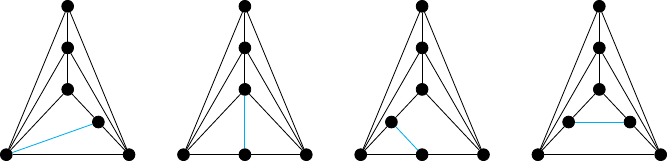}%
	\caption{$H_a$ is a subdivision of $K_5$. The four non-isomorphic positions for the crossing chord (shown in light blue) when endpoints are not along the same path of $H_a$.}
	\label{ReduceK5}
\end{figure}

In each of these cases (see \autoref{ReduceK5}), a new $K_{3,3}$ subdivision for $H_a$ exists with $x$ and $y$ branch vertices connected by $b$ (see \autoref{ReduceK5Resolved}). One can argue in each case that this new $H_a$ has a unique embedding of $H_a-b$ such that $a$ now hits multiple faces, contradicting Step~1. This argument requires some care. For example, in the left-most graph, clearly $a$ must have neighbors that cross the blue chord. These could simply be a vertex along the $(u,x)$-path and another along the $(u,v)$-path. However all of the neighbors of $a$ could not be along these paths by Step 5.

\begin{figure}[H]
	\centering
	\includegraphics[page=44]{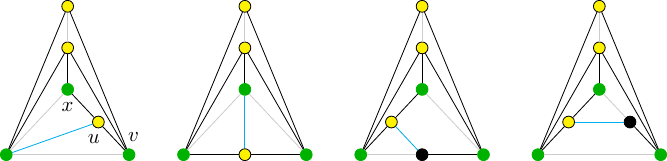}%
	\caption{In each case of \autoref{ReduceK5} a $K_{3,3}$ subdivision for $H_a$ exists with $x$ and $y$ branch vertices connected by $b$.}
	\label{ReduceK5Resolved}
\end{figure}

Now consider the ten non-isomorphic ways a crossing chord can appear in $F$ when 
$H_a$ is a subdivision of $K_{3,3}$.
Six of these cases (see \autoref{K33ReduceCase1}) can be dismissed using an argument similar to
the one given in the prior paragraphs.

\begin{figure}[H]
	\centering
	\includegraphics[page=45]{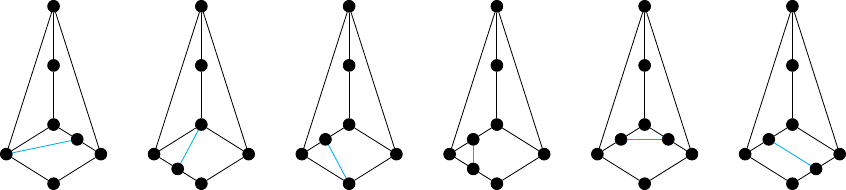}%
	\caption{Six of the ten non-isomorphic ways a crossing chord (shown in light blue) can appear in $F$ when $H_a$ is a subdivision of $K_{3,3}$.}
	\label{K33ReduceCase1}
\end{figure}

Each case of \autoref{K33ReduceCase1} has a new $K_{3,3}$ subdivision with $x$ and $y$ branch
vertices connected by $b$.  This new $H_a-b$ has $a$ that now hits multiple faces (by Step 5), contradicting Step 1
(see \autoref{K33ReduceCase1Resolved});

\begin{figure}[H]
	\centering
	\includegraphics[page=46]{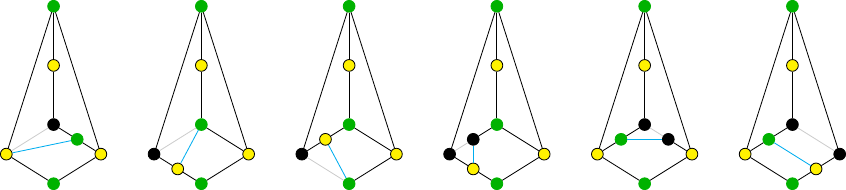}%
	\caption{Each case of \autoref{K33ReduceCase1} has a new $K_{3,3}$ subdivision with $x$ and $y$ branch vertices connected by $b$. This new $H_a-b$ has $a$ that now hits multiple faces, contradicting Step 1.}
	\label{K33ReduceCase1Resolved}
\end{figure}

The last four remaining cases are shown in \autoref{RemainingCases}.

\begin{figure}[H]
	\centering
	\includegraphics[page=47]{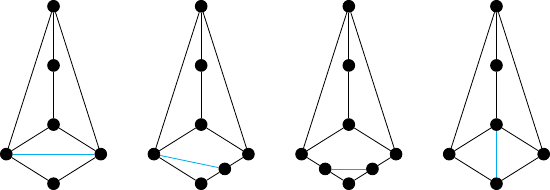}%
	\caption{The four remaining chord placements in $H_a-b$.}
	\label{RemainingCases}
\end{figure}

Consider left-most graph in \autoref{RemainingCases} in which $uv$ is the crossing chord.
This means that $a$ has a neighbor on the upper half (above the chord $uv$) of face $F$ and a neighbor on the lower half
(below the chord $uv$) of face $F$. Observe that
there must be an external path $P$ that connects the lower open interval 
$(v,u)$ of the exterior face (shown in red in \autoref{ONE})
to the upper open interval $(u,v)$ of the exterior face (shown in green in \autoref{ONE}); otherwise
the $uv$ chord could be drawn on the exterior face, reducing the number of crossings
produced when reinserting $a$ into the planar embedding of $C$; that is, the external path $P$ {\em blocks} the chord $uv$.
The resulting graph can be contracted to
the one shown on the right of \autoref{ONE}.  Adding vertex $a$ to this graph produces a $K_5$ minor of $G$ which implies
that $H_a - b = H_b -a$ in the original graph $G$, so $\{x,y\}$ is another $2$-cut of $G$, a contradiction.

\begin{figure}[H]
	\centering
	\includegraphics[page=48]{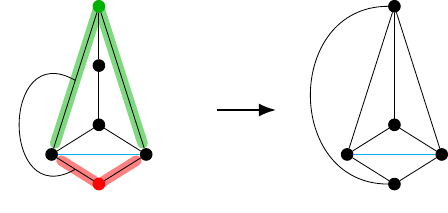}%
	\caption{An external path connecting the lower exterior face (red) to the upper exterior face (green) implies a $K_5$ minor of $G$.}
	\label{ONE}
\end{figure}

Consider next the center-left graph in \autoref{RemainingCases}.  Again let $uv$ be the crossing chord; let
$w$ be the branch vertex of $H_a-b$ opposite $u$ on face $F$ (as shown in \autoref{TWO}).  

\begin{figure}[H]
	\centering
	\includegraphics[page=49]{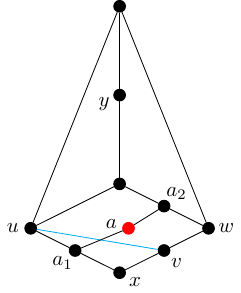}%
	\caption{By Step 5, the vertex $a$ must have neighbors above $uw$ and below the blue crossing chord.}
	\label{TWO}
\end{figure}

The vertex $a$ has neighbors
on $F$ above and below the chord $uv$ because $uv$ is a crossing chord.
By Step 5, $a$ has neighbors between $u$ and $w$ (on the upper part of $F$ --- \autoref{TWO} shows one case); otherwise
all of $a$'s neighbors would occur on the path of $H_a-b$ containing $x$.
Now the same analysis as given in the prior paragraph (\autoref{ONE}) applies.
We omit further details for this case. Similar reasoning applies to the center-right case shown 
in \autoref{RemainingCases}.

The final analysis regards the case in which the crossing chord is a vertical through $F$ connecting
branch vertices of $H_a-b$ (shown as the rightmost graph of \autoref{RemainingCases}).
By Step 5 not all of $a$ neighbors can occur on the path of $H_a-b$ containing $x$.
Thus there are three cases that remain; these are shown in \autoref{VERTICAL}.

\begin{figure}[H]
	\centering
	\includegraphics[page=50]{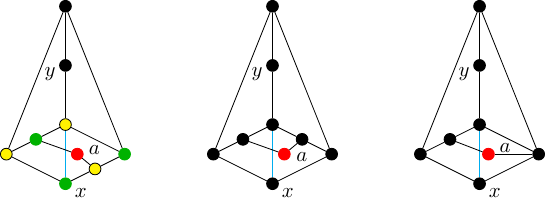}%
	\caption{The left graph contains a subdivision of $K_{3,3}$ avoiding $b$ and $y$. The other two contract to $P_7$ after adding $b$ and the edge $ab$.}
	\label{VERTICAL}
\end{figure}

Observe that the left most graph of \autoref{VERTICAL} contains a subdivision of $K_{3,3}$ avoiding
$b$ and $y$.  The other two graphs contain a minor of $P_7$ (see \autoref{P7}) after adding $b$ and the edge $ab$.
\end{proof}

The next corollary strengthens the statement of Proposition~\ref{ExistenceOfKuratowskis}
to show that vertices in the unique $2$-cut must be branch vertices of all of their Kuratowski witnesses.

\begin{corollary}\label{BranchVertices}
If Assumptions~\ref{assumptions} are satisfied, then $b$ is a branch vertex 
for any Kuratowski subgraph of $G$ avoiding $a$ (and vice versa, $a$ is a branch vertex 
for any Kuratowski subgraph of $G$ avoiding $b$).
\end{corollary}
\begin{proof} 
Note that $b$ must have at least three neighbors in $C$ since Proposition~\ref{ExistenceOfKuratowskis}
guarantees a Kuratowski subgraph avoiding $a$ (with vertices all in $C \cup \{a\}$) in which $b$ must appear as a branch vertex.
Suppose now that $H_a$ is an arbitrary Kuratowski subgraph of $G$ avoiding $a$.
Clearly $V(H_a) \subseteq C \cup \{b\}$.
If one of the edges in $G[C \cup \{b\}]$ incident to $b$ does not appear in $H_a$, then it is not in
$E(H_a) \cup E(H_b)$, contradicting Lemma~\ref{LastCaseStructure} part~(iii).
\end{proof}


\subsection{Two Kuratowski subgraphs have branch vertices that cover \texorpdfstring{$C \cup \{a,b\}$}{C∪\{a,b\}}}

With Corollary~\ref{BranchVertices}, establishing that $a$ and $b$ must be branch vertices 
of all of their Kuratowski witnesses, we are almost ready to prove the main objective of this subsection: the
existence of two Kuratowski subgraphs
whose branch vertices cover $C \cup \{a,b\}$ (Theorem~\ref{WeakWTheorem}).  First 
we focus on properties of vertices in $C$ that are not branch vertices of Kuratowski witnesses.

\begin{lemma}[Non-Branch Vertex Lemma]
\label{NonBranchLemma}
Assume Assumptions~\ref{assumptions} and
$H_a$ and $H_b$ are Kuratowski subgraphs of $G$ avoiding $a$ and $b$, respectively.
If $w\in C$ is a branch vertex of neither $H_a$ nor $H_b$:
\begin{itemize}
\item[(i)] $d_G(w) = 4$, with $w$ incident to two edges in $E(H_a)$ and two in $E(H_b)$.
\item[(ii)] If $e=wx \in E(H_a)$ then $x \in V(H_b) \cup \{b\}$, $H_b \not\subset G/e$, and
   $b$ is the only apex for $G/e$.
\item[(iii)] If $e=wx \in E(H_b)$ then $x \in V(H_a) \cup \{a\}$, $H_a \not\subset G/e$, and
   $a$ is the only apex for $G/e$.
\end{itemize}
\end{lemma}
\begin{proof} (i) Assume that $c \in C$ is not a branch vertex of $H_a$ or $H_a$.  So the degree of $c$ in $H_a$ and $H_b$ is two.
Lemma~\ref{MinDegree3} implies $d_G(w) \geq 3$.  If $d_G(w)>4$, then there is an edge of $G[C \cup \{a,b\}]$ that is not
covered by $H_a$ or $H_b$ contradicting Lemma~\ref{LastCaseStructure} part~(iii).  To prove claim (i), it suffices to prove
that $d_G(w)\neq 3$.  Assume, to the contrary, that $d_G(w)=3$. The pigeon-hole principle guarantees
an edge $wx \in E(H_a) \cap E(H_b)$ and so $x \not\in \{a,b\}$.
Note that $G/wx$ must have an apex $z$ in $H_a/wx \cap H_b/wx$, so $z \not\in \{a,b\}$.
However, $z$ must also separate $a$ from $b$ in $G[C \cup \{a,b\}]$ since otherwise
$L^+$ would still be a minor of $G/wx$.  But then either $\{a,z\}$ or $\{b,z\}$ is another $2$-cut
of $G$ (contradicting that $S$ is the only $2$-cut) or $z$ is the vertex resulting from the contraction of the edge $wx$,
contradicting Lemma~\ref{LastCaseStructure} part~(iv).

(ii-iii) By symmetry, it suffices to prove (ii).  Consider $e=wx \in E(H_a)$.  
If $x \not\in V(H_b) \cup \{b\}$, then $H_a$ and $H_b$ remain Kuratowski
subgraphs in $G/e$ implying that any apex $z$ for $G/e$ must be in $H_a \cap H_b \subset C$.  However, 
this means that
$z$ must also separate $a$ from $b$ in $G[C \cup \{a,b\}]$ since otherwise
$L^+$ would still exist in $G/e$.  But then either $\{a,z\}$ or $\{b,z\}$ is another $2$-cut
of $G$ (contradicting that $S$ is the only $2$-cut) or $z$ is the vertex resulting from the contraction of the edge $wx$,
contradicting Lemma~\ref{LastCaseStructure} part~(iv).  Therefore $x \in V(H_b) \cup \{b\}$.
The reader can check that similar reasoning applies if $H_b$ remains a Kuratowski subgraph of $G/e$
or if $b$ is not the only apex for $G/e$.
\end{proof}

\begin{theorem} 
\label{TK5withWVertices}
Assume Assumptions~\ref{assumptions}.
Choose any Kuratowski subgraphs $H_a$ and $H_b$ avoiding $a$ and $b$, respectively,
that also minimize $|E(H_a)| + |E(H_b)|$.
If $H_a$ is a subdivision of $K_5$ or $H_b$ is a subdivision of $K_5$, then any vertex in $C$ is a branch vertex of $H_a$ or a branch vertex of $H_b$.
\end{theorem}
\begin{proof} Without loss of generality, $H_a$ is a subdivision of $K_5$. Assume, to the contrary,
there is a vertex $w \in C$ that is a not a branch vertex either $H_a$ or $H_b$.
By Theorem~\ref{NonBranchLemma} part~(i), $deg_G(w)=4$ and $w \in V(H_a) \cap V(H_b)$.
Now $w$ has two neighbors in $H_b$, at least one of which is not vertex $a$.
Consider a neighbor $x$ of $w$ such that $wx \in E(H_b)$ and $x \neq a$.  
By Theorem~\ref{NonBranchLemma} part~(iii),  $x \in V(H_a)$.
There are three cases shown in \autoref{NonBranchTK5}.  

\begin{figure}[H]
	\centering
	\includegraphics[page=51]{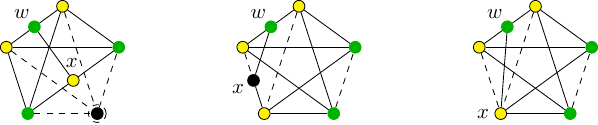}%
	\caption{Deleting dotted lines/vertex and adding the edge $wx$ reveals a subdivision of $K_{3,3}$ 
	with fewer edges than the original subdivision of $K_5$.}
	\label{NonBranchTK5}
\end{figure}

The figure shows
a subdivision of $K_5$ that represents $H_a$; $w$ is a non-branch vertex with
a neighbor $x$ such that $wx \in E(H_b) - E(H_a)$.  The vertex $x$ could
be in any of three positions.
The dotted lines indicates paths of $H_a$ that can be deleted
leaving a new subdivision of $K_{3,3}$ with fewer edges than $H_a$, in each case.  This new Kuratowski
subgraph avoids vertex $a$ also, so contradicts the choice of $H_a$.
So these cases essentially follow from a commonly rediscovered fact that a vertex minimal non-planar graph that is not just a subdivision of $K_5$ has a spanning $K_{3,3}$ subdivision.

In each case, a new choice of $H_a$ as a subdivision of $K_{3,3}$ has fewer edges than the current $H_a$;
this contradicts that the original choice of $H_a$ and $H_b$ minimized $|E(H_a)| + |E(H_b)|$.
\end{proof}

The next theorem is the main result of this subsection.

\begin{theorem} 
\label{WeakWTheorem}
Assume Assumptions~\ref{assumptions}.  There are Kuratowski subgraphs $H_a$ and $H_b$ avoiding $a$ and $b$ respectively,
such that 
any vertex in $C \cup \{a,b\}$ is a branch vertex of $H_a$ or a branch vertex of $H_b$.
\end{theorem}
\begin{proof}  Choose Kuratowski subgraphs $H_a$ and $H_b$
as follows:
\begin{itemize}
\item[(i)] $a \not\in V(H_a)$ and $b \not\in V(H_b)$,
\item[(ii)] maintaining (i), minimize $|E(H_a)| + |E(H_b)|$,
\item[(iii)] maintaining (i) and (ii), minimize $|W|$, where
$$W = (C \cup \{a,b\}) - \{ v : v \mbox{ is branch vertex of } H_a \mbox{ or } H_b\}.$$
\end{itemize}
It suffices to prove that this choice produces $W = \varnothing$.
Corollary~\ref{BranchVertices} implies $a,b \notin W$.

If $H_a$ is a subdivision of $K_5$ or $H_b$ is a subdivision of $K_5$, then Theorem~\ref{TK5withWVertices}
yields $W = \varnothing$.
So, we may assume $H_a$ and $H_b$ are subdivisions of $K_{3,3}$.

\begin{figure}[ht]
	\centering
	\includegraphics[page=52]{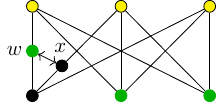}%
	\caption{A depiction of $H_a$ and the extra edge $wx \in E(H_b)$.  Contracting $wx$ preserves $H_b$ and non-planarity of $H_a$.}
	\label{WandXAlongCommon}
\end{figure}

Assume, to the contrary, that $W \neq \varnothing$.
Let $w$ be an arbitrary vertex in $W$.  If possible, choose  $H_a$ and $H_b$, subdivisions of $K_{3,3}$,
satisfying (i) - (iii), minimum $|W|$ and $w \in W$ so that $w \notin N(a)$. 
By Theorem~\ref{NonBranchLemma}, $d_G(w)=4$, $w \in V(H_a) \cap V(H_b)$,
and all of the neighbors of $w$ are also vertices in $(V(H_a) \cap V(H_b)) \cup \{a,b\}$.
Let $x,y$ be the neighbors of $w$ such that $wx, wy \in E(H_b)$.
We may assume that $x \notin \{a,b\}$ because $w$ has two neighbors in $H_b$ and $b \notin V(H_b)$.  
Note also the minimality of $H_a$ implies that
 $x,y$ are not on the same path of $H_a$ that contains $w$.

Now consider $H_a$.
If $w$ and $x$ are internal vertices of paths of $H_a$ that
intersect at a common branch vertex, then contracting $wx$ preserves $H_b$ and
also the non-planarity of $H_a$ (see \autoref{WandXAlongCommon}).  
This contradicts Lemma~\ref{Proposition3} part~(ii).

Therefore $x$ must be a branch vertex of $H_a$ or it is an internal vertex on a path of $H_a$ that does not
intersect at a common branch vertex with $w$.  
\begin{figure}[H]
	\centering
	\includegraphics[page=53]{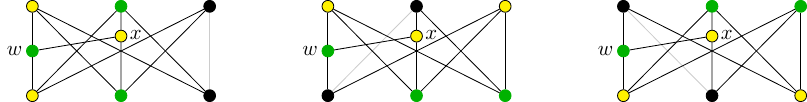}%
	\caption{For any branch vertex $t$ of $H_a$, there is a subdivision of $K_{3,3}$ in $H_a + wx$ in which $t$ is not a branch vertex.}
	\label{XtoNonBranchVertex}
\end{figure}
Assume that $x$ is not a branch vertex of $H_a$.  In this case,
\autoref{XtoNonBranchVertex} shows that for any branch vertex $t$ of $H_a$, there is a subdivision of $K_{3,3}$ in $H_a + wx$ in which $t$ is not a branch vertex.  Recall that $b$ is a branch vertex of $H_a$.
So, if $x$ is not a branch vertex of $H_a$, then there exists a subdivision of $K_{3,3}$ missing $a$ 
that does not have $b$ as a branch vertex, contradicting Corollary~\ref{BranchVertices}.

Consequently we may assume that $x$ is a branch vertex of $H_a$ different from the ones at the end of the path of $H_a$ containing $w$.  Without loss of generality $H_a$ appears as in see \autoref{XBranchVertex}, where the label
$v$ has been introduced on the branch vertex of $H_a$ along the path of $H_a$ containing $w$ that has opposite color
as $x$.

\begin{figure}[H]
	\centering
	\includegraphics[page=54]{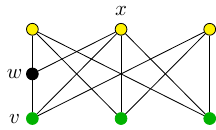}%
	\caption{$x$ is a branch vertex of $H_a$ (different from the ones at the end of the path of $H_a$ containing $w$).}
	\label{XBranchVertex}
\end{figure}

Now consider the subdivision of $K_{3,3}$, call it $J_v$, in $H_a + wx$ that remains after deleting the edges
along the $vx$-path (see \autoref{NewSubdivision}).

\begin{figure}[H]
	\centering
	\includegraphics[page=55]{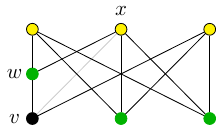}%
	\caption{The graph $J_v$ that is a subdivision of $K_{3,3}$ in $H_a + wx$.}
	\label{NewSubdivision}
\end{figure}
Note that $J_v$ covers all of the branch vertices of $H_a$ and does not use any vertices
outside of $H_a$.  By choice of $H_a$, there cannot be fewer edges
in $J_v$ than in $H_a$.  Consequently $J_v$ and $H_a$ must have the same number of edges.  In particular
the $vx$-path in $H_a$ is just the edge joining $v$ and $x$.
Because $H_a$ minimizes $|W|$ and the branch vertices of $J_v$ and $H_a$ only differ at $v$ and $w$,
the vertex $v$ is a branch vertex of neither $H_b$ nor $J_v$.
In particular the edge $vx$ must be an edge of $H_a$ and $H_b$.

If $w$ and $v$ are both neighbors of $a$, a $4$-cycle formed by $a,w,x$, and $v$ appears in $H_b$
(see \autoref{WVNotBothNeighborsOfA}).
However, since $w$ and $v$ are not branch vertices of $H_b$, this $4$-cycle is an impossible configuration
in $H_b$.

\begin{figure}[H]
	\centering
	\includegraphics[page=56]{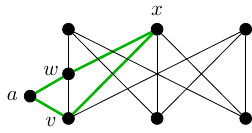}%
	\caption{If both $w$ and $v$ are neighbors of $a$, then an impossible $4$-cycle (green) appears in $H_b$.}
	\label{WVNotBothNeighborsOfA}
\end{figure}

By the choice of $w$ we conclude that $w \notin N_G(a)$ (if $w \in N_G(a)$ then replace $w$ with $v$).
It follows from earlier reasoning that, like $x$, the other neighbor of $w$ in $H_b$, namely $y$, must also be a branch vertex
of $H_a$.  

If $x$ and $y$ were branch vertices of $H_a$ with the same color, then there would be a subdivision of $K_{3,3}$
with fewer edges than $H_a$ that could have been chosen (as shown in \autoref{XYNOTSAMECOLORBRANCH}) contradicting
the choice of $H_a$.

\begin{figure}[H]
	\centering
	\includegraphics[page=57]{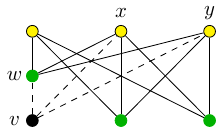}%
	\caption{Deleting edges along the three dotted paths and adding the two edges $wx,wy$ produces a subdivision of $K_{3,3}$ with fewer edges than $H_a$.}
	\label{XYNOTSAMECOLORBRANCH}
\end{figure}

So, without loss of generality, $H_a$ appears as shown in \autoref{WhasFull4}.  This
figure introduces labels for all of the branch vertices of $H_a$.
Observe that, like $v$, the vertex $u$ is not a branch vertex of $H_b$.  Also,
like the edge $vx$, the edge $uy$ must be an edge in $H_a$ and $H_b$.

\begin{figure}[H]
	\centering
	\includegraphics[page=58]{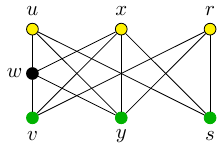}%
	\caption{In $H_a$, the vertex $w$ must have two neighbors, $x$ and $y$, that are branch vertices of opposite color in $H_a$ and not on the path of $H_a$ containing $w$.}
	\label{WhasFull4}
\end{figure}

If $u, v \in N_G(a)$, then the edges $au$,$uy$,$yw$,$wx$,$xv$,$va$ form a $6$-cycle in $H_b$, as shown
in \autoref{UorVNotNeighborOfA}.  Recall $a$ is a branch vertex of $H_b$ but 
$u,v$, and $w$ are not branch vertices of $H_b$.
If $x$ or $y$ are not branch vertices of $H_b$, then a cycle with at most two branch vertices of $H_b$
exists in $H_b$, an impossibility.
So, $a$, $x$ and $y$ are three branch vertices of $H_b$.  However the $6$-cycle 
induced by the edges $au$,$uy$,$yw$,$wx$,$xv$,$va$ from $H_b$ implies that $a$, $x$ and $y$ cannot be $2$-colored
as the branch vertices of a subdivision of $K_{3,3}$, contradicting that $H_b \cong TK_{3,3}$.

\begin{figure}[H]
	\centering
	\includegraphics[page=59]{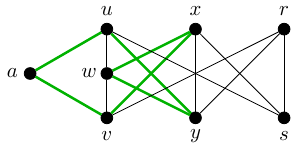}%
	\caption{If $u, v \in N_G(a)$, then an impossible $6$-cycle  (shown in green) emerges in $H_b$.}
	\label{UorVNotNeighborOfA}
\end{figure}

Consequently $u \notin N_G(a)$ or $v \notin N_G(a)$.  Without loss of generality, $v \notin N_G(a)$.
Now $d_G(v)=4$ so $v$ has another neighbor in $V(H_a) \cap V(H_b)$.  
Applying the same reasoning to $v$ as we have applied previously to $w$, we conclude
that the remaining unknown neighbor of $v$ must be a branch vertex of $H_a$.
Further applying this reasoning to $J_v$ (see \autoref{NewSubdivision}) reveals that
this neighbor of $v$ must be either $y$ or $s$.
However, if $y$ is a neighbor of $v$ in $H_b$, then the
four edges $vy,yw,wx,xv$ form an impossible $4$-cycle in $H_b$ with at most
two branch vertices (since $w$ and $v$ cannot be branch vertices of $H_b$).
Therefore, $s$ must be the final neighbor of $v$ in $H_b$ and $vs$ must be an edge of $H_b$ 
(see leftmost graph in \autoref{FinalFigure}).

\begin{figure}[ht]
	\centering
	\includegraphics[page=60]{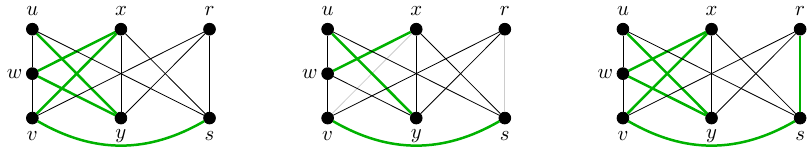}%
	\caption{$s \in N_G(v)$ (left) implies a subdivision of $K_{3,3}$ (middle) that implies a path in $H_b$ that covers the branch vertices of $H_a$ (path shown in green edges at right).}
	\label{FinalFigure}
\end{figure}

Because $s \in N_G(v)$ there is a subdivision of $K_{3,3}$, call it $J_r$,
that covers the branch vertices of $H_a$, contains only vertices from $H_a$,
but does not have $r$ as a branch vertex (see middle of \autoref{FinalFigure}).
Applying the same reasoning to $J_r$ as we applied before to $J_v$, 
we conclude that $r$ is a branch vertex of neither $J_r$ nor $H_b$.
Moreover the edge $rs$ is in $H_b$.  Consequently the path $uywxvsr$ has all its edges in $H_b$ and 
covers all of the branch vertices of $H_a$.  In particular the vertices of this path must be in $H_b$.
This is a contradiction because $b$ is in $V(H_a) - V(H_b)$ and,
by Corollary~\ref{BranchVertices}, it must also be a branch vertex of $H_a$.
This contradiction proves that $W=\varnothing$, as desired.
\end{proof}

\begin{figure}[htb]
	\centering
	\includegraphics[page=61,scale=0.675]{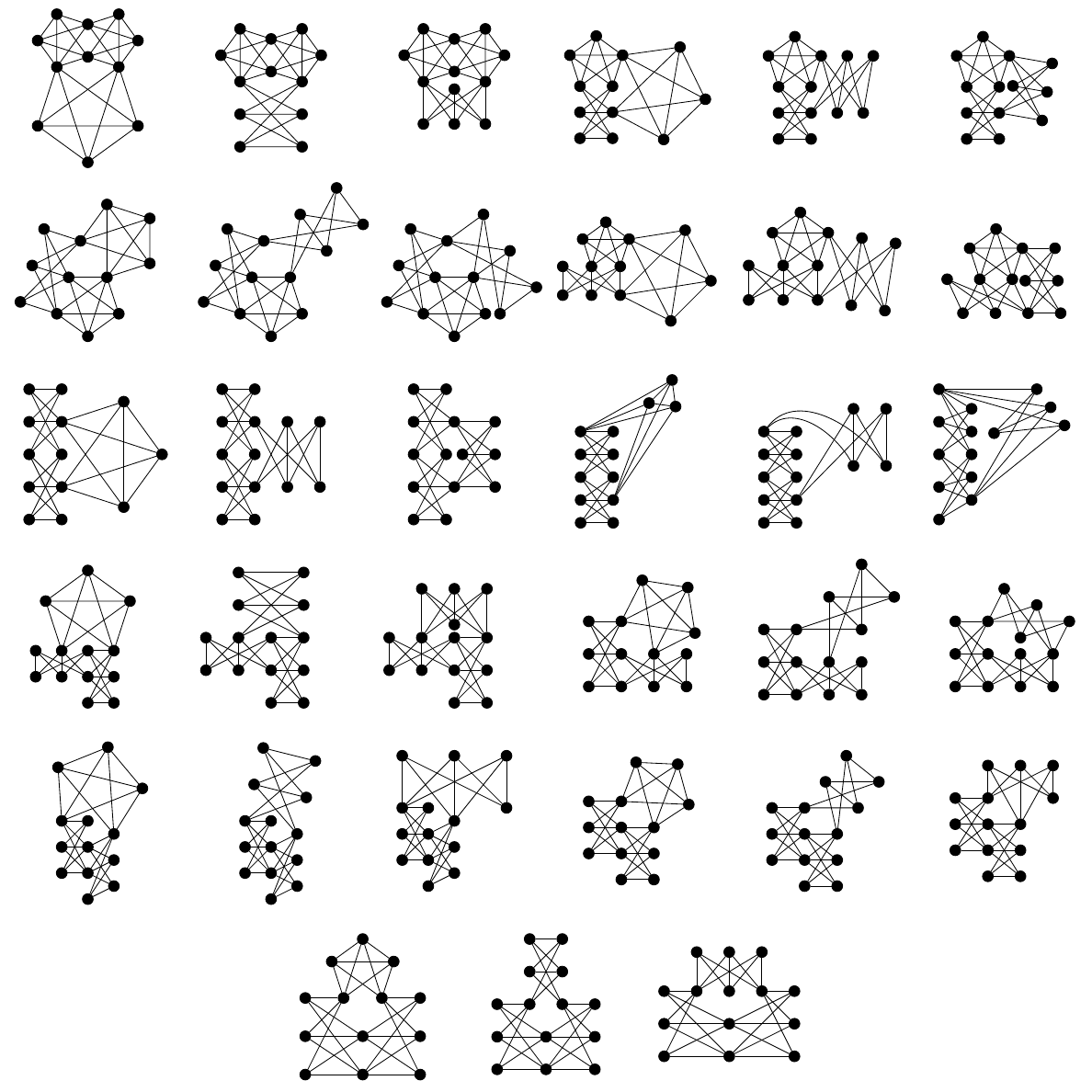}%
	\caption{The $33$ connectivity-$2$ apex obstructions that have a unique $2$-cut $\{a,b\}$, a planar heavy component $C$, and a $2$-cut separating $a$ from $b$ in $G[C \cup \{a,b\}]$.}
	\label{ThirtyThree}
\end{figure}

\subsection{The final computation}

In this subsection we outline how to show the list
of $72$ connectivity-$2$ apex obstructions satisfying Assumptions~\ref{assumptions} is complete.  While much of the case-analysis can be reduced by hand,
ultimately we confirmed the final list using computers.  We omit many details.  Much of the case analysis applies to small (order $\leq 10$) graphs and is routine, but it is sufficiently tedious that it precludes comprehensive presentation.

\begin{figure}[htbp]
	\centering
	\includegraphics[page=62,scale=0.675]{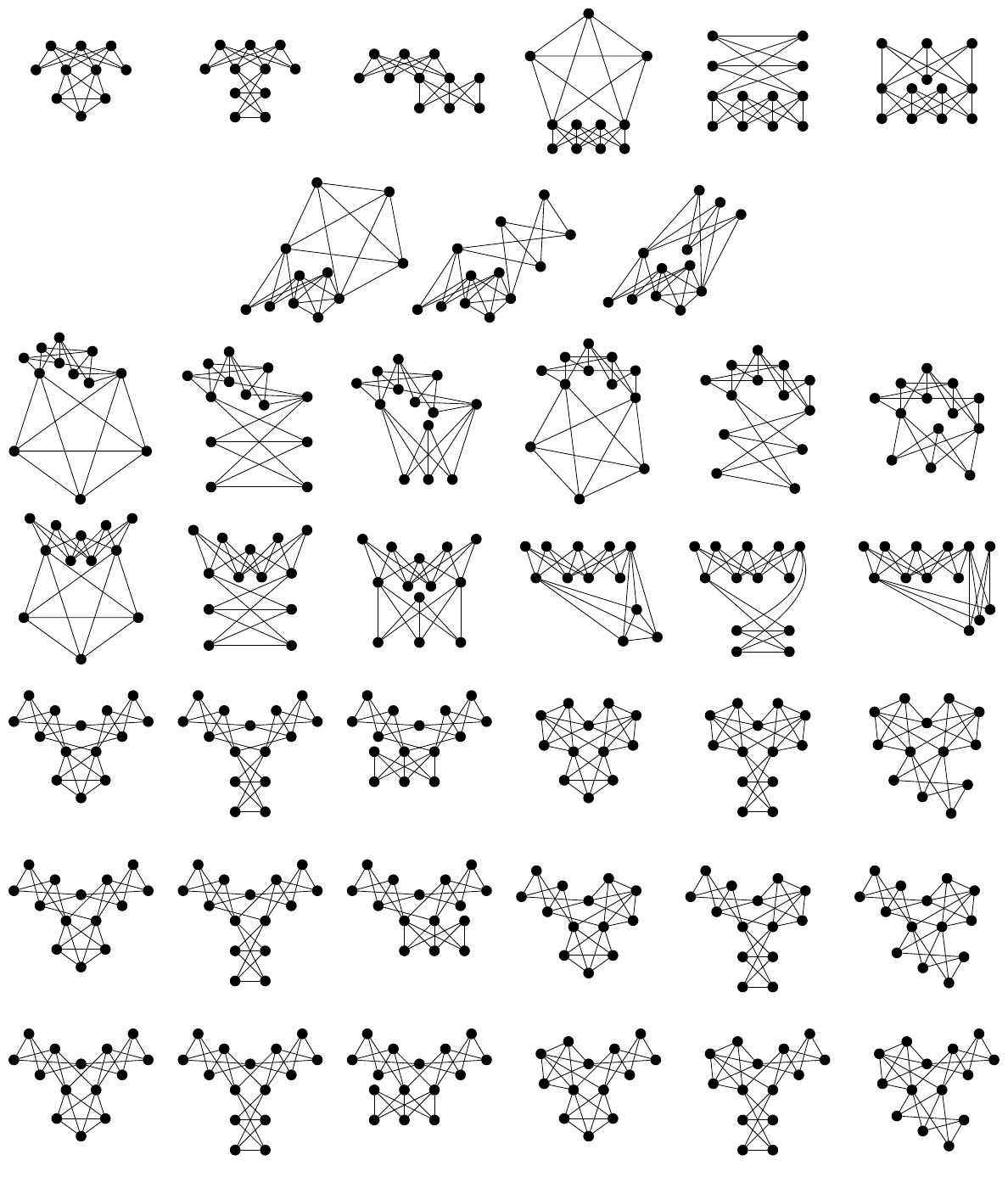}%
	\caption{The $39$ connectivity-$2$ apex obstructions that have a unique $2$-cut $\{a,b\}$, a planar heavy component $C$, and no $2$-cut separating $a$ from $b$ in $G[C \cup \{a,b\}]$.}
	\label{ThirtyNine}
\end{figure}

\begin{theorem}
Suppose that $G \in {\cal F}$, $\kappa(G) = 2$, $S=\{a,b\}$ is the unique $2$-cut of $G$,
$C$ is the heavy component of $G-S$ and $C$ is non-planar.
\begin{itemize}
\item If $G[C \cup \{a,b\}]$ has a $2$-cut separating $a$ from $b$, 
 then $G$ is isomorphic to a graph in \autoref{ThirtyThree}.
\item If $G[C \cup \{a,b\}]$ has no $2$-cut separating $a$ from $b$, 
 then $G$ is isomorphic to a graph in \autoref{ThirtyNine}.
\end{itemize}
\end{theorem}
\begin{proof}
By Theorem~\ref{WeakWTheorem}, there are Kuratowski witnesses $H_a$ and $H_b$ in $G$ avoiding $a$ and $b$ respectively,
such that 
any vertex in $C \cup \{a,b\}$ is a branch vertex of $H_a$ or a branch vertex of $H_b$.
Consequently, every vertex of $G$ is a branch vertex of $H_a$, $H_b$, or the Kuratowski witness in $L^+$.

Because $a$ and $b$ are branch vertices of $H_b$ and $H_a$ respectively, it follows that $|C| \leq 10$.
Indeed, the only way that $|C|=10$ is if $H_a$ and $H_b$ are subdivisions of $K_{3,3}$ with disjoint branch sets.
This case can be shown never to occur by examining possible subdivided edges of $H_a$ and $H_b$ that must involve
branch vertices from the other Kuratowski witness.  We omit the details.

Clearly $|C|\geq 4$ because $H_a$, for  one, has at least $5$ branch vertices.
If $|C|=4$, then it is easy to show $K_6$ is a minor of $G$, an impossibility.
Thus it suffices to consider cases in which $5 \leq |C| \leq 9$, $C$ is connected and planar.
There are $87,816$ non-isomorphic, connected planar graphs with order between $5$ and $9$ inclusive.
A computer check of all of these graphs (together with adding a light component) reveals the $72$ obstructions indicated.
\end{proof}

\section{A double apex graph interpretation}
\label{sec:DoubleApex}

In this section we discuss another interpretation of our characterization of connectivity-$2$ apex obstructions
because it may be of independent interest.

\begin{figure}[H]
	\centering
	\includegraphics[page=63,scale=1.0]{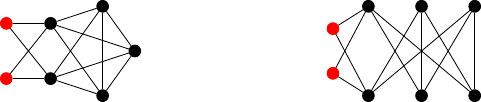}%
	\caption{The two rooted double apex graphs $K_5^*$ and $K_{3,3}^*$.}
	\label{TriviallyDoublyApex}
\end{figure}

A graph $H$ with vertices $a$ and $b$ is {\em double apex} (with respect to $a$ and $b$) if
$H-a$ and $H-b$ are non-planar, but $H-a-b$ is planar.  Vertices $a$ and $b$ are the {\em roots} of $H$.
For example, consider the rooted graphs obtained from $K_5$ and $K_{3,3}$ by replacing an edge with two
subdivided edges and making the degree-two vertices the roots; denote these
two graphs $K_5^*$ and $K_{3,3}^*$, respectively.  These graphs (shown in Figure~\ref{TriviallyDoublyApex} with the roots colored red)
are double apex.

Characterizing minor-minimal double apex graphs is a special type of `intertwine' problem in which
the goal is to determine minimal graphs that contain two Kuratowski subgraphs, each containing one root but avoiding the other.
Note that $K_5^*$ and $K_{3,3}^*$ are minor-minimal double apex graphs.
Recall that the Petersen family of graphs consists of the seven graphs: $K_6, K_{1,3,3},Y^-,K_{4,4}-e,M,P_7$, $P$.
Significant minor-minimal double apex graphs can be obtained from one of the seven Petersen family graphs
by removing an edge and making its endpoints the roots.  Note that not all such edge deletions
produce a minor-minimal double apex graph ($K_{1,3,3},Y^-,K_{4,4}-e,$ and $M$ have problematic edges).

One could restate 
Corollary~\ref{BranchVertices} in the language of double apex graphs this way:

\begin{theorem} \label{DoubleApexFormulation} Suppose that $H$ is a minor-minimal double apex graph with roots $a$ and $b$.
If $H + \{ab\} \not\cong Y^-,M,P_7,P$ and $H \not\cong K_5^*,K_{3,3}^*$, then $a$ (resp. $b$) is a branch vertex
of every Kuratowski subgraph in $H-b$ (resp. $H-a$). 
\end{theorem}

Using Theorem \ref{DoubleApexFormulation} one can prove 
(see the proof of Theorem \ref{WeakWTheorem}) that in any minor-minimal double apex graph $H$
satisfying $H + \{ab\} \not\cong Y^-,M,P_7,P$ and $H \not\cong K_5^*,K_{3,3}^*$, two Kuratowski subgraphs exist, one avoiding $a$ and the other avoiding $b$,
whose branch vertices cover the entire vertex set.  In this way all minor-minimal double apex graphs can
be enumerated.   Indeed it follows from our characterization of connectivity-$2$ apex obstructions, that
there are $57$ non-isomorphic (as rooted graphs!) minor-minimal double apex graphs.    
Three are disconnected graphs; these are rooted versions of $2K_5$, $2K_{3,3}$ and $K_5 + K_{3,3}$ in which
a root appears in each component.  Twelve can be obtained from Petersen family graphs by removing an edge.  For example,
there are four non-isomorphic rooted minor-minimal double apex graphs that can be obtained from $M$ by removing an edge.
The remaining $42$ non-isomorphic rooted minor-minimal double apex graphs (including $K_5^*$ and $K_{3,3}^*$) 
can be found by inspecting the augmented heavy components of the connectivity-$2$ apex obstructions: see \autoref{AppendixB} and consider
the $42$ `cards' shown there in which all three $2$-sums (with $K_5$, $K_{3,3}$, or $K_{3,3}+e$) produce an apex obstruction.

Assuming Theorem \ref{DoubleApexFormulation}, it is easy to reason that if $H$ is a minor-minimal double apex graph with roots $a$ and $b$ and
$H + \{ab\}$ is not isomorphic to a Petersen family graph, then a $2$-sum at the roots $a$ and $b$ of $H$ with $K_5$ or $K_{3,3}$ 
will produce a connectivity-$2$ apex obstruction.  Indeed Lemma~\ref{lightComponent} states that, under the right circumstances, a connectivity-$2$ apex obstruction must
be generated in this way.  So it is tempting to believe that knowing the $57$ 
non-isomorphic rooted minor-minimal double apex graphs suffices to generate all 
connectivity-$2$ apex obstructions via such $2$-sums with $K_5$ and $K_{3,3}$.
But these $2$-sums may generate the same graph more than one way.  More significantly,
because of non-planar heavy components or vertices that are not branch of any Kuratowski subgraph, 
there are $23$ apex obstructions that can not be generated as a $2$-sum in this way:
nineteen of the $21$ graphs in \autoref{NonPlanarC} and 
four graphs at the bottom of \autoref{Intersecting2Cuts}.
The two left-most graphs in the middle rows of Figure \ref{NonPlanarC} can be generated
from a $2$-sum applied to $K_5^*$ and $K_{3,3}^*$, respectively.
Incidentally it is worth mentioning that, among the $23$ apex obstructions that can not be generated as a $2$-sum but excepting the seven graphs that contain vertices that are not branch of any Kuratowski subgraph,
the remaining $16$ graphs require only the branch vertices of two Kuratowski subgraphs to cover the entire vertex set.


\subsection*{Acknowledgments}

We thank Csaba Biro and Tim Pervenecki for stimulating discussions and the referees
for their thoughtful comments and suggestions.  We are also very grateful to the editors for alerting
us to references \cites{HugoAyala,MR3733966,MikePierce}.
\clearpage




\clearpage
\noindent
\appendix
\section{graph6 codes for connectivity-2 apex obstructions}

Some general information on graph6 codes may be found here: 

\texorpdfstring{\protect{See \url{http://users.cecs.anu.edu.au/~bdm/data/formats.txt}}}{}
\label{graph6Appendix}

\subsection*{\autoref{NonPlanarC}}
\begin{verbatim}
I`KxuJBpw
J]oo?CJ@w^_
IJaK[\x\_
J]oo?SJ@wN_
Js?GZ`oBw^?
J?C^F?]FV@_
JWC]E?]FPF_
JWCXEC]FUD_
J@KEMJCNHu?
J`MAIIBNHu?
J`MAIIBNPt?
Ks?GOKH@b`mE
K??WooEwV@[K
Ko?WooEWR@ON
K?Ku?CKKUBWM
K]??WWKKME?]
K]?M@_IA_J_m
K]?M@_I@gQ_l
JwCWFFaFo|?
K]?GWB@KO^@y
Lr??WYOBGEK@?N
\end{verbatim}

\subsection*{\autoref{Disjoint2Cuts}}
\begin{verbatim}
IwC^F?^Fo
Jr?G[`_Bw^?
Kr?GOMOAWLKB
\end{verbatim}

\subsection*{\autoref{Intersecting2Cuts}}
\begin{verbatim}
KwCW?CB~FFb~
L]?GW?@?][ENB~
LF?GW?@?^[[MB~
M?ope???G@}?A^@n_
M?B_oo??G@~??~wM_
M??CZ_??G@~_bM[M_
No@_??B?ooB?f??NkWW
N?@_op__C?O?@N?vN?W
N?KsA@?OC?O??~BbEMG
NCO`@@?OC?O?ENDr@}?
LJ?GW?@?^{]Mb{
MBW?GK??G@{weM`}?
N?B@`b????_B}?BN@Z_
Oo@_??B?ooB???f?_FrEC
\end{verbatim}

\subsection*{\autoref{MoreIntersecting2Cuts}}
\begin{verbatim}
K`K?GN?N~pW|
L]?GW?@?XbxqB}
Ls?G?CBBBf`}^D
M?ope???N_?FA\@l_
M??Wv???G@_v}AwN?
LWCW?CBo@Fz`F{
MEo`?K??G@{G@nE]?
M?CV?W??G@`f{Pw]?
N??@`_KBE?W??NBp]?O
N??BB?[_C??Bw@FEbw_
M?KuE???G@}??~B]?
M??^?o??G@_n}@w]?
N???@_oBe?W?{?Bb`{_
LWCW?CBGEww]F{
M?NE?O??G@}G@}K]?
M??^?G??G@bN}Ow]?
No???oE@_oO?F`WRKE_
N??uE?GA?O?_{@?~EF_
L?CW?CBwFw[[F{
M?BE@o??G@~?MM@}?
M??F?w??G@bf~?w]?
N????oE@_o[?F`wQ[E?
N??EE?C@?GF?}??~FF?
\end{verbatim}

\subsection*{\autoref{ThirtyThree}}
\begin{verbatim}
JwC?G{]}^N?
K]?G?FKKo^`}
Ks?G?CNBrxm]
KWC?GKwuENB}
L?r@_?@?xbFao]
L??^??@@Wr^Aw]
K?CW?FbwvwB}
L?BE??@}@rFK@z
L??F??@FWz^_w]
KoCW?DbWsxb}
LEo`??@w?N_}EZ
L?CV??@BWZ]Bw]
LQQ@}???G@cnE]
M]??OGCA?F`_KdoX?
MF??OGCA?F`_wdwX?
LQQBKo??G@cnEm
M]??OGCA?K`KKdoY?
MF??OGCA?K`KwdwY?
LBa??CBFBWK]_}
M]??GGGA@`_]p@Aq_
M?w??KOC?B_uxa{B?
LW???CBNEwB{o{
M]?GO???@b_}A[p__
MF??W????F`mw[zA?
LSP??CBN@wO^O}
M???WZ?K?F@a{B{B?
MCaAA?_G?F_]VBNB?
Lo???CB^BwB{_}
M???@_oo?Bforara?
MS`A????@r_}@{]@_
L????CB~FwP{[w
M]??????E[EMB{B{?
Ms???????^`}]K\o?
\end{verbatim}

\subsection*{\autoref{ThirtyNine}}
\begin{verbatim}
I@NEE?~No
Jr??WWKczF?
Js??WWK[zf?
J@K?ENEnb{?
K]????NBuUEw
Ks????NBruMw
J?CXFFav`~?
K?r??CrKpyW]
KF???CNBvY[]
K??G[``{C|Lw
L]????N?oWxaKq
Ls????NAOKnH\c
K@?GXBPw?}xw
L]???OF?O[eqqK
LF???OF?O[{qyK
K???GNw}C}Lw
L]?????rHf@{Bw
Ls?????Bw^NK\g
K@C?GNgy?nxw
LBW?CA?@wNBswE
L?`o?A?AwVM[{E
MQQ@Go??G@wWH]Em?
N]??OGCA?C_KBBKdWL?
Ns??OGCA?C_KBB[dML?
K`?G]?o{]^F{
L]?G?CK?pxw]B{
Ls?G?CK?o^ne[{
MIa?X`L???_BKf_v?
NBW?GI?_?@_W@LEq[@_
N?w?GGOC?@_W@Lwq]@_
LoD_?CBECwk]F{
M?NE@_??G@`NLao]?
M?CV?W??G@aNzAw]?
MSP@x_K???_B_^O^?
N]??OOC@?E?E?{MBW`_
Ns??OOC@?E?E?{]BM`_
LWC?GN?EEoc}F{
M@JE?o??G@bFHqo]?
M??^??@E?BBLxEw]?
\end{verbatim}

\section{Augmented heavy components}
\label{AppendixB}
\begin{figure}[htbp]
	\hfill
	\includegraphics[page=1,scale=0.65]{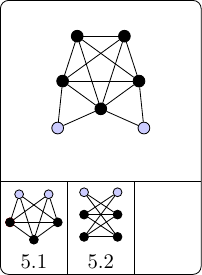}\hfill%
	\includegraphics[page=2,scale=0.65]{Augmented.pdf}\hfill%
	\includegraphics[page=3,scale=0.65]{Augmented.pdf}\hfill%
	\includegraphics[page=4,scale=0.65]{Augmented.pdf}\hfill%
	\includegraphics[page=5,scale=0.65]{Augmented.pdf}\hfill%
	\includegraphics[page=6,scale=0.65]{Augmented.pdf}\hfill\vspace{\baselineskip}%
	
	\includegraphics[page=7,scale=0.65]{Augmented.pdf}\hfill%
	\includegraphics[page=8,scale=0.65]{Augmented.pdf}\hfill%
	\includegraphics[page=9,scale=0.65]{Augmented.pdf}\hfill%
	\includegraphics[page=10,scale=0.65]{Augmented.pdf}\hfill%
	\includegraphics[page=11,scale=0.65]{Augmented.pdf}\hfill%
	\includegraphics[page=12,scale=0.65]{Augmented.pdf}\hfill%
	\includegraphics[page=13,scale=0.65]{Augmented.pdf}%
	\caption{Augmented heavy components appearing in \autoref{NonPlanarC} or \autoref{Disjoint2Cuts}}
\end{figure}
\begin{figure}[htbp]
	\hfill\hfill
	\includegraphics[page=14,scale=0.65]{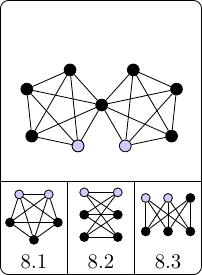}\hfill%
	\includegraphics[page=15,scale=0.65]{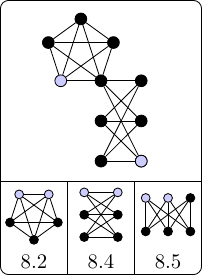}\hfill%
	\includegraphics[page=16,scale=0.65]{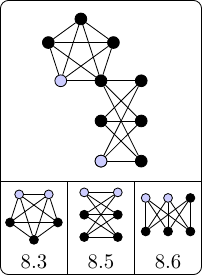}\hfill%
	\includegraphics[page=17,scale=0.65]{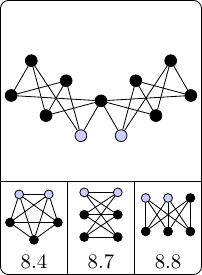}\hfill%
	\includegraphics[page=18,scale=0.65]{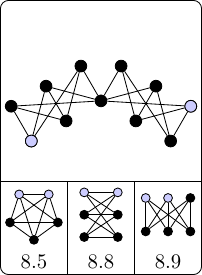}\hfill\hfill\vspace{\baselineskip}%
	
	\hfill\hfill\hfill
	\includegraphics[page=19,scale=0.65]{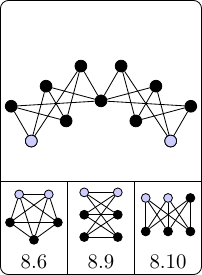}\hfill%
	\includegraphics[page=20,scale=0.65]{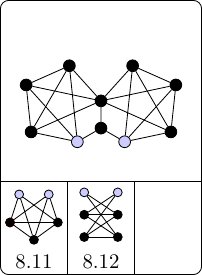}\hfill%
	\includegraphics[page=21,scale=0.65]{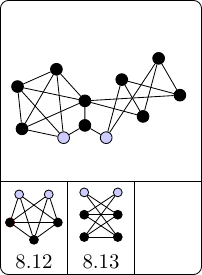}\hfill%
	\includegraphics[page=22,scale=0.65]{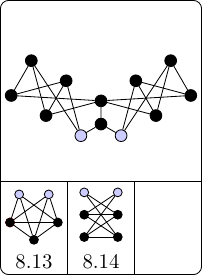}\hfill\hfill\hfill\vphantom{.}%
	
	\caption{Augmented heavy components appearing in \autoref{Intersecting2Cuts}}
\end{figure}
\begin{figure}[htbp]
	\hfill\hfill
	\includegraphics[page=23,scale=0.65]{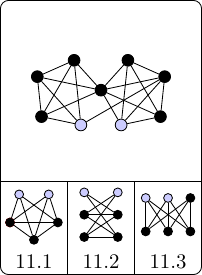}\hfill%
	\includegraphics[page=24,scale=0.65]{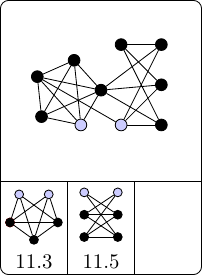}\hfill%
	\includegraphics[page=25,scale=0.65]{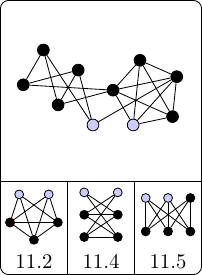}\hfill%
	\includegraphics[page=26,scale=0.65]{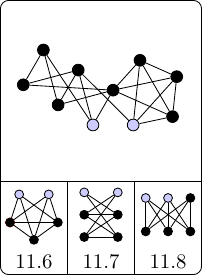}\hfill%
	\includegraphics[page=27,scale=0.65]{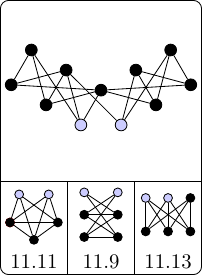}\hfill\hfill\vspace{\baselineskip}%
	
	\hfill\hfill
	\includegraphics[page=28,scale=0.65]{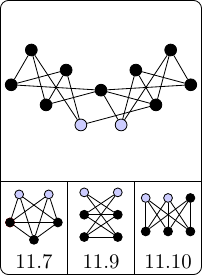}\hfill%
	\includegraphics[page=29,scale=0.65]{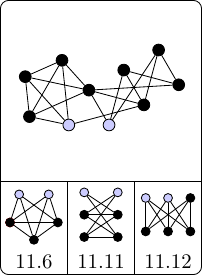}\hfill%
	\includegraphics[page=30,scale=0.65]{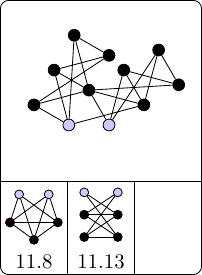}\hfill%
	\includegraphics[page=31,scale=0.65]{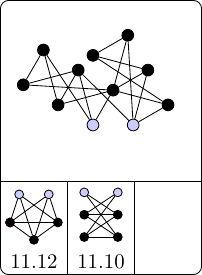}\hfill%
	\includegraphics[page=32,scale=0.65]{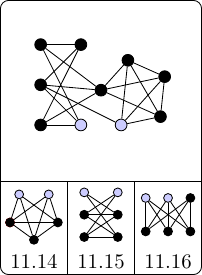}\hfill\hfill\vspace{\baselineskip}%
	
	\hfill\hfill
	\includegraphics[page=33,scale=0.65]{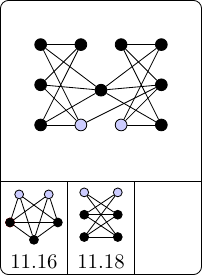}\hfill%
	\includegraphics[page=34,scale=0.65]{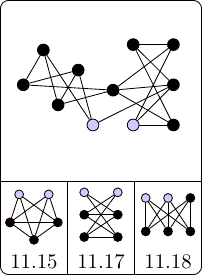}\hfill%
	\includegraphics[page=35,scale=0.65]{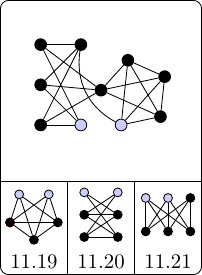}\hfill%
	\includegraphics[page=36,scale=0.65]{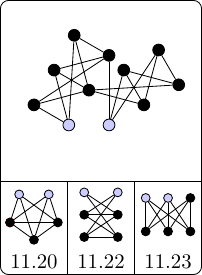}\hfill%
	\includegraphics[page=37,scale=0.65]{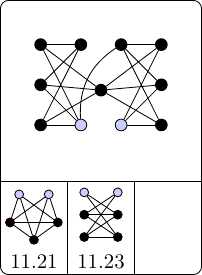}\hfill\hfill\vphantom{.}%
	
	\caption{Augmented heavy components appearing in \autoref{MoreIntersecting2Cuts}}
\end{figure}
\begin{figure}[htbp]
	\hfill
	\includegraphics[page=38,scale=0.65]{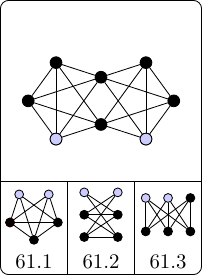}\hfill%
	\includegraphics[page=39,scale=0.65]{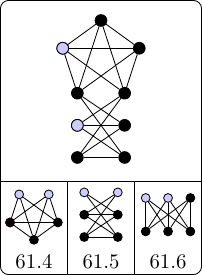}\hfill%
	\includegraphics[page=40,scale=0.65]{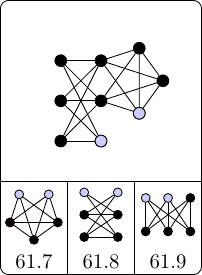}\hfill%
	\includegraphics[page=41,scale=0.65]{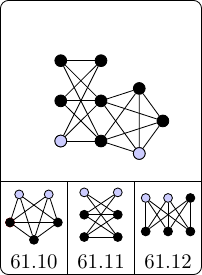}\hfill%
	\includegraphics[page=42,scale=0.65]{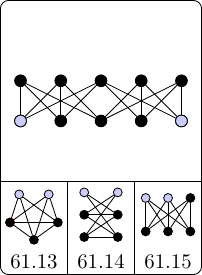}\hfill%
	\includegraphics[page=43,scale=0.65]{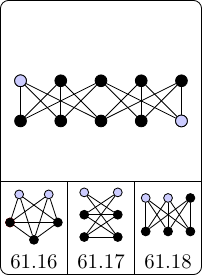}\vspace{\baselineskip}%

	\hfill\hfill
	\includegraphics[page=44,scale=0.65]{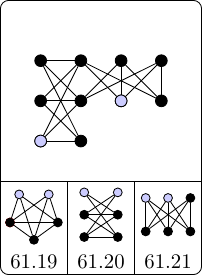}\hfill%
	\includegraphics[page=45,scale=0.65]{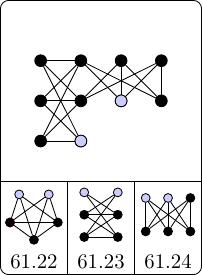}\hfill%
	\includegraphics[page=46,scale=0.65]{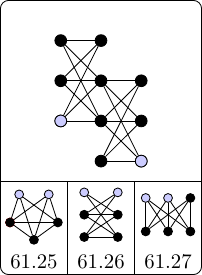}\hfill%
	\includegraphics[page=47,scale=0.65]{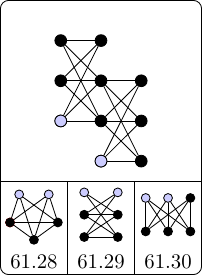}\hfill%
	\includegraphics[page=48,scale=0.65]{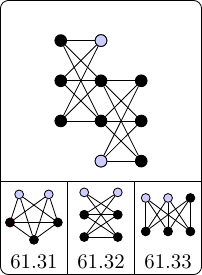}\hfill\hfill\vphantom{.}%
	
	\caption{Augmented heavy components appearing in \autoref{ThirtyThree}}
\end{figure}
\begin{figure}[htbp]
	\hfill
	\includegraphics[page=49,scale=0.65]{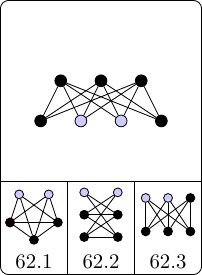}\hfill%
	\includegraphics[page=50,scale=0.65]{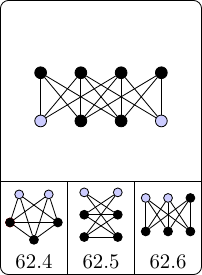}\hfill%
	\includegraphics[page=51,scale=0.65]{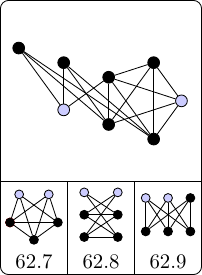}\hfill%
	\includegraphics[page=52,scale=0.65]{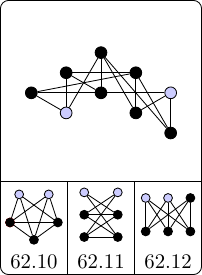}\hfill%
	\includegraphics[page=53,scale=0.65]{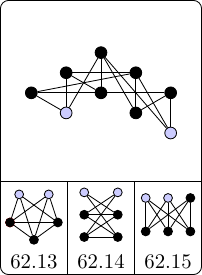}\hfill%
	\includegraphics[page=54,scale=0.65]{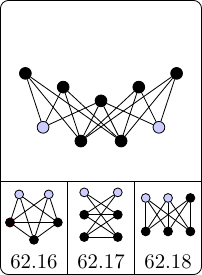}\hfill\vspace{\baselineskip}%
	
	\includegraphics[page=55,scale=0.65]{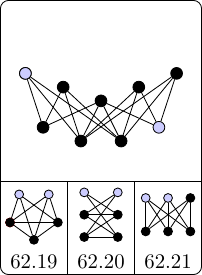}\hfill%
	\includegraphics[page=56,scale=0.65]{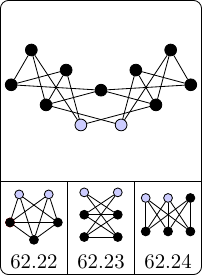}\hfill%
	\includegraphics[page=57,scale=0.65]{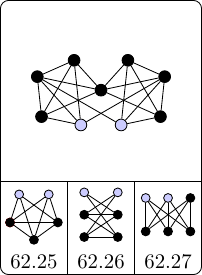}\hfill%
	\includegraphics[page=58,scale=0.65]{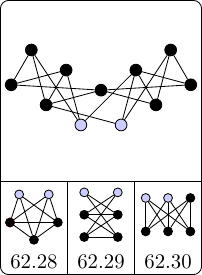}\hfill%
	\includegraphics[page=59,scale=0.65]{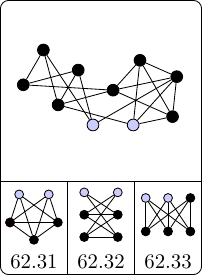}\hfill%
	\includegraphics[page=60,scale=0.65]{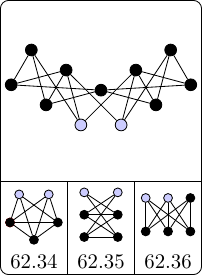}\hfill%
	\includegraphics[page=61,scale=0.65]{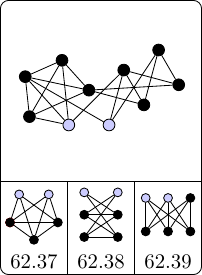}%
	
	\caption{Augmented heavy components appearing in \autoref{ThirtyNine}}
\end{figure}
\end{document}